\documentclass[11pt]{amsart}
\usepackage{amsmath, amssymb, amscd, mathrsfs, slashed, enumerate, url}

\usepackage[all,cmtip]{xy}
\usepackage{multicol}
\usepackage{indentfirst}
\usepackage{latexsym}
\usepackage{bm}
\usepackage{graphicx}
\usepackage{subfigure}
\usepackage{float}
\usepackage{verbatim}
\usepackage{comment}
\usepackage[top=1in, bottom=1in, left=1.25in, right=1.25in]{geometry}
\usepackage{epsfig,dsfont,amsthm,amsfonts,amsbsy}
\usepackage[colorlinks]{hyperref}

\newcommand{\gpp}{\mathfrak{g}_P}
\newcommand{\gppc}{\mathfrak{g}_P^{\mathbb{C}}}
\newcommand{\slc}{SL(2;\mathbb{C})}
\newcommand{\dA}{d_{\mathbb{A}}}
\newcommand{\LGF}{\mathcal{L}^{gf}}
\newcommand{\DKW}{\mathcal{D}_{(A,\Phi)}}
\newcommand{\DKWX}{\mathcal{D}_{(A^{\sharp},\Phi^{\sharp})}}
\newcommand{\GP}{\mathcal{G}_P}
\newcommand{\DKWA}{\mathcal{D}_{(A,\Phi),\alpha}}

\newcommand{\dAd}{d_{\mathbb{A}}^{\star}}
\newcommand{\HOC}{H^1_{\mathbb{A}}}
\newcommand{\HZC}{H^0_{\mathbb{A}}}
\newcommand{\CP}{\mathcal{C}_P}
\newcommand{\Hess}{Q_{\mathbb{A}}}
\newcommand{\EHess}{\widehat{Q}_{\mathbb{A}}}
\newcommand{\EHessr}{\widehat{Q}_{\mathbb{A}_{\rho}}}
\newcommand{\LKW}{\mathcal{L}_{(A,\Phi)}}
\newcommand{\Eng}{\mathcal{E}^{an}(t)}

\newcommand{\cs}{CS^{\mathbb{C}}}
\newcommand{\ro}{(A_{\rho},\Phi_{\rho})}
\newcommand{\Aro}{A_{\rho}}
\newcommand{\Pro}{\Phi_{\rho}}
\newcommand{\QAB}{\widetilde{\{a_t,b_t\}}}
\newcommand{\LTY}{L^2(Y_t)}

\newcommand{\lan}{\langle}
\newcommand{\ran}{\rangle}
\newcommand{\XT}{X^{\sharp T}}
\newcommand{\PT}{P^{\sharp T}}
\newcommand{\ML}{\mathcal{L}}
\newcommand{\MD}{\mathcal{D}}
\newcommand{\LamP}{y^{\lambda+\frac{1}{p}}}
\newcommand{\LamPP}{y^{\lambda+\frac{1}{p}}}
\newcommand{\LamMM}{y^{\lambda+\frac{1}{p}-1}}
\newcommand{\LaPP}{y^{\lambda+\frac{1}{p}+1}}
\newcommand{\LamM}{y^{\lambda+\frac{1}{p}-1}}
\newcommand{\FMS}{\mathcal{M}^{fr}_{p,\lambda}}

\newcommand{\LS}{\mathcal{L}^{\sharp T}}
\newcommand{\QS}{\mathcal{Q}^{\sharp T}}
\newcommand{\RF}{\mathbb{R}^4_+}
\newcommand{\MQ}{\mathcal{Q}}

\newcommand{\DZ}{d^{0}_{(A_0,\Phi_0)}}
\newcommand{\DZZ}{d^{0}_{(A,\Phi)}}
\newcommand{\DZSZ}{d^{0,\star}_{(A,\Phi)}}
\newcommand{\DZS}{d^{0,\star}_{(A_0,\Phi_0)}}
\newcommand{\DZZS}{d^{0,\star}_{(A,\Phi)}}
\newcommand{\QR}{Q_{\rho}}
\newcommand{\QRH}{\hat{Q}_{\rho}}
\newcommand{\SIS}{\mathcal{\sigma}}
\newcommand{\PS}{\mathcal{P}^{\sharp T}}
\newcommand{\RS}{\mathcal{\pi}^{\sharp T}}
\newcommand{\MK}{\mathcal{K}}
\newcommand{\MS}{\mathcal{S}}
\newcommand{\MM}{\mathcal{M}}
\newcommand{\MH}{\mathcal{H}}
\newcommand{\MHN}{\mathcal{H}_{(A,\phi,\phi_y)}}
\newcommand{\MSI}{\mathcal{M}_i}
\newcommand{\MSIS}{\mathcal{M}^{\star}_{P_i}}
\newcommand{\MSPS}{\mathcal{M}^{\star}_P}

\newcommand{\MSZ}{\mathcal{M}_X}
\newcommand{\CS}{\mathcal{C}_P}
\newcommand{\ox}{\vec{x}}
\newcommand{\Mk}{\kappa}
\newcommand{\FCS}{\mathcal{C}^{fr}_{p,\lambda}}
\newcommand{\Lone}{embeddinglemma}
\newcommand{\Ltwo}{Sobolevembedding}
\newcommand{\HRA}{\hookrightarrow}
\newcommand{\FBS}{\mathcal{B}_{p,\lambda}^{fr}}

\newcommand{\Ker}{\mathrm{Ker}}
\newcommand{\Vol}{\mathrm{Vol}}
\newcommand{\Coker}{\mathrm{Coker}}
\newcommand{\FGG}{\mathcal{G}^{fr}}
\newcommand{\FGGG}{\FGG_{p,\lambda}}
\newcommand{\BC}{\mathbb{C}}
\newcommand{\BA}{\mathbb{A}}
\newcommand{\GAS}{\Gamma_{(A,\Phi)}}

\newcommand{\CSM}{\mathrm{CS}}
\newcommand{\ECS}{\mathrm{ECS}}
\newcommand{\La}{\lambda}
\newcommand{\MEH}{\mathcal{EH}}
\newcommand{\MEHR}{\mathcal{EH}_\rho}
\newcommand{\MHR}{\mathcal{H}_\rho}
\newtheorem{theorem}{Theorem}[section]

\newtheorem{corollary}[theorem]{Corollary}

\newtheorem{definition}[theorem]{Definition}
\newtheorem{example}[theorem]{Example}

\newtheorem{lemma}[theorem]{Lemma}

\newtheorem{proposition}[theorem]{Proposition}
\newtheorem*{remark}{Remark}

\renewcommand{\Im}{\mathrm{Im}}
\newcommand{\Ind}{\mathrm{Ind}}
\newcommand{\MC}{\mathcal{C}}
\newcommand{\Tr}{\mathrm{Tr}}

\usepackage[utf8]{inputenc}
\usepackage[english]{babel}
\usepackage{fancyhdr}

\begin{document}

\title[A Gluing Theorem for the Kapustin-Witten Equations with a Nahm Pole]{A Gluing Theorem for the Kapustin-Witten Equations with a Nahm Pole}

\author[Siqi He]{Siqi He}
\address{Department of Mathematics, California Institute of Technology\\Pasadena, CA, 91106 }
\email{she@caltech.edu}

\begin{abstract}
In the present paper, we establish a gluing construction for the Nahm pole solutions to the Kapustin-Witten equations over manifolds with boundaries and cylindrical ends. Given two Nahm pole solutions with some convergence assumptions on the cylindrical ends, we prove that there exists an obstruction class for gluing the two solutions together along the cylindrical end. In addition, we establish a local Kuranishi model for this gluing picture. As an application, we show that over any compact four-manifold with $S^3$ or $T^3$ boundary, there exists a Nahm pole solution to the obstruction perturbed Kapustin-Witten equations. This is also the case for a four-manifold with hyperbolic boundary under some topological assumptions.
\end{abstract}

\maketitle

\begin{section}{Introduction}
In \cite{witten2011fivebranes}, Witten proposed a gauge theory approach to the Jones polynomial and Khovanov homology. Witten predicted that the coefficients of Jones polynomial should count certain solutions to the Kapustin-Witten equations over $\mathbb{R}^3\times (0,+\infty)$ with singular boundary conditions on $\mathbb{R}^3\times \{0\}$. See \cite{gaiotto2012knot} for a physics approach of this program.

Given a smooth 4-manifold $X$ with boundary, let $P$ denote a principal $SU(2)$ bundle over $X$ and let $\gpp$ be the adjoint bundle. Let $A$ be a connection over $P$ and $\Phi$ be a $\gpp$ valued one-form. The Kapustin-Witten equations are:
\begin{equation}
\begin{split}
    F_A-\Phi\wedge\Phi+\star d_A\Phi=0,\\
    d^{\star}_A\Phi=0.
    \label{KW}
\end{split}
\end{equation}

When the knot is empty, the singular boundary condition is called the Nahm pole boundary condition and in \cite{mazzeo2013nahm}, Mazzeo and Witten proved that there exists a unique Nahm pole solution to (\ref{KW}) which corresponds to the Jones polynomial of the empty knot. For a general 4-manifold $X$ with 3-manifold boundary $Z$, we hope to find ways to count the number of solutions to the Kapustin-Witten equations with the Nahm pole boundary condition over the boundary. This might lead to the discovery of some new invariants.

Therefore, a basic question to ask is whether there exists a solution to (\ref{KW}) with the Nahm pole boundary condition over a general 4-manifold with boundary? In \cite{he2015rotationally}, the author constructed some explicit solutions to the Kapustin-Witten equations over $S^3\times (0,+\infty)$. Kronheimer \cite{KrP} constructed some explicit solutions to the Kapustin-Witten equations over $Y^3\times (0,+\infty)$, where $Y^3$ is any hyperbolic closed 3-manifold.

Following Taubes \cite{taubes1982self} \cite{taubes1984self}, in order to prove the existence of solutions, we hope to establish a gluing theory for the Kapustin-Witten equations, such that the known Nahm pole model solutions can be glued to general 4-manifolds with boundary to obtain new Nahm pole solutions.

The main difference in gluing in the Nahm pole case compared to the gluing in the Yang-Mills case and the Seiberg-Witten case is that the Nahm pole boundary is not a classical non-degenerate elliptic boundary condition. However, it is a uniformly degenerate elliptic problem, as studied by R.Mazzeo \cite{rafe1991elliptic}. We mainly need the analytic tools developed in \cite{rafe1991elliptic} \cite{mazzeo2013nahm}.

For $i=1,2$, let $X_i$ be 4-manifolds with boundaries $Z_i$ and infinite cylindrical ends identified with $Y_i\times (0,+\infty)$. Let $(A_i,\Phi_i)$ be solutions to the Kapustin-Witten equations (\ref{KW}) over $X_i$ with Nahm pole boundary conditions over $Z_i$ and convergence to flat $\slc$ connections $(A_{\rho_i},\Phi_{\rho_i})$ over the cylindrical ends.

If $Y_1=Y_2$, $(A_{\rho_1},\Phi_{\rho_1})=(A_{\rho_2},\Phi_{\rho_2})$, we can define a new 4-manifold $X^{\sharp}$ and approximate solutions $(A^{\sharp},\Phi^{\sharp})$ by gluing together the cylindrical ends. See Figure \ref{fig:path3672}, where the shaded parts are glued together.

\begin{figure}[H]
    \centering
    \includegraphics[width=0.8\textwidth]{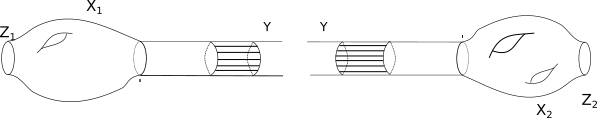}
    \caption{Gluing $X_1$, $X_2$ along the cylindrical ends}
    \label{fig:path3672}
\end{figure}

We prove the following theorem:

\begin{theorem}
Under the hypotheses above, if

(a)For some $p_0>2$, $\lim_{T\rightarrow +\infty}\|(A_i,\Phi_i)-(A_{\rho},\Phi_{\rho})\|_{L^{p_0}_1(Y_i\times\{T\})}=0$,

(b) $\rho$ is an acyclic flat $\slc$ connection,

then for $p\geq 2$ and $\lambda\in[1-\frac{1}{p},1)$ and sufficiently large $T$, we have:

(1) for some constant $\delta$, there exists a $y^{\lambda+\frac{1}{p}}H^{1,p}_0$ pair $(a,b)\in \Omega^1_{\XT}(\gpp)\times\Omega^1_{\XT}(\gpp)$ with $$\|(a,b)\|_{y^{\lambda+\frac{1}{p}-1}L^p_1}\leq Ce^{-\delta T},$$

(2) there exists an obstruction class $h\in H^2_{(A_1,\Phi_1)}(X_1)\times H^2_{(A_2,\Phi_2)}(X_2)$ such that $h=0$ if and only if $(A^{\sharp}+a,\Phi^{\sharp}+b)$ is a solution to the Kapustin-Witten equations (\ref{KW}).
\end{theorem}

In the statement of the theorem, $\rho$ acyclic means that $\rho$ is a regular point in the representation variety, and $H^2_{(A_i,\Phi_i)}$ means the cokernel of the linearization operator of the Kapustin-Witten equations over the point $(A_i,\Phi_i)$. Further, ${\LamP}H^{1,p}_0$ and ${\LamM}L^p_1$ are weighted norms which will be precisely introduced in Section 5. 

In addition, in Section 8, we also prove a gluing theorem when $\rho$ is reducible with a different weighted norm.

The statement and proof of Theorem 1.1 are analogous to the statement and proof of the gluing theorem for the ASD equation, due to C.Taubes \cite{taubes1982self}, \cite{taubes1984self}; cf. also \cite{donaldson2002floer}, \cite{donaldson1990geometry},
\cite{freed2012instantons}.

Moreover, for $p\in (2,4)$ and $\lambda\in [1-\frac{1}{p},1)$, denote by $\MM_i$ the moduli space of solutions to the Kapustin-Witten equations satisfying the assumption (a), (b) in Theorem 1.1 modulo the gauge action. We have the following Kuranishi model for the gluing picture.
\begin{theorem}
Let $(A_i,\Phi_i)$ be a connection pair over a manifold $X_i$ with a Nahm pole over $Z_i$. For sufficiently large $T$, there is a local Kuranishi model for an open set in the moduli space over $X^{\sharp}$: 

(1) There exists a neighborhood $N$ of $\{0\}\subset H^1_{(A_1,\Phi_1)}\times H^1_{(A_2,\Phi_2)}$ and a map $\Psi$ from $N$ to $H^2_{(A_1,\Phi_1)}\times H^2_{(A_2,\Phi_2)}$.

(2) There exists a map $\Theta$ which a homeomorphism from $\Psi^{-1}(0)$ to an open set $V\subset \mathcal{M}_{X^{\sharp}}.$
\end{theorem}

Here $H^k_{(A_i,\Phi_i)}$ is the $k$-th homology associated to the Kuranishi complex of $(A_i,\Phi_i)$ and $\mathcal{M}_{X^{\sharp}}$ is the moduli space of Nahm pole solutions to the Kapustin-Witten equations over $X^{\sharp}$. See also the Kuranishi model construction in Seiberg-Witten theory by T.Walpuski and A.Doan \cite{doan2017deformation}.

As for the model solutions, we don't know whether the obstruction class vanishes or not and right now we don't have any transversality results for the Kapustin-Witten equations. We just consider the obstruction class as a perturbation to the equation. See \cite{donaldson1983application} for the obstruction perturbation for ASD equations. We obtain the following theorem:

\begin{theorem}
Let $M$ be a smooth compact 4-manifold with boundary $Y$. Assume $Y$ is $S^3$, $T^3$ or any hyperbolic 3-manifold. When $Y$ is hyperbolic, we assume that the inclusion of $\pi_1(Y)$ into $\pi_1(M)$ is injective. For a real number $T_0$, we can glue $M$ to $Y\times (0,T_0]$ along $\partial M$ and $Y\times \{T_0\}$ to get a new manifold, which denote as $M_{T_0}$. For $T_0$ large enough, there exists a 
$SU(2)$ bundle $P$ and its adjoint bundle $\gpp$ over $M_{T_0}$ such that given any interior non-empty open neighborhood $U\subset M$, we have:

(1) There exist $h_1\in\Omega^2_{M_{T_0}}(\gpp),\;h_2\in\Omega^0_{M_{T_0}}(\gpp)$ supported on $U$,

(2) There exist a connection $A$ over $P$ and a $\gpp$-valued 1-form $\Phi$ such that $(A,\Phi)$ satisfies the Nahm pole boundary condition over $Y\times \{0\}\subset M_{T_0}$ and $(A,\Phi)$ is a solution to the following obstruction perturbed Kapustin-Witten equations over $M_{T_0}$: 
\begin{equation}
\begin{split}
    F_A-\Phi\wedge\Phi+\star d_A\Phi&=h_1,\\
    d_A^{\star}\Phi&=h_2.
\end{split}
\end{equation}
\end{theorem}

Here is the outline of the paper. In Section 2, we introduce some preliminaries on the Kapustin-Witten equations, including the Kuranishi complex and some examples of the Nahm pole solutions. In Section 3, we introduce a gauge fixing condition and the elliptic system associated to the equations. In Section 4, we study the gradient flow of the Kapustin-Witten equations, and the structure of the linearization operator over $Y\times \mathbb{R}$. In Section 5, we establish the Fredholm theory for the linearization operator over manifolds with boundaries and cylindrical ends. In Section 6, we build up a slicing theorem and Kuranishi model for the Nahm pole solutions. In Section 7, after assuming the solution over cylindrical ends is simple and $L^p_1$ converges to a flat $\slc$ connection over the cylindrical end for $p>2$, we prove that the solution will exponentially decay to the $\slc$ flat connection in the cylindrical ends. In Section 8, we describe the obstruction in the second homology group of the Kuranishi complex to the existence of solutions when gluing along the cylindrical ends. In Section 9, we build up a local Kuranishi model for the gluing picture. In Section 10, we apply the gluing theorem and get some existence results for the Nahm pole solutions to the perturbed equations. In Appendix 1, we introduce the $L^p$ version of Mazzeo's work for a uniformly degenerate elliptic operator. In Appendix 2, we introduce a proof of a Hardy type inequality for the weighted norm which is used to prove a slicing theorem.

\end{section}
\begin{section}{Preliminaries of the Kapustin-Witten Equations and the Nahm Pole Boundary condition}
In this section, we introduce some preliminaries on the Kapustin-Witten equations and the Nahm pole boundary condition.
\begin{subsection}{Kapustin-Witten Map}

Let $\hat{X}$ be a smooth compact connected four-manifold with two connected boundary components $Y$ and $Z$. Take $X$ to be the four-manifold obtained by gluing $\hat{X}$ and $Y\times [0,+\infty)$ along the common boundary $Y$ , that is $X:=\hat{X}\cup_Y (Y\times [0,+\infty)).$ For any positive real number $T$, we denote by $Y_T$ the slice $Y\times\{T\}\subset X$ and $X(T):=\hat{X}\cup_Y(Y\times(0,T))$. For simplicity, the metric we always consider on $X$ is cylindrical along a neighborhood of $Z$ and is the product metric over $Y\times [T,+\infty) $ for some $T$ big enough. This is illustrated in Figure \ref{fig:text4419}:

\begin{figure}[H]
    \centering
    \includegraphics[width=0.8\textwidth]{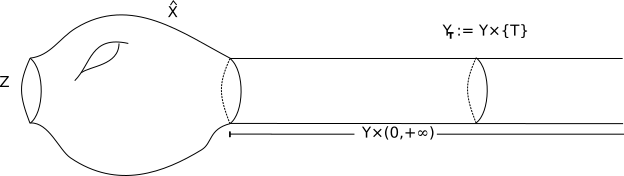}
    \caption{The shape of manifold we study}
    \label{fig:text4419}
\end{figure}

Now suppose $P$ is an $SU(2)$ bundle over $X$, $\gpp$  is the associated adjoint bundle and $\mathcal{A}_P$ is the set of all the $SU(2)$ connections on $P$. We define the configuration space as follows: $\mathcal{C}_P:=\mathcal{A}_P\times \Omega^1(\gpp),$ $\mathcal{C}^{'}_P:=\Omega^2(\gpp)\times \Omega^0(\gpp).$

The gauge-equivariant Kapustin-Witten map is the map $KW:\mathcal{C}_P\rightarrow \mathcal{C}^{'}_P$:
\begin{equation}
    \begin{split}
        KW(A,\Phi):&=
        \left( \begin{array}{l}
            F_A-\Phi\wedge\Phi+\star  d_A\Phi\\
            d_A^{\star}\Phi
        \end{array} 
        \right).
    \end{split}
\end{equation}

To be more explicit, denote by $\mathcal{G}_P$ the gauge group of $P$. Then, the action of $g\in \mathcal{G}_P$ on $(A,\Phi)\in \mathcal{A}_P\times \Omega^1(\gpp)$ is given by
$$
g(A,\Phi)=(A-(d_Ag)g^{-1},g\Phi g^{-1}).
$$

Under this action, the Kapustin-Witten map is gauge equivariant, i.e.,
$$KW(A-(d_Ag)g^{-1},g\Phi g^{-1})=gKW(A,\Phi)g^{-1}.$$



\end{subsection}
\begin{subsection}{Nahm Pole Boundary condition}
In \cite{witten2011fivebranes}, Witten proposed a gauge theoretic approach to Jones polynomial. A key objective of this program is to study the solutions to the Kapustin-Witten equations (\ref{KW}) satisfying the Nahm pole boundary condition.

To begin with, we introduce the Nahm pole boundary condition.

Given a 4-manifold $X$, with 3-dimensional boundary $Z$, a $SU(2)$ bundle $P$ over $X$ and the associated adjoint bundle $\gpp$, for integers $a=1,2,3$, take $\{\mathfrak{e}_a\}$ to be any unit orthogonal basis of $TZ$, the tangent bundle of $Z$, take $\{\mathfrak{e}^{\star}_a\}$ to be its dual and take 
$\{\mathfrak{t}_a\}$ to be section of 
the adjoint bundle $\gpp$ with the 
relation $[\mathfrak{t}_a,\mathfrak{t
}_b]=2\epsilon_{abc}\mathfrak{t}_c$. 
Identify a neighborhood of $Z$ with 
$Z\times (0,1)$, denote the boundary of $W$ by
$\partial W$ and identify it with 
$Z\times \{0\}$. We denote by $y$ as the 
coordinate on $(0,1)$.

\begin{definition}
A connection pair $(A,\Phi)\in\CP$ over X satisfies the Nahm pole boundary condition if there exist 
$\{\mathfrak{e}_a\}$, $\{\mathfrak{t}_a\}$ as above such that the 
expansion of $(A,\Phi)$ in $y\rightarrow 0$ of $Z\times (0,1)$ 
will be $A\sim A_0+yA_1+...$ and $\Phi\sim 
\frac{\sum_{a=1}^3e^{\star}_at_a}{y}+\Phi_0+y\Phi_1+....$. In addition, we call $(A,\Phi)\in\CP$ a Nahm pole solution if $(A,\Phi)$ is a solution to the Kapustin-Witten equations \rm{(\ref{KW})}.
\label{NahmPole}
\end{definition}

\begin{remark}
The definition of Nahm pole boundary condition depends on a choice of frame $\{e_a\}$, but locally up to gauge we can fixe a frame.
\end{remark}

In fact, a Nahm pole solution to the Kapustin-Witten equation will have more restrictions on the expansion, as pointed out in \cite{mazzeo2013nahm}.

\begin{proposition}{\rm{\cite{mazzeo2013nahm}} }For a Nahm pole solution $(A,\Phi)$ to the Kapustin-Witten equation, we have

(1) $\Phi_0=0.$

(2) Using $\sum e_a^{\star}t_a$ to identify $\gpp|_Y$ with $TY$, $A_0$ is the Levi-Civita connection of $Z$.
\end{proposition}

\begin{subsection}{Examples of Nahm Pole Solutions}

Here are some examples of solutions to the Kapustin-Witten equations satisfying the Nahm pole boundary condition.

\begin{example}
\rm{(Nahm \cite{nahm1980simple})Nahm pole solutions on $T^3\times \mathbb{R}^+$.
\label{NahmpoleT3}
Take the trivial $SU(2)$ bundle and denote $(A,\Phi)=(0,\frac{\sum \mathfrak{t}^idx^i }{y})$. Then $F_A=0$ and $\Phi\wedge\Phi=\frac{\sum[\mathfrak{t}^i,\mathfrak{t}^j]dx^i\wedge dx^j}{2y^2}$. In addition, $d_A\Phi=-\frac{\sum \mathfrak{t}^idy\wedge dx^i}{y^2}$. Therefore, $(A,\Phi)$ is a Nahm pole solution to the Kapustin-Witten equations $F_A-\Phi\wedge\Phi+\star d_A\Phi=0$.}
\end{example}

\begin{example}
\rm{
Nahm pole solutions on $S^3\times \mathbb{R}^+$.
\label{NahmpoleS3}
Equip $S^3$ with the round metric and take $\omega$ be 
Maurer–Cartan 1-form of $S^3$. Then, if $y$ is the 
coordinate of $\mathbb{R}^+$, denote 
\begin{equation}
(A,\Phi)=(\frac{6e^{2y}}{e^{4y}+4e^{2y}+1}\omega,\frac{6(e^{2y}+1)e^{2y}}{(e^{4y}+4e^{2y}+1)(e^{2y}-1)}\omega),
\label{S1}
\end{equation}

Theorem 6.2 in \cite{he2015rotationally} shows that $(A,\Phi)$ is a Nahm-Pole solution to the Kapustin-Witten equations. In addition, the solutions (\ref{S1}) will converge to the unique flat $\slc$ connection in the cylindrical end of $S^3\times \mathbb{R}^{+}$.}
\end{example}

\begin{example}
\rm{
(Kronheimer \cite{KrP}) Nahm pole solutions on $Y^3\times \mathbb{R}^+$ where $Y^3$ is any hyperbolic three manifold.
\label{NahmpoleY3}

Let $Y^3$ be a hyperbolic three manifold equipped with the hyperbolic metric $h$. Consider the associated $PSL(2;\mathbb{C})$ representation of $\pi_1(Y)$. By Culler's theorem \cite{Culler86lifting}, this lifts to $SL(2;\mathbb{C})$ and determines a flat $\slc$ connection $\nabla^{flat}$. Denote by $\nabla^{lc}$ the Levi-Civita connection and by $A^{lc}$ the connection form. Take $i\omega:=\nabla^{flat}-\nabla^{lc}$. Then locally, $\omega=\sum \mathfrak{t_i}e_i^{\star}$ where $\{e_i^{\star}\}$ is an orthogonal basis of $T^{\star}Y$ and $\{\mathfrak{t}_a\}$ are sections of 
the adjoint bundle $\gpp$ with the 
relation $[\mathfrak{t}_a,\mathfrak{t
}_b]=2\epsilon_{abc}\mathfrak{t}_c$. We also have $\star_Y\omega=F_{\nabla^{lc}}$. Therefore, by the Bianchi identity, we obtain $\nabla^{lc}(\star_Y\omega)=0$.

Combining $F_{flat}=0$ and the relation $\nabla^{flat}-\nabla^{lc}=i\omega$, we obtain $F_{\nabla^{lc}+i\omega}=0$. Hence $F_{lc}=\omega\wedge\omega$, $\nabla^{lc}\omega=0$.

Take $y$ to be the coordinate of $\mathbb{R}^+$ in $Y^3\times \mathbb{R}^+$, set $$f(y):=\frac{e^{2y}+1}{e^{2y}-1},$$ and take 
\begin{equation}
(A,\Phi)=(A^{lc},f(y)\omega).
\label{S3}
\end{equation} Clearly, $f(y)\rightarrow 1$ as $y\rightarrow +\infty$ and $f(y)\sim \frac{1}{y}$ as $y\rightarrow 0$.

Let us check that the solution satisfies the Kapustin-Witten equations over $Y^3\times (0,+\infty)$. We compute
\begin{equation*}
    \begin{split}
       &F_{A^{lc}}-\Phi\wedge\Phi
       =(1-f^2)F_{A^{lc}},\\
    &d_{A^{lc}}\Phi=d_{A^{lc}}\omega+f'(y)dy\wedge\omega
    =f'(y)dy\wedge\omega.
\end{split}    
\end{equation*}
and
\begin{equation*}
    \begin{split}
        &d_{A^{lc}}\star_{Y\times(0,+\infty)}(f(y)\omega)
        =d_{A^{lc}}(f(y)(\star_{Y}\omega)\wedge dy)
        =f(y) d_{A^{lc}}(\star\omega)\wedge dy
        =f(y) (d_{A^{lc}}F_{A^{lc}}\wedge dy)
        =0.
    \end{split}
\end{equation*}
Combining this with the previous equations and using the relation $1-f^2+f'=0$, we see that $KW(A,\Phi)=0.$

Since $f(y)\rightarrow 1$ as $y\rightarrow +\infty$, \rm{($\ref{S3}$)} converges to the $\slc$ flat connection $\rho$. }
\end{example}

\begin{example}
\rm{
Nahm Pole solutions on the unit disc $D^4$.}

This is an example from \cite{he2015rotationally} of the Nahm pole solution to the Kapustin-Witten equations (\ref{KW}) over a compact manifold with boundary. Identify the quaternions $\mathbb{H}$ with $\mathbb{R}^4$,  $x=x_1+x_2 I+x_3J+x_4K\in \mathbb{H}$ and let $D^4$ be the unit disc of $\mathbb{H}$. Now define:
$$(A,\Phi)=(Im(\frac{3}{|x|^4+4|x|^2+1}\bar{x}dx),\;Im(\frac{3(|x|^2+1)}{(|x|^4+4|x|^2+1)(|x|^2-1)}\bar{x}dx)).$$

It is shown in \cite{he2015rotationally} that this solution is a Nahm pole solution to the Kapustin-Witten equations over $D^4$.
\end{example}

\begin{example}{\rm{(S.Brown, H.Panagopoulos and M.Prasad \cite{BPP82})Two-sided Nahm pole solutions on $(-\frac{\pi}{2},\frac{\pi}{2})\times T^3$}

Consider the trivial $SU(2)$ bundle P over $(-\frac{\pi}{2},\frac{\pi}{2})\times T^3$ and let $\{\mathfrak{t}_a\}$ to be elements in $\gpp$ with the 
relation $[\mathfrak{t}_a,\mathfrak{t
}_b]=2\epsilon_{abc}\mathfrak{t}_c$ and $dx_i$ to be three orthogonal basis of cotangent bundle of $T^3$. Now define: \begin{equation}
(A,\Phi)=(0,\frac{1}{cos(y)}dx_1\mathfrak{t}_1+\frac{1}{cos(y)}dx_2\mathfrak{t}_3+\frac{sin(y)}{cos(y)}dx_3\mathfrak{t}_3)
\end{equation}

Then it is easy to check that $KW(A,\Phi)=0.$}

\end{example}

\begin{remark}
All these solutions over manifolds with cylindrical ends decay exponentially to flat $\slc$ connections. The $\Phi$ terms in these examples do not have a $dy$ component on the cylindrical ends.
\end{remark}
\end{subsection}

\end{subsection}
\begin{subsection}{Kuranishi Complex}
Now we will present the Kuranishi complex associated to the Kapustin-Witten equations (\ref{KW}). See also \cite{BMPhd} some similar computations for the Vafa-Witten equations

Given a connection pair $(A,\Phi)\in \CS$ satisfying  (\ref{KW}), the complex associated to $(A,\Phi)$ is:
\begin{equation}
0\rightarrow \Omega^0(\gpp)\xrightarrow{d^0_{(A,\Phi)}}\Omega^1(\gpp)\times\Omega^1(\gpp)\xrightarrow{d^1_{A,\Phi}}\Omega^2(\gpp)\times\Omega^0(\gpp)\rightarrow 0,
\label{complex}
\end{equation}
where $d^0_{(A,\Phi)}$ is the infinitesimal gauge transformation and $d^1_{(A,\Phi)}$ is the linearization of $KW$ at the pair $(A,\Phi)$.

To be more explicit, denote the Lie algebra of $\mathcal{G}_P$ by $\Omega^0(\gpp)$. Then, the corresponding infinitesimal action of $\xi\in\Omega^0(\gpp)$ will be:
\begin{equation}
\begin{split}
    &d^0_{(A,\Phi)}(\xi):\Omega^0(\gpp)\rightarrow \Omega^1(\gpp)\times \Omega^1(\gpp),\\
    &d^0_{(A,\Phi)}(\xi)=\left(\begin{array}{l}
    -d_A\xi    \\
    \big[\xi,\Phi\big]
    \end{array}\right).
\end{split}
\end{equation}
The linearization of the Kapustin-Witten equations at a point $(A,\Phi)$ is given by:
\begin{equation}
\begin{split}
    &d^1_{(A,\Phi)}:\Omega^1(\gpp)\times\Omega^1(\gpp)\rightarrow\Omega^2(\gpp)\times\Omega^0(\gpp),\\
    &d^1_{(A,\Phi)}
    \left(\begin{array}{l}
    a\\
    b
    \end{array}\right)=
    \left(\begin{array}{l}
    d_{A}a-[\Phi,b]+\star (d_Ab+[\Phi,a])\\
    -\star[a,\star\Phi]+d^{\star}_Ab
    \end{array}\right).
    \label{linearoperator}
\end{split}
\end{equation}

The following result about this Kuranishi complex is classical:

\begin{proposition}
The connection pair $(A,\Phi)\in\mathcal{C}_P$ satisfies $KW(A,\Phi)=0$ if and only if $d^1_{(A,\Phi)}\circ d^0_{(A,\Phi)}(\xi)=0$ $\forall \xi\in \Omega^0(\gpp)$.
\end{proposition}
\proof
The $\Omega^2(\gpp)$ component of the image of $d^1_{(A,\Phi)}\circ d^0_{(A,\Phi)}(\xi)$ equals:
\begin{equation*}
    \begin{split}
        &-d_Ad_A\xi+[\Phi,[\Phi,\xi]]+\star(d_A[\xi,\Phi]+[\Phi,-d_A\xi])\\
        =&-[F_A,\xi]+[\Phi\wedge\Phi,\xi]-\star[d_A\Phi,\xi]\\
        =&-[F_A-\Phi\wedge\Phi+\star d_A\Phi,\xi].
    \end{split}
\end{equation*}

While the $\Omega^0(\gpp)$ component is:
\begin{equation*}
    \begin{split}
        &d_A^{\star}[\xi,\Phi]-\star[-d_A\xi,\star\Phi]\\
        =&-\star d_A[\xi,\star\Phi]+\star[d_A\xi,\star\Phi]\\
        =&-[\xi,\star d_A\star\Phi]\\
        =&[\xi,d^{\star}_A\Phi].
    \end{split}
\end{equation*}
The statement follows immediately.
\qed

Therefore, $d^0_{(A,\Phi)}$ and $d^1_{(A,\Phi)}$ in (\ref{complex}) will form a complex. We can define the homology groups: $$H^0_{(A,\Phi)}=\Ker\; d^0_{(A,\Phi)},\;  H^1_{(A,\Phi)}=\Ker\; d^1_{(A,\Phi)}\slash \Im\; d^0_{(A,\Phi)},\text{ and }H^2_{(A,\Phi)}=\Coker\; d^1_{(A,\Phi)}.$$

We denote the isotropy group of connection pair $(A,\Phi)$ as $\Gamma_{(A,\Phi)}=\{u\in\mathcal{G}|u(A,\Phi)=(A,\Phi)\}.$ 
Recall that $H^0_{(A,\Phi)}$ is the Lie algebra of the stabilizer of $(A,\Phi)$ and $H^1_{(A,\Phi)}$ is the formal tangent space.

The formal dual Kuranishi complex with respect to the $L^2$ norm and Dirichlet boundary condition is:
\begin{equation}
0\rightarrow\Omega^2(\gpp)\times\Omega^0(\gpp)\xrightarrow{d^{1,\star}_{(A,\Phi)}}\Omega^1(\gpp)\times\Omega^1(\gpp)\xrightarrow{d^{0,\star}_{(A,\Phi)}}\Omega^0(\gpp)\rightarrow 0,
\label{dcomplex}
\end{equation}

Where
\begin{equation}
\begin{split}
&d^{1,\star}_{(A,\Phi)}:\Omega^2(\gpp)\times\Omega^0(\gpp)\rightarrow\Omega^1(\gpp)\times\Omega^1(\gpp),\\
&d^{1,\star}_{(A,\Phi)}\left(\begin{array}{l}
\alpha\\
\beta
\end{array}\right)=
    \left(\begin{array}{l}
    d_{A}^{\star}\alpha+\star[\Phi,\alpha]-[\Phi,\beta]\\
    -\star d_A\alpha+\star[\Phi,\star\alpha]+d_A\beta
    \end{array}\right).
\end{split}
\label{lstar}
\end{equation}

and 
\begin{equation}
    \begin{split}
        &d^{0,\star}_{(A,\Phi)}:\Omega^1(\gpp)\times\Omega^1(\gpp)\rightarrow\Omega^0(\gpp),\\
        &d^{0,\star}_{(A,\Phi)}\left(\begin{array}{l}
a\\
b
\end{array}\right)=
    \left(\begin{array}{l}
    -d_A^{\star}a+\star[\Phi,\star b]
    \end{array}\right).
    \end{split}
\end{equation}

The Kapustin-Witten map also has the following structure:
\begin{proposition}
\label{qe}
The map $KW$ has an exact quadratic expansion:
$$
KW(A+a,\Phi+b)=KW(A,\Phi)+d^1_{(A,\Phi)}(a,b)+\{(a,b),(a,b)\},
$$
where 
\begin{equation}
\begin{split}
    \{(a,b),(a,b)\}=\left(
    \begin{array}{l}
        a\wedge a-b\wedge b+\star[a,b] \\
        -\star[a,\star b]
    \end{array}
    \right).
\end{split}
\end{equation}
\label{structure}
\end{proposition}
\proof We have the following direct computation:
\begin{equation*}
    \begin{split}
        &F_{A+a}-(\Phi+b)\wedge(\Phi+b)+\star d_{A+a}(\Phi+b)\oplus d^{\star}_{A+a}(\Phi+b)\\
        =&F_A+d_Aa+a\wedge a-\Phi\wedge\Phi-[\Phi,b]-b\wedge b\\
        &+\star d_A\Phi+\star d_A b+\star[\Phi,a]+\star[a,b]\\
        &\oplus d_A^{\star}\Phi+d_A^{\star}b-\star[a,\star\Phi]-\star[a,\star b]\\
        =&KW(A,\Phi)+d^1_{(A,\Phi)}(a,b)+\{(a,b),(a,b)\}.
    \end{split}
\end{equation*}
\qed

\end{subsection}
\end{section}
\begin{section}{Gauge Fixing, Elliptic System and Inner Regularity}
Recall that we consider a smooth 4-manifold $X$ with boundary 3-manifold $Z$ and cylindrical ends identified with $Y\times (0,+\infty)$ and with an SU(2) bundle $P$ over X. In this section we will discuss the properties of solutions to the Kapustin-Witten equations (\ref{KW}) away from the boundary $Z$.

Suppose that $(A_0,\Phi_0)$ is a fixed reference connection pair in $\CP$ and write $(A,\Phi)=(A_0,\Phi_0)+(a,b)$. Our Sobolev norms used in this section are defined in the usual way: for example, for $a\in \Omega^1(\gpp)$, we write 
$$\|a\|_{L^p_k(X)}:=(\sum^k_{j=0}\|\nabla^j_{A_0}a\|^p_{L^p(X)})^{\frac{1}{p}},$$
and for a pair $(a,b)\in \Omega^1(\gpp)\oplus \Omega^1(\gpp)$, we write 
$$\|(a,b)\|_{L^p_k(X)}:=(\|a\|^p_{L^p_k(X)}+\|b\|^p_{L^p_k(X)})^{\frac{1}{p}},$$
for any $1\leq p\leq \infty$ and non-negative integer $k$.

\begin{subsection}{Gauge Fixing Condition}
For gauge-invariant equations, in order to use elliptic PDE theory, we need to define a suitable gauge fixing condition. 

Given a reference connection pair $(A_0,\Phi_0)$, in our situation, we can considered the traditional Coulomb gauge or another gauge differing by lower order terms which is associated by the operator $d^{0,\star}_{(A_0,\Phi_0)}$ in (\ref{dcomplex}). 

Denote \begin{equation}
\LGF_{(A_0,\Phi_0)}:=d^{0,\star}_{(A_0,\Phi_0)}
\label{KWGF1}
\end{equation} and we have the following definition:

\begin{definition} Let $(A,\Phi)\in \CP$ and denote $(a,b):=(A,\Phi)-(A_0,\Phi_0)$. We say $(A,\Phi)$ is in the Coulomb gauge relative to $(A_0,\Phi_0)$ if $d_{A_0}^{\star}a=0.$ In addition, we say $(A,\Phi)$ is in the Kapustin-Witten gauge relative to $(A_0,\Phi_0)$ if $(a,b)$ satisfies $\LGF_{(A_0,\Phi_0)}(a,b)=0$ or  $$d^{\star}_{A_0}a-\star[\Phi_0,\star b]=0.$$

\end{definition}

For the Coulomb gauge fixing, there are some known results in \cite{wehrheim2004uhlenbeck} for compact manifolds with boundary:
\begin{proposition}{\rm{\cite{wehrheim2004uhlenbeck} Theorem 8.1}}\label{gaugefixingboundary}
Suppose $U$ is a compact submanifold of $X$, $P$ is the $SU(2)$ bundle over $X$. Fixed a reference connection pair $A_0$, there exists a constant $C$ depending on $A_0$ such that if $(A,\Phi)\in\CP$ and for $p> 2$,  
$$\|A-A_0\|_{L^{p}_1(U)}\leq C$$
then there exists a gauge transformation $u\in\GP$ such that 
\begin{equation}
\begin{split}
    &d_{A_0}^{\star}(u(A)-A_0)=0,\\
    &\star(u(A)-A_0)|_{\partial U}=0.
\end{split}
\end{equation}
\end{proposition}
We also prove a simple result on the existence of the Kapustin-Witten gauge representatives, working over a closed base manifold.

\begin{proposition}
Let $M$ be a closed 3 or 4-dimensional manifold and let P be an SU(2) bundle over $M$, fix $(A_0,\Phi_0)\in \mathcal{C}_P:=\mathcal{A}_P\times \gpp$. There exists a constant $c(A_0,\Phi_0)$ such that if $(A,\Phi)\in \mathcal{C}_P$ and for some $p>2$,
$$\|(A-A_0,\Phi-\Phi_0)\|_{L^p_1(M)}\leq c(A_0,\Phi_0),$$
then there is a gauge transformation $u\in \mathcal{G}_P$ such that $u(A,\Phi)$ is in the Kapustin-Witten gauge relative to $(A_0,\Phi_0)$. 
\label{KWG}
\end{proposition}
\proof
Denote $a:=A-A_0$ and $b:=\Phi-\Phi_0$, so by definition, for $u\in \mathcal{G}_P$, the gauge group action on $(A,\Phi)$ will be:
$$u(A_0+a)-A_0=uau^{-1}-(d_{A_0}u) u^{-1},$$
$$u(\Phi_0+b)-\Phi_0=u\Phi_0 u^{-1}+ubu^{-1}-\Phi_0.$$
The equation to be solved for $u\in \mathcal{G}_P$, is
\begin{equation}
d^{\star}_{A_0}(uau^{-1}-(d_{A_0}u)u^{-1})-\star[\Phi_0,\star(u\Phi_0u^{-1}-\Phi_0)]-\star[\Phi_0,\star ubu^{-1}]=0.
\end{equation}

We write $u=exp(\chi)=e^{\chi}$ for a section $\chi\in\gpp$, and define 
$$G(\chi,a,b):=d^{\star}_{A_0}(e^{\chi}ae^{-\chi}-(d_{A_0}e^{\chi})e^{-\chi})-\star[\Phi_0,\star(e^{\chi}\Phi_0e^{-\chi}-\Phi_0)]-\star[\Phi_0,\star e^{\chi}be^{-\chi}].$$

To solve the equation $G(\chi,a,b)=0$, we  use the implicit function theorem.
We extend the domain of $G$ to sections $\chi \in L^p_2(M)$ and bundle valued 1-forms $a,b\in L^p_1(M)$. Since for $p>2$, $L^p_2(M)$ sections are continuous in 3-dimensions and 4-dimensions, we get $G(\chi,a,b)$ in $L^p(M)$. The derivative of $G$ at $\chi=0$, $a=0$, $b=0$ is 
\begin{equation*}
    \begin{split}
        &DG(\xi,\alpha,\beta)\\
        =&d_{A_0}^{\star}\alpha-d_{A_0}^{\star}d_{A_0}\xi-\star[\Phi_0,\star[\xi,\Phi_0]]-\star[\Phi_0,\star\beta]\\
        =&-d_{A_0}^{\star}d_{A_0}\xi-\star[\Phi_0,\star[\xi,\Phi_0]]+d^{\star}_{A_0}\alpha-\star[\Phi_0,\star\beta].
    \end{split}
\end{equation*}

Denote $H(\xi):=-d_{A_0}^{\star}d_{A_0}\xi-\star[\Phi_0,\star[\xi,\Phi_0]]$ and $I(\alpha,\beta):=d^{\star}_{A_0}\alpha-\star[\Phi_0,\star\beta],$ then $$DG(\xi,\alpha,\beta)=H(\xi)+I(\alpha,\beta).$$

Denote $\mathbb{A}_0=A_0+i\Phi_0$ and define  $$d_{\mathbb{A}_0}^{\star}(\alpha+i\beta):=d^{\star}_{A_0}\alpha-\star[\Phi_0,\star \beta]+i(d^{\star}_{A_0}\beta+\star[\Phi,\star \alpha]).$$ Obviously, $$I(\alpha,\beta)=Re( d_{\mathbb{A}_0}^{\star}(\alpha+i\beta)).$$

If we show that the operator $H$ is surjective to the image of $I$, the implicit function theorem will give a small solution $\chi$ to the equation $G(\chi,a,b)=0$. Thus we will study the cokernel of the operator $H$.

If $\eta\in \Coker\; H$, we have
$$\lan H(\xi),\eta\ran =0\text{ for all }\xi,$$
by taking $\xi=\eta$, we obtain 
\begin{equation*}
    \begin{split}
        \lan H(\eta),\eta\ran=-\|d_{A_0}\eta\|_{L^2(M)}-\|[\eta,\Phi_0]\|_{L^2(M)}.
    \end{split}
\end{equation*}
Therefore, any element $\eta$ in the cokernel of $H$ satisfies $\|d_{A_0}\eta\|=0$ and $\|[\eta,\Phi_0]\|=0$, which implies $$d_{\mathbb{A}_0}(\eta)=0.$$ 

If for some $\alpha_0,\beta_0$, $I(\alpha_0,\beta_0)$ is not in the image of $H$, we have $d_{\mathbb{A}_0}I(\alpha_0,\beta_0)=0$ and we know that:
\begin{equation*}
\begin{split}
    0=&\lan d_{\mathbb{A}_0}Re(d_{\mathbb{A}_0}^{\star}(\alpha_0+i\beta_0)),\alpha_0+i\beta_0\ran\\
    =&\lan Re(d_{\mathbb{A}_0}^{\star}(\alpha_0+i\beta_0)),d_{\mathbb{A}_0}^{\star}(\alpha_0+i\beta_0)\ran\\
    =&\lan Re(d_{\mathbb{A}_0}^{\star}(\alpha_0+i\beta_0)),Re(d_{\mathbb{A}_0}^{\star}(\alpha_0+i\beta_0)) \ran.
\end{split}
\end{equation*}
By taking the real part of the inner product, we get $I(\alpha_0,\beta_0)=0$. 

Therefore, $H$ is surjective to the image of $I$ and by implicit function theorem, we prove the result.
\qed
\end{subsection}

\begin{subsection}{Elliptic System}
Now, we will use the gauge fixing condition to get an elliptic system associated to the Kapustin-Witten equations. This is also considered similarly for the Vafa-Witten equations in \cite{BMPhd}. We denote $\mathcal{L}_{(A_0,\Phi_0)}:=d^1_{(A_0,\Phi_0)}.$ Here the $d^1_{(A_0,\Phi_0)}$ is the linearization of the Kapustin-Witten map in (\ref{complex}).

Given $(A_0,\Phi_0)\in \mathcal{C}_P$, by Proposition \ref{qe}, for $(a,b)\in \Omega^1(\gpp)\times \Omega^1(\gpp)$, the equation $KW(A_0+a,\Phi_0+b)=\psi_0$ is equivalent to
$$\mathcal{L}_{(A_0,\Phi_0)}(a,b)+\{(a,b),(a,b)\}=\psi_0-KW(A_0,\Phi_0).$$

To make the equation elliptic in the interior, it is natural to add the gauge fixing condition $\LGF_{(A_0,\Phi_0)}(a,b)=0$ or $d_{A_0}^{\star}a=0$.  

By adding the Kapustin-Witten gauge, we define the Kapustin-Witten operator $\mathcal{D}_{(A_0,\Phi_0)}$:
\begin{equation}
\begin{split}
    &\mathcal{D}_{(A_0,\Phi_0)}:\Omega^1(\gpp)\times \Omega^1(\gpp)\rightarrow\Omega^2(\gpp)\times\Omega^0(\gpp)\times\Omega^0(\gpp)\\
    &\mathcal{D}_{(A_0,\Phi_0)} :=\mathcal{L}_{(A_0,\Phi_0)}+\LGF_{(A_0,\Phi_0)}
\end{split}
\label{KWO}
\end{equation}

Denote $\psi=\psi_0-KW(A_0,\Phi_0)$, 
then the elliptic system can be rewritten as :
\begin{equation}
\mathcal{D}_{(A_0,\Phi_0)}(a,b)+\{(a,b),(a,b)\}=\psi.
\label{ES}
\end{equation}

Similarly, for the Coulomb gauge, we can denote $\hat{\mathcal{D}}_{(A_0,\Phi_0)}:=\mathcal{L}_{(A_0,\Phi_0)}+d_{A_0}^{\star}$, then we get another elliptic system:
\begin{equation}
\hat{\mathcal{D}}_{(A_0,\Phi_0)}(a,b)+\{(a,b),(a,b)\}=\psi.
\label{ES1}
\end{equation}

\end{subsection}

\begin{subsection}{Interior Regularity of the Elliptic System}
Local interior estimates for the elliptic system (\ref{ES}) are considered in \cite{feehan1998rm} in the context of $PU(2)$ monopoles.

\begin{theorem}{\rm{\cite{feehan1998rm}}}\label{interiorregularity}
Take a bounded open set $\Omega$ in the interior part of $X$ and let $P$ to be a principal $SU(2)$ bundle over X. Suppose that $P|_{\Omega}$ is trivial and $\Gamma$ is a smooth flat connection. Suppose that $(a,b)$ is an $L^2_1(\Omega)$ solution to the elliptic system (\ref{ES}) over $\Omega$ and take the back ground pair $(A_0,\Phi_0)=(\Gamma,0)$, where $\psi$ is in $L^2_{k}(\Omega)$ for $k\geq 1$ an integer. There exist an constant $\epsilon=\epsilon(\Omega)$ such that 
for any precompact open subset $\Omega'\subseteq \Omega$, if $\|(a,b)\|_{L^4(\Omega)}\leq \epsilon$
we have $(a,b)\in L^2_{k+1}(\Omega')$ and there is a universal polynomial $Q_k(x,y)$, with positive real coefficients, depending at most on $k,\Omega,\Omega'$ with $Q_k(0,0)=0$ and 

$$\|(a,b)\|_{L^2_{k+1}(\Omega')}\leq Q_k(\|\psi\|_{L^2_k(\Omega)},\|(a,b)\|_{L^2(\Omega)}).$$

In addition, if $(\psi,\tau)$ is in $C^{\infty}(\Omega)$ then $(a,b)$ is in $C^{\infty}(\Omega')$ and if $(\psi,\tau)=0$, then

$$\|(a,b)\|_{L^2_{k+1}(\Omega')} \leq C\|(a,b)\|_{L^2(\Omega)}.$$
\end{theorem}
\end{subsection}
By the previous theorem, we can get interior regularity for the solutions with Nahm pole boundary conditions:

\begin{corollary}
If $(A,\Phi)$ is a solution to the Kapustin-Witten equations (\ref{KW}) over $X$, for $\Omega\subset$ a bounded open set, if $\|(A,\Phi)\|_{L^p_1(\Omega)}$ is bounded, then for any proper open subset $\Omega'\subset\Omega$, $(A,\Phi)$ is smooth over $\Omega'$. 
\end{corollary}
\proof
Applying Proposition \ref{gaugefixingboundary} and Theorem \ref{interiorregularity}, the corollary comes out immediately.

\qed
\end{section}

\begin{section}{Gradient Flow}
In this section, we will discuss the gradient flow associated to Kapustin-Witten equations over a cylinder $X:=Y\times \mathbb{R}$. See Taubes \cite{taubes2013compactness} for a general computation of the topological twitsted equations. Denote the coordinate in $\mathbb{R}$ as $y$, then we use the product metric on $X$ and the volume form we specify is $\Vol_Y\wedge dy$, where $\Vol_Y$ is a volume form over $Y$. In this section, we denote by $\star_4$ the 4-dimensional Hodge star operator of $\Vol_Y\wedge dy$ and denote by $\star$ the 3 dimensional Hodge star operator with respect to $\Vol_Y$. 

\begin{subsection}{Generalized Gradient Flow Equations}
To begin, suppose $P$ is an $SU(2)$ bundle over $X$ and $A$ is a given connection on $P$. Using parallel transport along the slice of $Y$ into $Y\times \mathbb{R}$, we can consider $A$ as a map from $\mathbb{R}$ to connections on $P$. Similarly, the field $\Phi$ can be written as $\Phi=\phi+\phi_ydy$. Here $\phi$ is a map from $\mathbb{R}$ to $\Omega^1(\gpp)$ and $\phi_y$ is a map from $\mathbb{R}$ to the section of the adjoint bundle $\gpp$.

If $(A,\Phi)=(A,\phi+\phi_ydy)$ obeys the Kapustin-Witten equations ($\ref{KW}$), then we compute
\begin{equation}
\begin{split}
F_A-\Phi\wedge\Phi&=-\frac{d}{dy}Ady+F_A+[\phi_y,\phi]dy-\phi\wedge\phi,\\
\star_4 d_A\Phi&=\star d_A\phi_y-\star\frac{d}{dy}\phi+\star d_A\phi\wedge dy,\\
\star_4 d_A\star_4\Phi&=\frac{d}{dy}\phi_y+\star d_A\star\phi=0.
\end{split}
\end{equation}

Thus the Kapustin-Witten equations ($\ref{KW}$) are reduced to the following flow equations:

\begin{equation}
\begin{split}
    &\frac{d}{dy}A-\star d_A\phi-[\phi_y,\phi]=0,\\
&\frac{d}{dy}\phi-d_A\phi_y-\star(F_A-\phi\wedge\phi)=0,\\
&\frac{d}{dy}\phi_y- d_A^{\star}\phi=0.
\label{BigGradientFlow}
\end{split}
\end{equation}

These gradient flow equations are closely related to the complex Chern-Simons functional. Denote $\mathbb{A}:=A+i\phi$, then the complex Chern-Simons functional $\CSM^{\mathbb{C}}(\mathbb{A})$  on $3$-manifold $Y^3$ is: 
\begin{equation}
\CSM^{\mathbb{C}}(\mathbb{A}):=\int Tr(\mathbb{A}\wedge d\mathbb{A}+\frac{2}{3}\mathbb{A}\wedge\mathbb{A}\wedge\mathbb{A}).
\label{CCS}
\end{equation}

Now, we define the following functional for the flow equations ($\ref{BigGradientFlow}$)
\begin{definition}
Use the notation above, the extended Chern-Simons functional $\ECS$ is denoted as follows:
\begin{equation}
\ECS(A,\phi,\phi_y)=\frac{1}{2}\Im (\CSM^{\mathbb{C}}(\mathbb{A}))+\int _YTr(\phi\wedge\star d_A\phi_y),
\label{ECS}
\end{equation}
where the $\Im$ is taking the imaginary part of the complex Chern-Simons functional.
\end{definition}
 
Then we have the following proposition:

\begin{proposition}
Equation (\ref{BigGradientFlow}) is the gradient flow for the extended Chern-Simons functional.
\end{proposition}
\proof
Take $A=A_0+a$, $\phi=\phi_0+b$, $\phi_y=(\phi_y)_0+c$, then the linearization of $\frac{1}{2}\Im(\CSM^{\mathbb{C}}(\mathbb{A}))$ is 
$$\int b\wedge (F_{A_0}-\phi_0\wedge\phi_0)+\int a\wedge d_A\phi.$$ In addition, the linearization
of $\int_YTr(\phi\wedge\star d_A\phi_y)$ is $$\int_Y Tr(b\wedge\star d_{A_0}(\phi_y)_0)+\int_Y Tr(a\wedge\star[(\phi_y)_0,\phi_0])+\int_Y Tr(c\wedge\star d_{A_0}^\star\phi_0 ).$$ 

Therefore, the gradient of $\ECS$ at $(A_0,\phi_0,(\phi_y)_0)$ is
\begin{equation}
\nabla \ECS{(A_0,\phi_0,(\phi_y)_0)}= (-\star d_{A_0}\phi_0-[(\phi_y)_0,\phi_0],\;-d_{A_0}(\phi_y)_0-\star(F_{A_0}-\phi_0\wedge\phi_0),\;-d^{\star}_{A_0}\phi_0),
\label{GECS}
\end{equation}
where the minus sign is coming from the inner product we take for $s,s'\in\Omega^0(\gpp)$ is $-\Tr(ss')$.

The result follows immediately.
\qed

In addition, we can compute the Hessian operator $\MHN$ of $-\ECS$ at point $(A,\phi,\phi_y)$:

\begin{equation}
\begin{split}
        &\MHN:\Omega_Y^1(\gpp)\times\Omega_Y^1(\gpp)\times\Omega_Y^0(\gpp)\rightarrow \Omega_Y^1(\gpp)\times\Omega_Y^1(\gpp)\times \Omega_Y^0(\gpp),\\
    &\MHN\left(\begin{array}{l}
    a\\
    b\\
    c
    \end{array}
    \right)
    =\left(\begin{array}{l}
    \star d_A b+\star[\phi,a]+[\phi_y,b]-[\phi,c]\\
    \star d_Aa-\star[\phi,b]+d_A c-[\phi_y,a]\\
    d_A^{\star}b-\star[a,\star\phi]
    \end{array}
    \right).
    \label{HECS}
\end{split}
\end{equation}

Then we have the following expansion of the extended Chern-Simons functional $ECS$:

\begin{proposition}
For $(a,b,c)\in \Omega_Y^1(\gpp)\times\Omega_Y^1(\gpp)\times\Omega_Y^0(\gpp)$, we have the following expansions: 
(1) For the $\ECS$, we have
\begin{equation}
\begin{split}
    &\ECS(A+a,\phi+b,\phi_y+c)-\ECS(a,b,c)\\
    =&\int_Y(\lan (a,b,c),  \nabla\ECS{(A,\phi,\phi_y)}\ran-\frac{1}{2}\lan (a,b,c),\MHN(a,b,c)\ran-\{a,b,c\}^3),
\end{split}
\end{equation}
where $\nabla\ECS{(A,\phi,\phi_y)}$ is the gradient of $\ECS$ defined in (\ref{GECS}), $\MHN$ is the Hessian operator (\ref{HECS}) and $\{a,b,c\}^3$ are cubic terms of $a,b,c$.

To be explicit, if we denote $\mathbb{B}=a+ib$, then the cubic terms are 
\begin{equation}
\{a,b,c\}^3=\frac{1}{3}\Im(\mathbb{B}\wedge\mathbb{B}\wedge\mathbb{B})+b\wedge\star[a,c].
\end{equation}

(2) For $\nabla\ECS$, we have 
\begin{equation}
\begin{split}
&\nabla\ECS(A+a,\phi+b,\phi_y+c)-\nabla\ECS(A,\phi,\phi_y)\\
=&-\MHN(a,b,c)-\{a,b,c\}^2,
\end{split}
\end{equation}
where 
\begin{equation}
\begin{split}
\{a,b,c\}^2=\left(\begin{matrix}
    \star[a,b]+[c,b]\\
    [a,c]+\star(a\wedge a-\star b\wedge b)\\
    -\star[a,\star b]\\
    \end{matrix}
    \right).
\end{split}
\end{equation}
\label{balalala}
\end{proposition}
\proof
By a direct computation, we can verify these results.
\qed

As we have the gauge action, we can formally defined the extended Hession operator for $\ECS$:

The extended Hession operator $\MEH$ at the point $(A,\phi,\phi_y)$ is defined as:
\begin{equation}
\begin{split}
    &\MEH_{(A,\phi,\phi_y)}:\Omega_Y^1(\gpp)\times\Omega_Y^1(\gpp)\times\Omega_Y^0(\gpp)\times\Omega_Y^0(\gpp)\rightarrow \Omega_Y^1(\gpp)\times\Omega_Y^1(\gpp)\times\Omega_Y^0(\gpp)\times\Omega_Y^0(\gpp),\\
    &\MEH_{(A,\phi,\phi_y)}\left(\begin{array}{l}
    a_1\\
    b_1\\
    a_0\\
    b_0
    \end{array}
    \right)
    =\left(\begin{array}{l}
    \star d_A b_1+\star[\phi,a_1]+d_Aa_0-[\phi,b_0]+[\phi_y,b_1]\\
    \star d_Aa_1-\star[\phi,b_1]+d_A b_0+[\phi,a_0]-[\phi_y,a_1]\\
    d_A^{\star}a_1+\star[b_1,\star\phi]+[\phi_y,b_0]\\
    d_A^{\star}b_1-\star[a_1,\star\phi]+[\phi_y,a_0]
    \end{array}
    \right).
\end{split}
\end{equation}

We have the following proposition of these two Hessian operators:
\begin{proposition}
(1) $\MEH_{(A,\phi,\phi_y)}(a_1,b_1,0,b_0)=\MHN(a_1,b_1,b_0).$

(2) $\MEH$ and $\MH$ are self-adjoint operators.
\label{SASA}
\end{proposition}
\proof
By a direct computation, we can verify these results.
\qed

In some case of 4-manifold with boundary and cylinderical ends, we can have some simplification of the flow equations (\ref{BigGradientFlow}). Let $X$ to be a 4-manifold with boundary $Z$ and cylindrical end which identified with $Y\times(0,+\infty)$ and we denote $y$ to be the coordinate of $(0,+\infty)$.

\begin{definition}\label{SimpleS}
Let $(A,\Phi)$ to be a solution to the Kapustin-Witten equations (\ref{KW}) over $X$. Over the cylindrical end $Y\times (0,+\infty)$, let $\phi_y$ be the $dy$ component of $\phi$. The solution $(A,\Phi)$ is called \textbf{simple} if there exists $T_0$ such that the restriction of $\phi_y$ over $Y\times (T_0,+\infty)$ is zero.
\end{definition}
Here is an identity due to Taubes:
\begin{lemma}{\rm{(\cite{taubes2013compactness} Page 36)}}
If $(A,\phi,\phi_y)$ satisfies (\ref{BigGradientFlow}), we have the following identity:
$$\frac{1}{2}(-\frac{\partial^2}{\partial y^2}|\phi_y|^2+d^{\star}d|\phi_y|^2)+|\frac{\partial}{\partial y}\phi_y|+|d_A\phi_y|+2|[\phi_y,\phi]|^2=0.$$
\end{lemma}

With this identity, we have an immediate corollary by the maximum principle:
\begin{corollary}
Let $(A,\Phi=\phi+\phi_y dy)$ be a solution to the Kapustin-Witten equations or equivalently the flow equations (\ref{BigGradientFlow}) over $Y\times I$ where $I\subset \mathbb{R}$. We have the following:

(1) Over $Y\times \mathbb{R}$, if $sup_Y|\phi_y|$ has limit zero in the non-compact directions of $\mathbb{R}$, then $\phi_y=0$ over $Y\times \mathbb{R}$.

(2) Over $Y\times (0,+\infty)$, let $y$ be the coordinate of $(0,+\infty)$. If $(A,\Phi)$ satisfies the Nahm pole boundary condition over $Y\times \{0\}\subset Y\times (0,+\infty)$ and converges under $C^0$ norm to a flat $\slc$ connection when $y\rightarrow +\infty$, we have $\phi_y=0$.
\end{corollary}
\proof
(1) is an immediately corollary of the previous lemma and maximal principle.

For any $(A,\Phi=\phi+\phi_y dy)$ in the assumption of (2), by the definition of Nahm pole boundary condition and flat $\slc$ connection, we know $$\lim_{y\rightarrow 0}\sup|\phi_y|_{Y^3\times \{y\}}=0,\; \lim_{y\rightarrow +\infty}\sup|\phi_y|_{Y^3\times \{y\}}=0.$$ (2) also follows from an application of maximal principal.  
\qed

Therefore, we have the following simplification of the gradient flow equation:
\begin{corollary}
When $\phi_y=0$, (\ref{BigGradientFlow}) reduces to simple gradient flow equations:
\begin{equation}
    \begin{split}
        &\frac{d}{dy}A-\star d_A\phi=0,\\
&\frac{d}{dy}\phi-\star(F_A-\phi\wedge\phi)=0,\\
&d_A^{\star}\phi=0.
    \end{split}
    \label{SimpleGF}
\end{equation}
This will be the gradient flow to the functional $Im(CS^{\mathbb{C}}(\mathbb{A}))$ along with the stability condition $d_A^{\star}\phi=0$.
\end{corollary}

\begin{proposition}
If $(A,\phi)$ satisfies the first two equations of \rm{(\ref{SimpleGF})}:
\begin{equation*}
    \begin{split}
        &\frac{d}{dy}A-\star d_A\phi=0,\\
&\frac{d}{dy}\phi-\star(F_A-\phi\wedge\phi)=0,
    \end{split}
\end{equation*}
then $\frac{d}{dy}(d_A\star\phi)=0$.
\end{proposition}
\proof
We compute:
\begin{equation*}
    \begin{split}
        \frac{d}{dy}(d_A\star\phi)&=d\star \frac{d}{dy}\phi+[\frac{d}{dy}A,\star\phi]-[\star\frac{d}{dy}\phi,A]\\
        &=d_A(F_A-\phi\wedge\phi)+\star d_A\phi\wedge\star \phi-\star\phi\wedge\star d_A\phi\\
        &=0.
    \end{split}
\end{equation*}
\qed

\end{subsection}

\begin{subsection}{Acyclic Connection of the Characteristic Variety }

Now, consider the behavior of the complex Chern-Simons functional in a neighborhood of an $SL(2;\mathbb{C})$ flat connection. For a 3 manifold $Y^3$, let $\rho:\pi_1(Y^3)\rightarrow \slc$ be an $\slc$ representation of $Y^3$'s fundamental group and let $\mathbb{A}=A+i\phi$ be the flat $\slc$ connection associated to $\rho$. For simplicity, we only consider simple solutions (Definition \ref{SimpleS}) to the Kapustin-Witten equations.

A flat $SL(2;\mathbb{C})$ connection will satisfy the equations:
\begin{equation*}
    \begin{split}
        &F_A-\phi\wedge\phi=0,\\
        &d_A\phi=0.
    \end{split}
\end{equation*}

Denote $\gppc$ as the complexification of $\gpp$, then a flat $\slc$ connection will bring in a twisted de Rham complex:
\begin{equation}
    0\rightarrow\Omega^0(\gppc)\xrightarrow{d_{\mathbb{A}}} \Omega^1(\gppc)\xrightarrow{d_{\mathbb{A}}}\Omega^2(\gppc)\xrightarrow{d_{\mathbb{A}}}\Omega^3(\gppc)\rightarrow 0,
\end{equation}
with the homology groups:
$$
H^k_{\mathbb{A}}:=\frac{\Ker(\dA:\Omega^k(\gppc)\rightarrow\Omega^{k+1}(\gppc))}{\Im(\dA:\Omega^{k-1}(\gppc)\rightarrow\Omega^{k}(\gppc))}.
$$

In addition, we have the natural identification given by the real part and imaginary part of the bundle: $$\Omega^k(\gppc)\cong\Omega^k(\gpp)\oplus\Omega^k(\gpp).$$
To be explicit, given $a+ib\in \Omega^k(\gppc)$, we have the following maps:
\begin{equation*}
\begin{split}
    &\dA:\Omega^k(\gppc)\rightarrow \Omega^{k+1}(\gppc),\\
    &\dA(a+ib)=d_Aa-[\phi,b]+i(d_Ab+[\phi,a]),\\
    &\dAd:\Omega^{k+1}(\gppc)\rightarrow\Omega^k(\gppc),\\
    &\dAd(a+ib)=d_A^{\star}a-\star[\phi,\star b]+i(d_A^{\star}b+\star[\phi,\star a]).
\end{split}
\end{equation*}

\begin{remark}
Given $s,s'\in\Omega^0(\gppc)$, under the identification of $\Omega^0(\gppc)\cong \Omega^0(\gpp)\oplus\Omega^0(\gpp)$, there exist $s_1,s_2,s_1',s_2'\in\Omega^0(\gpp)$ such that $s=s_1+is_2$ and $s'=s_1'+s_2'$. The inner product we take is $\lan s,s'\ran =-\Tr(s\bar{s}')=\lan s_1s_1'\ran+\lan s_2s_2'\ran$ and $\dAd$ is the adjoint of $\dA$ with respect to this inner product. This explains the sign of $\dAd(a+ib)$.
\end{remark}

Now we will discuss the Hessian for the complex Chern-Simons functional.

For the first two equations of (\ref{SimpleGF}), we have:
\begin{equation}
    \begin{split}
        &\frac{d}{dy}A-\star d_A\phi=0,\\
&\frac{d}{dy}\phi-\star(F_A-\phi\wedge\phi)=0.
    \end{split}
    \label{24}
\end{equation}

By a direct computation, we define the Hessian for the functional at an $\slc$ connection $\mathbb{A}=A+i\phi$ as:
\begin{equation}
    \begin{split}
        &\Hess:\Omega_Y^1(\gppc)\rightarrow\Omega_Y^1(\gppc),\\
        &\Hess\left(\begin{array}{l}
             a\\
             b 
        \end{array}
        \right)
        =\left(\begin{array}{l}
             \star d_A b+\star[\phi,a]\\
             \star d_A a-\star[\phi,b] 
        \end{array}
        \right),
    \end{split}
    \label{HessO}
\end{equation}
which is a $y$-independent linearization of $\rm{(\ref{24})}$.

In addition, it is not hard to see that the Hessian operator is horizontal, which means $\Hess$ can also be considered as an operator:
$$\Hess:\Ker\;\dAd\rightarrow \Ker\;\dAd.$$

However, it is much eaiser to consider operators with gauge fixing conditions. We define the extended Hessian $\EHess$ as follows:
\begin{equation}
\begin{split}
    &\EHess:\Omega_Y^1(\gpp)\times\Omega_Y^1(\gpp)\times\Omega_Y^0(\gpp)\times\Omega_Y^0(\gpp)\rightarrow \Omega_Y^1(\gpp)\times\Omega_Y^1(\gpp)\times\Omega_Y^0(\gpp)\times\Omega_Y^0(\gpp),\\
    &\EHess\left(\begin{array}{l}
    a_1\\
    b_1\\
    a_0\\
    b_0
    \end{array}
    \right)
    =\left(\begin{array}{l}
    \star d_A b_1+\star[\phi,a_1]+d_Aa_0-[\phi,b_0]\\
    \star d_Aa_1-\star[\phi,b_1]+d_A b_0+[\phi,a_0]\\
    d_A^{\star}a_1+\star[b_1,\star\phi]\\
    d_A^{\star}b_1-\star[a_1,\star\phi]
    \end{array}
    \right).
\end{split}
\label{EHessO}
\end{equation}

Now we have the following proposition about the  Hessian operator:

\begin{proposition}\label{kernelhess}
(1) $\EHess=\MEH_{(A,\phi,0)}$ and it is a self-adjoint operator.

(2) $\EHess$ is an isomorphism if and only if $\HOC=0$ and $\HZC=0$. 
\end{proposition}
\proof 
(1) is an immediate corollary of Proposition \ref{SASA}.

For (2), using the Hodge theorem, we can decompose the 1-form as
$\Omega^1(\gppc)=\Ker\;\dAd\oplus \Im\;\dA$. Under this decomposition, the extended Hessian operator can be separated into two parts, which we denote as $\EHess=\Hess\oplus S_{\mathbb{A}}.$ Here $S_{\mathbb{A}}$ is defined as follows:

\begin{equation}
    \begin{split}
        &S_{\mathbb{A}}:\Im\; d_{\mathbb{A}}\times\Omega^0(\gpp) \rightarrow \Im\; d_{\mathbb{A}}\times\Omega^0(\gpp),\\
        &S_{\mathbb{A}}\left(\begin{array}{l}
a_1\\
b_1\\
a_0\\
b_0
\end{array}\right)=
    \left(\begin{array}{l}
    d_Aa_0-[\phi,b_0]\\
    d_Ab_0+[\phi,a_0]\\
    d_A^{\star}a_1-\star[\phi,\star b_1]\\
    d_A^{\star}b_1+\star[\phi,\star a_1]
    \end{array}\right).
    \end{split}
\end{equation}

By Hodge theory, we know that $\Ker(\Hess)=\Ker\; d_{\mathbb{A}}\cap \Ker\; d_{\mathbb{A}}^{\star}=H^1_{\mathbb{A}}$ and $\Ker(S_{\mathbb{A}})=H^0_{\mathbb{A}}$
\qed

Therefore, we have the following terminology:
\begin{definition}
The flat connection $\mathbb{A}$ is called non-degenerate if $\HOC$ is zero and acyclic if $\HOC$ and $\HZC$ is zero.
\end{definition}

Now, we will discuss the relation of the extended Hessian and the linearization of the Kapustin-Witten map. Let $(A,\Phi)$ be a solution to the Kapustin-Witten equations. Recall the linearization
\begin{equation}
\begin{split}
\mathcal{L}_{(A,\Phi)}:\Omega_X^1(\gpp)\times\Omega_X^1(\gpp)\rightarrow\Omega_X^2(\gpp)\times\Omega^0(\gpp),\\
\LKW
    \left(\begin{array}{l}
    a\\
    \phi
    \end{array}\right)=
    \left(\begin{array}{l}
    d_{A}a-[\Phi,b]+\star (d_Ab+[\Phi,a])\\
    -\star[a,\star\Phi]+d^{\star}_Ab
    \end{array}\right),
\end{split}
\end{equation}

as well as the gauge fixing operator
\begin{equation}
    \begin{split}
        &\LGF_{(A,\Phi)}:\Omega^1(\gpp)\times\Omega^1(\gpp)\rightarrow\Omega^0(\gpp),\\
        &\LGF_{(A,\Phi)}\left(\begin{array}{l}
a\\
b
\end{array}\right)=
    \begin{array}{l}
    d_A^{\star}a+\star[b,\star\Phi].
    \end{array}
    \end{split}
\end{equation}

and define the following operator $$\DKW:=\LKW+\LGF_{(A,\Phi)}:\Omega^1(\gpp)\times\Omega^1(\gpp)\rightarrow\Omega^2(\gpp)\times\Omega^0(\gpp)\times\Omega^0(\gpp).$$

Let $I$ to be an interval of $\mathbb{R}$ and denote $y$ as the coordinate of $I$. Over $Y\times I$, we have the following identifications:
\begin{equation}
    \begin{split}
        \Omega^1_X(\gpp)&\cong\Omega^0_Y(\gpp)\oplus\Omega^1_Y(\gpp),\\ 
        \alpha_0dy+\alpha_1&\rightarrow \alpha_0 \oplus \alpha_1,\\
        \Omega^2_X(\gpp)&\cong\Omega^1_Y(\gpp)\oplus\Omega^1_Y(\gpp),\\ 
        \alpha_1dy+\alpha_2&\rightarrow \alpha_1 \oplus \star\alpha_2.\\
    \end{split}
    \label{identification}
\end{equation}

Take $(a,b)\in \Omega^1_X(\gpp)\times\Omega^1_X(\gpp)$, under the previous identification, we denote $a=a_0dy+a_1$, $b=b_0dy+b_1$, then we have the following relation of the operator $\DKW$ and the extended Hessian $\MEH$:
\begin{proposition}
For any connection $(A,\Phi=\phi+\phi_ydy)$ over $Y\times I$, if we choose a gauge such that $A$ don't have $dy$ component, then we have
$$\ML_{(A,\Phi)}=-\frac{d}{dy}+\MHN,$$
$$\DKW=-\frac{d}{dy}+\MEH_{(A,\phi,\phi_y)}.$$
\label{Hessian}
\end{proposition}
\proof
Over $Y\times \mathbb{R}$ with volume form $\Vol_Y\wedge dy$, under the identification (\ref{identification}), take $a=a_0dy+a_1$ and $b=b_0dy+b_1$, we have the following computation:
\begin{equation*}
    \begin{split}
        &d_A a-[\Phi,b]=(-\frac{d}{dy}a_1+d_Aa_0-[\phi,b_0]+[\phi_y,b_1])dy+(d_Aa_1-[\phi,b_1]),
    \end{split}
\end{equation*}
and
\begin{equation*}
    \begin{split}
        &d_A b+[\Phi,a]=(-\frac{d}{dy}b_1+d_Ab_0+[\phi,a_0]-[\phi_y,a_1])dy+(d_Ab_1+[\phi,a_1]).
    \end{split}
\end{equation*}

Therefore, we have
\begin{equation*}
    \begin{split}
        &d_A a-[\phi,b]+\star_4(d_A b+[\phi,a])\\
        =&(-\frac{d}{dy}a_1+d_Aa_0-[\phi,b_0]+[\phi_y,b_1]+\star( d_Ab_1+[\phi,a_1]))dy\\
        &+\star(-\frac{d}{dy}b_1+d_Ab_0+[\phi,a_0]-[\phi_y,a_1]+\star(d_Aa_1-[\phi,b_1])).
    \end{split}
\end{equation*}
By our assumption, $\Phi$ doesn't have $dy$ component, we have the following computation:
\begin{equation*}
    \begin{split}
        d_A^{\star_4}b-\star_4[a,\star_4\Phi]
        =d_A^{\star}b_1-\star[a_1,\star\Phi]-\frac{d}{dy}b_0+[\phi_y,a_0]
        =d_A^{\star}b_1+\star[\Phi,\star a_1]-\frac{d}{dy}b_0+[\phi_y,a_0].
    \end{split}
\end{equation*}

Similarly, we have 
\begin{equation*}
    d_A^{\star_4}a+\star_4[b,\star_4\Phi]
        =d_A^{\star}a_1+\star[b_1,\star\Phi]-\frac{d}{dy}a_0-[\phi_y,b_0]
        =d_A^{\star}a_1-\star[\Phi,\star b_1]-\frac{d}{dy}a_0-[\phi_y,b_0].
\end{equation*}

The result follows immediately from our computation.

\qed
\end{subsection}

\end{section}

\begin{section}{Fredholm Theory}
In this section, we will introduce the Fredholm theory for the Kapustin-Witten equations with Nahm Pole boundary condition over manifold with boundary and cylindrical end. See \cite{donaldson2002floer}, \cite{floer1988instanton} for the Fredholm theory of manifold with cylindrical end and see also
\cite{mazzeo2013nahm},\cite{rafe1991elliptic},  for the Fredholm theory on compact manifold with boundary. 
\begin{subsection}{Sobolev Theory on a Manifold with Boundary}
In this subsection, we review Sobolev theory on a manifold with boundary and cylindrical end.

We begin with a suitable functional space. Take $\hat{X}$ to be a compact four-manifold with two boundary components $Y$ and $Z$. Take $X$ to be the four-manifold gluing $\hat{X}$ and $Y\times [0,+\infty)$ along the common boundary $Y$ , $X=\hat{X}\cup_Y Y\times [0,+\infty).$ Denote by $E$ a bundle over $X$ and fix a background connection $\nabla$, let $L^p_k(X,E)$ be the completion of smooth $E$-valued functions on $X$ with respect to the norm
$$
\|f\|_{L^p_k}=\big{(}\sum^k_{i=0}\int_X\vert\nabla^i f\vert^p\big{)}^{\frac{1}{p}},
$$
where $\nabla^k f$ is the symmetric tensor product of $\nabla f$.

For a manifold with boundary, we have the following Sobolev embedding theorem
\begin{proposition}{\rm{(\cite{aubin1982nonlinear}, Thm 2.30, \cite{wehrheim2004uhlenbeck} Appendix)}}
For a compact 4-manifold $\hat{X}$ (with boundary), $k \geq l$, $q\geq p$ and the indices $p$, $q$ are related by
$$k-\frac{4}{p}\geq l-\frac{4}{q},$$
then there is a constant $C_{\hat{X},p,q}$, such that for any section $f$ of a unitary bundle over $\hat{X}$, we have 
$$
\|f\|_{L^q_l}\leq C_{\hat{X},p,q}, \|f\|_{L^p_k}.
$$
\end{proposition}

As a manifold with cylindrical ends has finite geometry, we have the parallel Sobolev embedding theorem for a manifold with boundary and cylindrical ends.
\begin{corollary}
If X is a 4-manifold with boundary and cylindrical ends,  $k \geq l$, $q\geq p$ and the indices $p$,$q$ are related by
$$k-\frac{4}{p}\geq l-\frac{4}{q},$$
then there is a constant $C_{X,p,q}$, such that for any section f of a unitary bundle over $X$,we have 
$$
\|f\|_{L^q_l}\leq \|f\|_{L^p_k}.
$$
\end{corollary}
\proof After identifying the cylindrical ends with $Y\times [0,+\infty)$, we can take open covers $\{U_i\}$ as follows: $U_0:=X/(Y\times [1,+\infty))$ and for $i\geq 1$, $U_i:=Y\times (i-1,i+1)$.

Given a function $f$, let $f_i$ be the restriction of the function to the open cover $U_i$, for $p\leq q$, we have the following inequality,
\begin{equation*}
    \begin{split}
        \|f\|_{L^q_l}&\leq C(\sum^{+\infty}_{i=0}\|f_i\|^q_{L^q_l})^{\frac{1}{q}}\leq C(\sum^{+\infty}_{i=0}\|f_i\|^p_{L^q_l})^{\frac{1}{p}}\leq C(\sum^{+\infty}_{i=0}\|f_i\|^p_{L^p_k})^{\frac{1}{p}}\leq C\|f\|_{L^p_k}.
    \end{split}
\end{equation*}
\qed
\end{subsection}

\begin{subsection}{Elliptic Weight and Nahm Pole Model}
In this subsection, we will discuss the elliptic weights and Fredhlom property of the Kapustin-Witten operator, which is first introduced in \cite{mazzeo2013nahm},\cite{rafe1991elliptic} .

Take $X$ to be a manifold with boundary and cylindrical end, choose a cylindrical neighborhood of $X$ which we will denote as $Y\times (0,1]\subset X$, $Y\times \{0\}=\partial X$. 

Now we shall study the action of $\mathcal{L}$ on the weighted Sobolev space, and so we start by giving the definition of these. 

Choose a smooth function $y:X\rightarrow \mathbb{R}$, which is smaller than 1 and equals the distance function $d(x,\partial X)$ in a neighborhood of $\partial X$.

For any $\lambda\in \mathbb{R}$, we can define the following weighted Sobolev space 
$$y^{\lambda}L^p(X,E):=\{y^\lambda f\vert f\in L^p(X,E)\}.$$ It is easy to see that a suitable norm on this space will be
$$\|f\|_{y^{\lambda}L^p(X,E)}:=(\int_X y^{-\lambda p}|f|^p dx)^{\frac{1}{p}}.$$

Next, we have the following edged Sobolev space which was introduced in \cite{rafe1991elliptic}. Using a local coordinate on $X$ and let $y$ to be the coordinates locally orthogonal to the boundary, we have 
\begin{equation}
H_0^{k,p}(X)=\{f\in L^p(X)\vert(y\partial_{\vec{x}})^{\alpha}(y\partial_y)^j f\in L^p(X),\forall j+\vert\alpha\vert\leq k\}.
\label{edgespapce}
\end{equation}
To be explicit, a suitable norm on this space will be 
$$
\|f\|_{H_0^{k,p}(X)}:=(\int_X\sum_{\forall j+\vert\alpha\vert\leq k}\vert (y\partial_{\vec{x}})^{\alpha}(y\partial_y)^j f\vert^p)^{\frac{1}{p}}.
$$

We define the weighted edge Sobolev space as follows: 
\begin{equation}
y^{\lambda}H^{k,p}_0=\{f=y^{\lambda}f_1\vert f_1\in H^{k,p}_0\}.
\label{weightedgespace}
\end{equation}

Mazzeo and Witten in \cite{mazzeo2013nahm} have the following theorem for the Fredholm property of the Kapustin-Witten operator $\DKW$ with suitable weighted edge Sobolev spaces:

\begin{theorem}{\rm{(\cite{mazzeo2013nahm} Proposition 5.2})}
Let X be a manifold with boudary and $(A,\Phi)$ be a Nahm pole solution to the Kapustin-Witten equations. Let $\DKW$ be the Kapustin-Witten operator \rm{($\ref{KWO}$)} to the Nahm pole solution, and suppose that $\lambda\in (-1,1)$, then the operator
$$\DKW:\LamP H^{1,2}_0(X)\rightarrow \LamM L^2(X)$$ is a Fredholm operator.
\end{theorem}

We also have the following modification of the Theorem for $p\geq 2$ also due to R.Mazzeo:

\begin{theorem}
Let X be a manifold with boudary and $(A,\Phi)$ be a Nahm pole solution to the Kapustin-Witten equations. Let $\DKW$ be the Kapustin-Witten operator \rm{($\ref{KWO}$)} to the Nahm pole solution, and suppose that $\lambda\in (-1,1)$ and $p\geq 2$, then the operator
$$\DKW:\LamP H^{1,p}_0(X)\rightarrow \LamM L^p(X)$$ is a Fredholm operator.
\label{FredholmCompactSpace}.
\end{theorem}
\proof
See Appendix 1.
\qed

Now, we will introduce some basic properties of these weighted Sobolev spaces that will be used in this paper.
\begin{proposition}
$y^lH^{k,p}_0=\{f\in L^p|f \in y^lL^p, \nabla f\in y^{l-1}L^p, \cdots \nabla^k f\in y^{l-k}L^p\}.$
\label{spacelooklike}
\end{proposition}
\proof
By (\ref{edgespapce}), for $g\in H^{k,p}_0$, we know that $g\in L^p,\nabla g\in y^{-1}L^p,\cdots,\nabla^k g\in y^{-k}L^p $. Therefore, by the definition of the weighted edge Sobolev space (\ref{weightedgespace}), for $f\in y^lH^{k,p}_0$, there exist a $g\in H_0^{k,p}$ such that $f=y^l g$. For any positive integers $m,n\leq k$, by the Leibniz rule, we have 
\begin{equation}
\begin{split}
    \nabla_x^m\nabla_y^n f=\sum^l_{i=0}y^{l-i}\nabla^m_x \nabla^{n-i}_y g.
\end{split}
\end{equation}
By the definition of $g$, we have $\nabla^m_x \nabla^{n-i}_y g\in y^{i-m-n}L^p$, therefore $y^{l-i}\nabla^m_x \nabla^{n-i}_y g\in y^{l-m-n}L^p.$ So for $m,n$, we have $\nabla_x^m\nabla_y^n f\in y^{l-m-n}L^p.$

Therefore, for any integer $j$ with $0 \leq j\leq k$, we have $\nabla^j f\in y^{l-j}L^p$.

\qed

In addition, we have the following properties for these spaces.

\begin{proposition}

$(1)$(Different Weight Relation)For any positive integer $p$, if $\lambda_2> \lambda_1$, then
$$y^{\lambda_2}L^p(X)\hookrightarrow y^{\lambda_1}L^p(X).$$ 
Here the $\hookrightarrow$ means the inclusion of Banach space.

$(2)$(Embedding to Usual Sobolev Space) $y^{\lambda} H^{k,p}_0\hookrightarrow y^{\lambda-k}L^p_k,$

$(3)$(H\"{o}lder inequality) $\forall \lambda,\lambda_1,\lambda_2\in\mathbb{R}$ with $\lambda\leq\lambda_1+\lambda_2$, and positive real numbers $p,q,r$, we have 
$$\|fg\|_{y^{\lambda}L^r(X)}\leq\|f\|_{y^{\lambda_1}L^p(X)}\|g\|_{y^{\lambda_2}L^q(X)},$$

\label{embeddinglemma}
\end{proposition}
\proof
(1) Given a function $f$ over a manifold $X$, we have the following inequality:
\begin{equation*}
    \begin{split}
        \|f\|_{y^{\lambda_1}L^p}&=(\int y^{-\lambda_1 p}|f|^p)^{\frac{1}{p}}\\
        &=(\int y^{-\lambda_2 p}y^{(\lambda_2-\lambda_1)p}|f|^p)^{\frac{1}{p}}\\
        &\leq(\int y^{-\lambda_2 p}|f|^p)^{\frac{1}{p}}(\text{here we use $\lambda_2> \lambda_1$}).
    \end{split}
\end{equation*}

(2) Given $f\in y^lH^{k,p}_0$, by Prop \ref{spacelooklike}, for any $0\leq j\leq k$, $\nabla^j f\in y^{l-j}L^p$ and by (1), we have the embedding $y^{l-j}L^p\HRA y^{l-k}L^p.$ The result follows immediately.
    
(3) Given two functions $f$ and $g$ over $X$, we have the following inequality:
\begin{equation*}
    \begin{split}
    \|fg\|_{y^\lambda L^r(X)}&=(\int f^rg^ry^{-r\lambda})^{\frac{1}{r}}\\
    &\leq(\int (fy^{-\lambda_1})^r(gy^{-\lambda_2})^r)^{\frac{1}{r}}\\
    &\leq(\int (fy^{-\lambda_1})^p)^{\frac{1}{p}}(\int (gy^{-\lambda_2})^{q})^{\frac{1}{q}}\\
    &\leq \|f\|_{y^{\lambda_1}L^p(X)}\|g\|_{y^{\lambda_2}L^q(X)}.
    \end{split}
\end{equation*}
\qed

Using these inequalities, we have the following two corollaries:
\begin{corollary}
(1) For any $\lambda \in \mathbb{R}$, $p\geq 2$, we have $$ \LamP H^{1,p}_0(X)\HRA \LamM L^p_1(X)\HRA \LamM L^{2p}(X),$$

(2) For $\La \geq 1-\frac{1}{p}$, we have $$\|fg\|_{y^{\lambda+\frac{1}{p}+1}L^p(X)}\leq\|f\|_{\LamP L^{2p}(X)}\|g\|_{{\LamP}L^{2p}(X)}.$$

(3)For $\La \geq 1-\frac{1}{p}$, $p\geq 2$, we have $\|fg\|_{\LamM L^p(X)}\leq \|f\|_{\LamM L^p_1(X)}\|g\|_{\LamM L^p_1(X)}.$
\label{Sobolevembedding}
\end{corollary}
\proof 
For (1) this is immediate corollary of Proposition \ref{embeddinglemma} combining with the Sobolev embedding $L^p_1\HRA L^{2p}$ with $p\geq 2$. 

(2) For $\lambda\geq 1-\frac{1}{p}$, we have $\lambda+\frac{1}{p}+1\leq \lambda+\frac{1}{p}+\lambda+\frac{1}{p}$. By Proposition \ref{embeddinglemma}, the Holder inequality implies the result.

(3) For $\lambda\geq 1-\frac{1}{p}$, we have $\lambda+\frac{1}{p}-1\leq \lambda+\frac{1}{p}-1+\lambda+\frac{1}{p}-1$. Using Proposition \ref{embeddinglemma} and Sobolev embedding $L^p_1\HRA L^{2p}$, the statement follows immediately.
\qed

\end{subsection}

\begin{subsection}{Fredholmness on Infinite Cylinder}
In this section, we introduce the Fredholm theory for the Kapustin-Witten operator $\DKW$ ($\ref{KWO}$) over the four-manifold $W:=Y^3\times (-\infty,+\infty)$.  

Consider a smooth solution $(A,\Phi)$ to the Kapustin-Witten equations over $W$, which converges in $L^p_1$ norm to acyclic flat $SL(2,\mathbb{C})$ connection over both sides, we have the following proposition:

\begin{proposition}
Under the assumption above, the operator $\DKW:L^p_1(W)\rightarrow L^p(W)$ is a Fredholm operator.
\label{cylindricalendfredholm}
\end{proposition}
\proof
By Proposition \ref{Hessian}, we have $\DKW=-\frac{d}{dy}+\MEH$. As we assume that $(A,\Phi)$ converges to acyclic flat connections over both sides, for some $p>2$, this is a classical result of \cite{lockhart1985elliptic} Theorem 1.3. 
See \cite{floer1988instanton} Proposition 2b.1 for the Yang-Mills case and also  \cite{kronheimer2007monopoles} Proposition 14.2.1 for the $p=2$ version.

\qed
\end{subsection}

\begin{subsection}{The Kapustin-Witten Operator with Acyclic Limit}
As before, denote by $\hat{X}$ a compact four-manifold with two boundary components $Y$ and $Z$. Take $X$ to be the four-manifold obtained by gluing $\hat{X}$ and $Y\times [0,+\infty)$ along the common boundary $Y$ , that is $X=\hat{X}\cup_Y Y\times [0,+\infty).$ Given an $SU(2)$-bundle $P$ over $X$, and a solution $(A,\Phi)$ to the Kapustin-Witten equations ($\ref{KW}$), with Nahm pole boundary condition on $Z$ and which converges in $L^p_1$ norm to acyclic connections over the cylindrical ends for some $p>2$, we have the following theorem:

\begin{proposition}
Under the assumption as above, the Kapustin-Witten operator \rm{($\ref{KWO}$)} $$\DKW:y^{\lambda+\frac{1}{p}}H^{1,p}_0(X)\rightarrow y^{\lambda-1+\frac{1}{p}}L^p(X)$$ is a Fredholm operator.
\label{FredholmOperator}
\end{proposition}
\proof

We will use the parametrix method to prove this theorem. For simplicity, we denote $\DKW$ by $D$ in this proof. To be more explicit, we hope to find two operators $$P:y^{\lambda-1+\frac{1}{p}}L^p(X)\rightarrow y^{\lambda+\frac{1}{p}}H^{1,p}_0(X)$$ and $$R:y^{\lambda+\frac{1}{p}}H^{1,p}_0(X)\rightarrow y^{\lambda-1+\frac{1}{p}}L^p(X) $$ such that $S^l(\rho):=DP(\rho)-\rho,\;S^r(\rho):=RD(\rho)-\rho,$ are two compact operators.

Choose $U_0:=Y\times (T,+\infty)$ and let $U_1$ be a compact cylindrical neighborhood of $\partial X$. By Proposition $\ref{cylindricalendfredholm}$, there exist $P_0, R_0$ such that over $U_0$, we have $DP_0(\rho)=\rho$ and $R_0D(\rho)=\rho$. By the compactness of $\hat{X}$, we can take a finite cover $\{ U_i|i=0 \cdots  n \}$. In each open set $U_i$, by Theorem $\ref{FredholmCompactSpace}$ and the elliptic operator property on the inner open set of the manifold, there exist $P_i,R_i$ and compact operators $S^l_i$ and $S^r_i$ such that $DP_i(\rho)=\rho+S^l_i(\rho),\;R_iD(\rho)=\rho+S^r_i(\rho).$

Denote $X=\bigcup_{i=0}^nU_i$, take $S^l_0(\rho):=0$ and $S^r_0(\rho):=0$. We take a partition of unity $\{ \beta_i \}$ to these covers and define operators 
$P(\rho):=\sum^n_{i=0}\beta_iP_i(\rho)$ and $R(\rho):=\sum^n_{i=0}\beta_iR_i(\rho).$

We have 
\begin{equation}
\begin{split}
    DP(\rho)=&\sum^n_{i=0}\nabla\beta_i\star P_i(\rho)+\sum^n_{i=0}\beta_i DP_i(\rho)\\
    =&\sum^n_{i=0}\nabla\beta_i\star P_i(\rho)+\rho+\sum^n_{i=0}\beta_iS^l_i(\rho).
\end{split}
\end{equation}

Here, $\sum^n_{i=0}S^l_i(\rho)$ is a finite sum of compact operators and it will be compact. 

For the terms $\sum^n_{i=0}\nabla\beta_i\star P_i(\rho)$, recall that 
$$P_i:y^{\lambda-1+\frac{1}{p}}L^p(X)\rightarrow y^{\lambda+\frac{1}{p}}H^{1,p}_0(X).$$

In addition, we know $\nabla\beta_i$ is supported over $\hat{X}$. For functions supported on $\hat{X}$, the norm $y^{\lambda+\frac{1}{p}}H^{1,p}_0(X)$ is equivalent to $L^p_1(X)$ and $y^{\lambda-1+\frac{1}{p}}L^p(X)$ is equivalent to $L^p(X)$. By the compactness of the Sobolev embedding of $L^p_1(X)$ into $L^p(X)$, we know that $\sum^n_{i=0}\nabla\beta_i\star P_i(\rho)$ is also a compact operator.

For the right inverse, we have the following computation:
\begin{equation}
\begin{split}
    RD(\rho)&=\sum^n_{i=0}\beta_iR_iD(\rho)=\sum^n_{i=0}\beta_i(\rho+S^r_i(\rho))=\rho+\sum^n_{i=0}\beta_i S^r_i(\rho).
\end{split}
\end{equation}

Here $\sum^n_{i=0}\beta_i S^r_i(\rho)$ is a finite sum of compact operator thus it is a compact operator.

To summarize, we proved that $D$ is a Fredholm operator.
\qed
$$$$

\end{subsection}
\begin{subsection}{Reducible Limit Connection}
  Given $(A,\Phi)$ a solution to the Kapustin-Witten equations (\ref{KW}) over $W:=Y\times (-\infty,+\infty)$, take $t$ as the coordinate for $(-\infty,+\infty)$. Assume that $(A,\Phi)$ $L^2_2$ converges to a non-degenerate $\slc$  flat connection $(A_{\rho_1},\Phi_{\rho_1})$ when $y\rightarrow -\infty$ and a non-degenerate $\slc$ flat connection $(A_{\rho_2},\Phi_{\rho_2})$ when $y\rightarrow +\infty.$ For $i=1,2$, if either of $(A_{\rho_i},\Phi_{\rho_i})$ is reducible, Proposition \ref{cylindricalendfredholm} is not true since zero can be in the spectrum of the extended Hessian operator $\hat{Q}_{(A_{\rho_i},\Phi_{\rho_i})}$ (\ref{EHessO}).

 Therefore, we hope to use a weight to get rid of the $0$ spectrum and we need to introduce the exponential weight in the cylindrical end. For any real positive number $\alpha$ and norm U, given an arbitrary smooth function $h$ which equals $e^{\alpha t}$ over every cylindrical end $[T,+\infty)\times Y$ and $Y\times(-\infty,-T]$ when $T$ is big enough, we can define the weighted norm by
$$\|f\|_{U_{\alpha}}=\|hf\|_{U}.$$

To be explicit, for a $f$, we denote $\|f\|_{\LamP H^{k,p}_{0,\alpha}}:=\|h f\|_{\LamP H^{k,p}_{0}}$ and 
$\|f\|_{\LamM L^p_{\alpha}}:=\|hf\|_{\LamM L^p}.$

Our operator $\DKW$ can naturally defined over these weighted spaces: $${\DKWA}:\LamP H^{1,p}_{0,\alpha}(W)\rightarrow \LamM L^p_{\alpha}(W).$$

We have the following result:
\begin{proposition}
Under the assumption as above, we can choose $\alpha$ such that the operator $\DKWA:L^p_{1,\alpha}(W)\rightarrow L^p_{\alpha}(W)$ is a Fredhom operator.
\label{RF}
\end{proposition}
\proof
Considering the operator $\DKWA$ over a tube, acting on the 
weighted Sobolev space with $e^{\alpha t}$ as weight function. This is equivalent to an operator $D'$ acting on an unweighted space with the relation
$$D'=\DKWA-\alpha.$$

Therefore, when $\alpha$ is not in the spectrum of the extended Hessian operator over the limit flat connections  $\hat{Q}_{(A_{\rho_i},\Phi_{\rho_i})}$, this is a classical result of \cite{lockhart1985elliptic} Theorem 1.3.
\qed

\begin{subsubsection}{Sobolev Theory for weighted space}
Recall that we denote by $\hat{X}$ a compact four-manifold with two boundary components $Y$ and $Z$. Take $X$ to be the four-manifold obtained by gluing $\hat{X}$ and $Y\times [0,+\infty)$ along the common boundary $Y$ , that is $X=\hat{X}\cup_Y Y\times [0,+\infty).$ Take an $SU(2)$-bundle $P$ over $X$, and a solution $(A,\Phi)$ to the Kapustin-Witten equations ($\ref{KW}$), with Nahm pole boundary condition on $Z$ which converges to a reducible $\slc$-connections over the cylindrical ends.

Over this space $X$, for any real number $\alpha$ and norm $U$ fix a smooth weight function $h$ which approximates $e^{\alpha t}$ over the cylindrical end $[T,+\infty)\times Y.$ When $T$ is large enough, we can define the weighted norm $U_{\alpha}$ as:
$$\|f\|_{U_{\alpha}}:=\|hf\|_{U}.$$

Similarily, we have the Sobolev embedding theorem for the weighted norms.

\begin{proposition}
If X is a 4-manifold with boundary and cylindrical ends,  $k \geq l$, $q\geq p$ and the indices $p$,$q$ are related by
$$k-\frac{4}{p}\geq l-\frac{4}{q},$$ for any given weighted function, there exists a constant $C$, such that for any section f of a unitary bundle over $X$, we have 
$$
\|f\|_{L^l_{q,\alpha}(X)}\leq C\|f\|_{L^k_{p,\alpha}(X)}.
$$
\end{proposition}
\proof This is immediate from the definition of the weighted norm and the usual Sobolev embedding theorem.
\qed
In addition, after fixing a weight function, we have the following inequalities for these weighted edge norms:
\begin{proposition}
(1) For any $\lambda\in\mathbb{R}$, we have $$ {\LamP}H^{1,p}_{0,\alpha}(X)\hookrightarrow {\LamM}L^p_{1,\alpha}(X),$$

(2) For $\lambda\geq1-\frac{1}{p}$
, $\alpha> 0$, we have the following inequality, $$\|fg\|_{{\LamM}L^p_{\alpha}(X)}\leq\|f\|_{\LamM L^{2p}_{\alpha}(X)}\|g\|_{y^{\lambda_0-\frac{1}{2}}L^{2p}_{\alpha}(X)}.$$
\label{SoI}
\end{proposition}
\proof
The statement in (1) is immediately using Corollary \ref{Sobolevembedding} and the definition of weighted Sobolev space.
By Corollary \ref{Sobolevembedding}, we know $$\|fg\|_{\LamM L^p_{\alpha}(X)}\leq\|f\|_{\LamM L^{2p}(X)}\|g\|_{\LamM L^{2p}_{\alpha}(X)}.$$

In addition, as we assume $\alpha> 0$, we know that $$\|f\|_{\LamM L^{2p}(X)}\leq \|f\|_{\LamM L^{2p}_{\alpha}(X)}.$$
The statement in (2) follows immediately.
\qed

\end{subsubsection}
\begin{subsubsection}{Fredholm Property for the Reducible Limit}

\begin{proposition}
Under the assumpution above, there exist $\alpha$ such that the operator $${\DKWA}:y^{\lambda+\frac{1}{p}} H^{1,p}_{0,\alpha}(X)\rightarrow y^{\lambda-1+\frac{1}{p}} L^p_{\alpha}(X)$$ is a Fredhlom operator.
\label{reducefredholm}
\end{proposition}
\proof
The main difference between the reducible limit and the acyclic limit is the behavior of the operator $\DKW$ over the cylindrical end. In the reducible case, we use Proposition $\ref{RF}$ over the cylindrical ends to get a  parametrix and the results follow similarly as Theorem \ref{FredholmOperator}.
\qed
\end{subsubsection}

\end{subsection}

\begin{subsection}{The Index}
Now we will give an explicit computation of the index for a manifold with cylindrical end.

For a compact manifold with boundary, in \cite{mazzeo2013nahm}, Mazzeo and Witten have a computation of the index of $\DKW$ where $(A,\Phi)$ is a Nahm pole solution to the Kapustin-Witten equations (\ref{KW}):

For our case, let $X$ be a manifold with boundary $Z$ and cylindrical end which is identified with $Y\times [0,+\infty)$. We finish a parallel computation for the index of a manifold with boundary and cylindrical end. Under the previous assumption, by Proposition \ref{FredholmOperator}, \ref{reducefredholm}, we can define the index for these Fredholm operator.

If $(A,\Phi)$ has an acyclic limit, we denote by $\Ind_X(P)$ the index of $\DKW$ in the setting of Theorem \ref{FredholmOperator}:
\begin{equation}
\Ind_X(P):=\dim\; \Ker \DKW- \dim\; \Coker \DKW.
\end{equation}

For compact manifold with Nahm pole boundary condition, Mazzeo and Witten has the following computations:
\begin{proposition}[{\cite{mazzeo2013nahm}}\;Proposition 4.2]
Let $X$ be a compact manifold with boundary, let $P$ be an SU(2) bundle over $X$, let $(A,\Phi)$ be a solution to the Kapustin-Witten equations with acyclic limit which satisfies the Dirichlet boundary condition over the boundary, then 
$$\Ind P=-3\chi(X).$$
\end{proposition}

After a modification of their proof, we have the following computation for solutions have irreducible limits:

\begin{proposition}
Let $X$ be a manifold with boundary and cylinderical ends, let $P$ be an SU(2) bundle over $X$, let $(A,\Phi)$ be a solution to the Kapustin-Witten equations with acyclic limit which satisfies the Dirichlet boundary condition over the boundary, then 
$$\Ind P=-3\chi(X).$$
\label{indexdirichilet}
\end{proposition}
\proof Let $(A,\Phi)$ be a solution to the Kapustin-Witten equation with Dirichlet boundary condition, then the index of the operator $\DKW$ corresponds to the relative boundary condition for the Gauss-Bonnet operator, which the index is equals to $-3\chi(X,\partial X)$ and by Poincare duality, $\chi(X,\partial X)=\chi(X).$

\qed

\begin{proposition}[{\cite{mazzeo2013nahm}}\;Proposition 4.3]
Under the same assumption as Proposition \ref{indexdirichilet}, let $(A,\Phi)$ be a solution to the Kapustin-Witten equations (\ref{KW}) 
with the Nahm pole boundary condition over $\partial X$, then

$$\Ind P=-3\chi(X).$$
Here $\chi(X)$ is the Euler characteristic of $X$.
\label{indexcompute}
\end{proposition}
\proof First consider the special case of $(0,1]\times Z$ with the product metric. Let $(A_0,\Phi_0)$ be a connection pair satisfying the Nahm pole boundary condition over $\{0\}\times Z$ and regular over $\{1\}\times Z$. Let $\mathcal{D_{N,R}}$ be the elliptic operator corresponding to this. By \cite{mazzeo2013nahm} (3.12), as $\mathcal{D_{N,R}}$ is pseudo skew-Hermitian, we have Ind $\mathcal{D_{N,R}}$=0. 

Now consider a general $(A,\Phi)$ over $X$ satisfying the Nahm pole boundary condition, we choose a tubular neighborhood of $X$ near the boundary and identify it with $Y\times (0,1]$. Let $\mathcal{D_{N,R}}$ be the restriction of the operator $\DKW$ to $Y\times (0,1]$ and let $\mathcal{D_R}$ be the restriction of $\DKW$ over the complement of $Y\times (0,1]$ in $X$. By a standard excision theorem of index \cite{booss2012elliptic},\;Prop 10.4, we obtain 

$$\Ind\; \DKW=\Ind\;\mathcal{D_{N,R}}+\Ind\; \mathcal{D_R}.$$

By the previous argument, we have $\Ind\;\mathcal{D_{N,R}}=0$ and by Proposition \ref{indexdirichilet}, we know $\Ind\;\mathcal{D_{R}}=\chi(X).$ Thus we have  $\Ind\; \DKW=\chi(X).$

\qed

If $(A,\Phi)$ is a reducible but non-degenerate limit as in the case of Proposition \ref{reducefredholm}, we denote by $\Ind_X(P,\alpha)$ the index of $\DKWA$ with respect to weight $\alpha$:
\begin{equation}
\Ind_X(P,\alpha):=\dim\; \Ker \DKWA- \dim\; \Coker \DKWA.
\end{equation}

Take $\alpha^+$ to be a real number that is slightly bigger than $0$ and below the positive spectrum of $\EHessr$ and $\alpha^-$ to be a real number that is slightly smaller than $0$ and above the negative spectrum of $\EHessr$, then we have two indices:

\begin{equation}
\Ind^+_X(P):=\Ind(P,\alpha^+),\; \Ind^-_X(P):=\Ind(P,\alpha^-).
\end{equation}

For $i=1,2$, suppose $X_i$ is a manifold with boundary $Z_i$ and cylindrical end $Y_i\times [0,+\infty)$ and $Y_1=Y_2$, $P_i$ is an $SU(2)$ bundle over $X_i$. Let $(A_i,\Phi_i)\in \mathcal{C}_{P_i}$ converge to the same $\slc$ flat connection $\rho$, then we can glue these two manifolds along the common boundary $Y$ to form $X^{\sharp}$ and a bundle $P^{\sharp}$. As it converges to the same flat connection $\rho$, we can define a new pair $(A^{\sharp},\Phi^{\sharp})$ on $P^{\sharp}$. 

We denote by
$\Ind_{X^{\sharp}}(P^{\sharp})$ the index of the operator $\DKWX$ on $X^{\sharp}$. 

If $(A_i,\Phi_i)$ both have acyclic limits, then we have the following gluing relation of these indices:

\begin{proposition}
$\Ind_{X^{\sharp}}(P^{\sharp})=\Ind_{X_1}(P_1)+\Ind_{X_2}(P_2).$

\label{glueindex}
\end{proposition}
\proof See  
\rm{(\cite{donaldson2002floer} Proposition 3.9)}, the same arguments also work in our case.
\qed

If $(A_i,\Phi_i)$ is reducible but non-degenerate, we can choose $\alpha_1=-\alpha_2> 0$ whose absolute value is smaller than the smallest absolute value of eigenvalues of $\EHessr$. We denote $\Ind_{X_1}^{+}(P_1)$ to be the index corresponding to the weight $\alpha_1$ and $\Ind_{X_2}^{-}(P_2)$ to be the index corresponding to the weight $\alpha_2$.

We have following Proposition:
\begin{proposition}
$\Ind_{X^{\sharp}}(P^{\sharp})=\Ind_{X_1}^{+}(P_1)+\Ind_{X_2}^{-}(P_2).$

\label{addindex}
\end{proposition}
\proof See  
\rm{(\cite{donaldson2002floer} Proposition 3.9)}, the same arguments also work in our case.
\qed

Moreover, for the index over the same bundle with small positive and negative weights, we obtain:
\begin{proposition}{\rm{(\cite{donaldson2002floer} Proposition 3.10)}}
$\Ind_X^+(P)-\Ind_X^-(P)=-\dim\;ker \EHessr$. 
\label{diff}
\end{proposition}

Now we will do some explicit computation of indices. Consider the model case $W$ to be the 'flask' manifold obtained by gluing a punctured 4-sphere to a tube $S^3\times \mathbb{R}$. Consider the trivial bundle and the trivial connection over this space, we have the following Lemma for indices:
\begin{lemma}
The indices for $W$ are $$\Ind^+_W=-6,\;\Ind^-_W=0.$$
\end{lemma}
\proof
Obviously, $W\sharp W$ is diffeomorphic to $S^4$. In addition, by Proposition \ref{indexcompute}, we know the index of operators $\DKW$ over $S^4$ is -6. Therefore, by the gluing relation of the index Proposition \ref{addindex}, we have 
$$\Ind_W^++\Ind_W^-=-6.$$

In addition, by Proposition \ref{diff}, we have 
$$\Ind_W^+=\Ind_W^--6.$$
Combining these two index formulas, we get the result we want.
\qed

\begin{corollary}
If the cylindrical end of X has the form $S^3\times [0,+\infty)$, we have 
\begin{equation}
\Ind_X^+(P)=-3\chi(X)-3,\; \Ind_X^-(P)=-3\chi(X)+3.
\end{equation}

\end{corollary}
\proof
Denote by $\bar{X}$ a smooth compactification of $X$ over the tube $S^3\times [0,+\infty)$ and denote by $\bar{P}$ the extension of the bundle P.

By Proposition \ref{addindex}, we have 
$$\Ind_X^+(P)+\Ind^-_W=\Ind_{\bar{X}}(\bar{P}),$$
$$\Ind_X^-(P)+\Ind^+_W=\Ind_{\bar{X}}(\bar{P}).$$
In addition, by Proposition \ref{indexcompute}, we know that $$\Ind_{\bar{X}}(\bar{P})=-3\chi(\bar{X})=-3\chi(X)-3.$$ We get the result we want.
\qed

\end{subsection}
\end{section}

\begin{section}{Moduli Theory}
In this section, we will introduce the moduli theory for the solutions to the Kapustin-Witten equations (\ref{KW}) with Nahm pole boundary condition.
\begin{subsection}{Framed Moduli Space}
In this section, we will give suitable norms to define the moduli space.

Let $X$ to be a smooth 4-manifold with 3-manifold boundary $Z$ and cylindrical end which is identified with $Y\times (0,+\infty)$. Now suppose $P$ is an $SU(2)$ bundle over $X$, $\gpp$ is the associated adjoint bundle, $\mathcal{A}_P$ is the set of all $SU(2)$ connections on $P$ and $\CP:=\mathcal{A}_P\times \Omega^1(\gpp)$. 

For $i=1,2,3$, fix an orthogonal frame $\{e_i\}\in T^{\star}Y$, choose a reference connection pair $(A_0,\Phi_0)\in \CP$ to the Kapustin-Witten equations ($\ref{KW}$). We require that $(A_0,\Phi_0)$ satisfies the Nahm pole boundary condition for this frame, thus there exists $\{t_i\}\in \gpp$ such that the expansion of $\Phi$ when $y\rightarrow 0$ is $\Phi\sim\frac{\sum_{i=1}^3e_it^i}{y}$. In addition, fix a flat acyclic $\slc$ representation $\rho$. If we denote $Y_T:=\{T\}\times Y\subset (0,+\infty)\times Y$, we assume that $(A_0,\Phi_0)$ convergence to $\rho$ in $L^p_1$ norm for some $p>2$.

Now, we will introduce a suitable configuration space we hope to study. 

Given real numbers $p,\lambda$ with $p> 2$ and $\lambda\in [1-\frac{1}{p},1)$ , we have the following definition:
\begin{definition}
Given a smooth Nahm pole solution $(A_0,\Phi_0)$, we define the framed configuration space $\FCS$ as follows:
\begin{equation}
\FCS:=\{(A_0,\Phi_0)+(a,b)\;|\;(a,b)\in \LamPP H^{1,p}_0(\Omega^1(\gpp)\times\Omega^1(\gpp))\}.
\label{configurationframe}
\end{equation}
Here $\LamPP H^{1,p}_0(\Omega^1(\gpp)\times\Omega^1(\gpp))$ is the completion of smooth 1-forms with respect to the norm $\LamPP H^{1,p}_0$.
\end{definition}

We have some basic properties of the framed configuration space:
\begin{proposition}

(1) Any $(A,\Phi)\in \FCS$ satisfies the Nahm pole boundary condition.

(2) Assume $X$ has non-vanishing boundary and cylindrical end which identified with $Y\times (0,+\infty)$. Let $P$ be an $SU(2)$ bundle over it, let $(A_1,\Phi_1)$ be a connection pair satisfying the Nahm pole boundary condition, denote $\frac{\sum_{j=1}^3e_jt^1_j}{y}$ is the leading part of $\Phi_1$. If $(A_1,\Phi_1-\frac{\sum_{j=1}^3e_jt^1_j}{y})\in H^{1,p}_0(X),$ then there exists a global gauge transformation $g\in\mathcal{G}$ such that $g(A_1,\Phi_1)\in\FCS.$
\label{propositionframe}
\end{proposition}
\proof
For (1), as for $\lambda\in [1-\frac{1}{p},1)$ ,$p>  2$, any differential form which blows-up as $y^{-1}$ is not contained in $\LamPP H^{1,p}_0$, the result follows immediately.

For (2), for $i=0,1$, $(A_i,\Phi_i)$ both satisfies the Nahm pole boundary condition (Definition \ref{NahmPole}). Then there exists orthogonal basis ${e_j}\in T^{\star}Y$ and ${t^i_j}\in\gpp$ for $j=1,2,3$ such that $[t^i_{j_1},t^i_{j_2}]=\epsilon_{j_1j_2j_3}t^i_{j_3}$ where $\epsilon_{j_1j_2j_3}$ is the Kronecker symbol of $j_1,j_2,j_3$. In addition, let the asymptotic expansion of $\Phi_i$ at $y=0$ to be $\frac{\sum_{j=1}^3e_jt_j^i}{y}+\mathcal{O}(y).$ By the commutation relation of $t_j^i$, there exists a $\hat{g}:Z\rightarrow SU(2)$ such that $\hat{g}(\frac{\sum_{j=1}^3e_jt_j^1}{y})\hat{g}^{-1}=\frac{\sum_{j=1}^3e_jt_j^0}{y}$. By the Hopf theorem, the homotopy type of maps from $Y^3$ to $SU(2)$ is totally determined by the degree. Given $\hat{g}:Z\rightarrow SU(2)$, consider $Y_{T_0}=Y\times \{T_0\}\subset X$, we can choose a band to connect $Z$ and $Y$ which is homeomorphic to $Z\sharp Y_{T_0}$. We can choose a map $\hat{g}':Y_{T_0}\rightarrow SU(2)$ whose degree equals minus degree of $g$ and extend these two maps to $Z\sharp Y_{T_0}$ and denote as $\tilde{g}$. Then $\tilde{g}$ has degree zero and can be extended to whole $X$, which we denote as $g$.

By the assumption that $(A_i,\Phi_i)$ are smooth, we have $g(A_1,\Phi_1)-(A_0,\Phi_0)\in \LamPP H^{1,p}_0.$ 
\qed

\begin{remark}
For a general compact 4-dimensional manifold with boundary, Proposition \ref{propositionframe} is not true. For $D^4$, the 4-dimensional unit disc, a choice of frame gives a map from $S^3\rightarrow SO(3)$ which is $\pi_3(SO(3))$. Two frames corresponding to different elements in $\pi_3(SO(3))$ can not be globally gauge equivalent. 
\end{remark}

As we fixed a base connection $(A_0,\Phi_0)$ to define the framed configuration space, we can also consider the gauge group that preserves the frame.
\begin{definition}
The framed gauge group $\mathcal{G}^{fr}$ is defined as follows:
\begin{equation}
\FGG:=\{g\in Aut(P)\;|\;g|_{Y}=1\}.
\end{equation}
\end{definition}

Given $g\in \FGG$, the action of $g$ on $(A_0,\Phi_0)$ will be 
\begin{equation}
g(A_0,\Phi_0)=(A_0-d_{A_0}g\;g^{-1},g\Phi_0 g^{-1}).
\end{equation}
Then, we have 
\begin{equation}
g(A_0,\Phi_0)-(A_0,\Phi_0)=(-d_{A_0}g\;g^{-1},\; [g, \Phi_0]g^{-1}).
\end{equation}
We consider the following weighted frame gauge group $\FGGG$:
\begin{equation}
\FGGG=\{g\in\FGG\;|\;d_{A_0}g\;g^{-1}\in\LamPP H^{1,p}_0(\Omega^1),\; [g, \Phi_0]g^{-1}\in \LamPP H^{1,p}_0(\Omega^1) \}.
\label{framegaugegroup}
\end{equation}

For convenience, we denote $d^0(\xi):=\DZ(\xi)=(d_{A_0}\xi,[\Phi_0,\xi])$ and we have the following lemma on the weighted frame gauge group $\FGGG$.
\begin{lemma}
$\FGGG=\{g\in \FGG\;|\;\nabla_0 g\in {\LamPP} L^p,\; \nabla_0^2 g\in {\LamMM} L^p,\; [\Phi_0,g]\in {\LamP} L^p,\; \nabla_0[\Phi_0,g]\in {\LamM} L^p\}.$ 
\label{Liegroup}
\end{lemma}
\proof Obviously $\{g\in \FGG\;|\;d^0 g\in {\LamPP} L^p,\;\nabla_0(d^0) g\in {\LamMM} L^p \}\subset \FGGG$. For the other side, we argue as follows: take $\alpha:=d_{A_0}gg^{-1}$, then $\alpha\in \LamPP H_0^{1,p}\hookrightarrow \LamMM L^{2p}.$ By $d_{A_0}g=\alpha g$, we have $\nabla_0g\in \LamMM L^{2p}$ since the pointwise norm of $g$ is 1. In addition, $\nabla_0d_{A_0}g=(\nabla_0\alpha)g+\alpha\nabla_0 g.$ As $\nabla_0\alpha\in \LamMM L^p$ and the pointwise norm of $g$ is $1$, we have $\nabla_0\alpha g\in \LamMM L^p$. 
As $\alpha\in \LamMM L^{2p}$, $\nabla_0 g\in \LamMM L^{2p}$, we have $\alpha\nabla_0 g\in \LamMM L^p$. For $[g, \Phi_0]g^{-1}\in \LamPP H^{1,p}_0(\Omega^1)$, we have $[g, \Phi_0]g^{-1}\in \LamPP L^p$. Letting $\beta=[g,\Phi_0]g^{-1}$, then $\beta\in \LamP L^p$ implies $\beta g=[g,\Phi_0]\in \LamP L^p.$ In addition, $\beta \in \LamP H^{1,p}_0$ implies $\nabla_0\beta\in \LamP L^p$ and $\beta\in \LamM L^{2p}$. In addition, we have $\nabla_0g\in \LamM L^p$, thus $\nabla_0[g,\Phi_0]=\nabla_0(\beta\;g)=(\nabla_0\beta)g+\beta\nabla_0 g\in \LamM L^p.$

Therefore,
$\FGGG\subset \{g\in \FGG\;|\;\nabla_0 g\in {\LamPP} L^p,\;\nabla_0^2 g\in {\LamMM} L^p \}.$

\qed

Thus we can rewrite $\FGGG$ as $\FGGG=\{g\in \FGG\;|\;d^0 g\in {\LamPP} L^p,\; \nabla_0d^0 g\in {\LamMM} L^p\}.$

\begin{lemma}
The space $y^{\lambda+\frac{1}{p}+1}H_0^{2,p}(\gpp)$ is an algebra and $\LamPP H_0^{1,p}(\gpp)$ is a module over this algebra.
\label{Liealgebra}
\end{lemma}
\proof
For the algebra statement, we only need to prove $u_1u_2\in {\LaPP}H_0^{2,p}(\gpp)$, or equivalently, $u_1u_2\in {\LaPP}L^p$, $\nabla_0(u_1u_2)\in \LamPP L^p$, $\nabla_0^2(u_1u_2)\in \LamMM L^p$.

Since $u_i\in \LaPP H^{2,p}_0$, we have $u_i\in \LaPP L^p$, $\nabla_0u_i\in\LamP L^p$ and $\nabla_0^2u_i\in \LamM L^p.$

By Proposition \ref{\Lone}, $u_i\in \LaPP H^{1,p}_0\HRA \LamPP L^p_1\HRA \LamPP L^{2p}.$ By Corollary \ref{\Ltwo}, we have $u_1u_2\in \LaPP L^p.$ In addition, we know $\nabla_0 u_i\in\LamP H^{1,p}_0\HRA\LamMM L^{2p}$ and $u_i\in \LamPP L^{2p}$. As $\lambda\geq 1-\frac{1}{p}$, we have $\lambda+\frac{1}{p}+\lambda+\frac{1}{p}-1\geq \lambda+\frac{1}{p}$. By the H\"older inequality in Proposition \ref{\Lone}, we have $\nabla_0u_1 u_2\in \LamPP L^p$. 

For $u_i\in \LaPP H^{2,p}_0$, as $p>  2$ and $\lambda\geq 1-\frac{1}{p}$, we have $u_i\in \LamM L^p_2\HRA C^0$ and $\nabla^2_0u\in \LamMM L^p$. Therefore, we have $\nabla^2_0u_1 u_2\in \LamMM L^p.$

The module statement can also be proved in a similar way.
\qed

Fix a base point $p_0\in X$ and define a system of neighborhoods of the identity in $\GP$ as 
\begin{equation}
U_{\epsilon}=\{g\in \FGG\;|\;\|d^0 g\|_{{\LamPP L^p}}\leq \epsilon,\;\|\nabla_0 d^0 g\|_{{\LamMM L^p}}\leq \epsilon,\;|g(p_0)-1|\leq\epsilon \}.
\end{equation}

This topology is independent of the base point $p_0$. 

With the previous lemma, we can establish a Lie group structure on $\FGG$:
\begin{corollary}
(1) $\FGGG$ is a Lie group with Lie algebra
$$Lie(\FGGG)=\LaPP H_0^{2,p}(\gpp).$$
(2) $\FGGG$ acts smoothly on $\FCS.$
\end{corollary}
\proof
This result follows immediately from Lemma $\ref{Liegroup}$ and Lemma $\ref{Liealgebra}.$
\qed

Now, we have the following proposition of the framed configuration space $\FCS$ and the framed gauge group $\FGGG$:
\begin{proposition}
For any $(A,\Phi)\in \FCS$, we have $KW(A,\Phi)\in \LamMM L^p(\Omega^2\oplus\Omega^0),$
\end{proposition}
\proof 
By the definition of $\FCS$, there exists $(a,b)\in \LamPP H^{1,p}_0(\Omega^1(\gpp)\times\Omega^1(\gpp))$ such that $(A,\Phi)=(A_0,\Phi_0)+(a,b).$

By Proposition \ref{qe}, we have $KW(A,\Phi)=KW(A_0,\Phi_0)+\mathcal{L}^1(a,b)+\{(a,b),\;(a,b)\}.$ Here $\{(a,b),\;(a,b)\}$ is a quadratic term. By theorem \ref{FredholmCompactSpace}, we have $\mathcal{L}^1(a,b)\in \LamMM L^p(\Omega^2\times\Omega^0).$ By the embedding $\LamPP H_0^{1,p}\HRA \LamMM L^p_1\HRA \LamMM L^{2p}$, we have $\{(a,b),(a,b)\}\in \LamMM L^{p}.$ 
\qed

Now we will study the behavior of the gauge group (\ref{framegaugegroup}) over the cylindrical end. We have the following proposition which describes the limit behavior of the group $\FGGG$. We need the hypothesis that $\rho$ is acyclic. Let $(A_\rho,\Phi_\rho)$ be the flat $\slc$ connection associate to $\rho$. 

Recall that $d^0_{(A,\Phi)}(\xi)=(d_{A}\xi,[\Phi,\xi])$ and $\rho$ acyclic implies $\Ker\;d^0_{(A_\rho,\Phi_\rho)}=0$ and the connection $d_{A_\rho}$ itself may still be a reducible $SU(2)$ connection.

We have the following lemma over the cylindrical end:

\begin{lemma}
Suppose $X$ is a manifold with boundary and cylindrical end which is identified with $Y\times (0,+\infty)$, for $(A_0,\Phi_0)$ a reference connection which $L^p_1$ converges to $(A_\rho,\Phi_\rho)$ over the cylindrical end for $p>2$, then for $T$ is large enough, we have 

(1) $\DZ:L^p_2(Y\times (T-1,T+1))\rightarrow L^p_1(Y\times (T-1,T+1))$ is injective,

(2) $\|\xi\|_{L^p_2(Y\times (T-1,T+1))}\leq C\|\DZ\xi\|_{L^p_1(Y\times (T-1,T+1))}$.
\label{irreducibletubular}
\end{lemma}
\proof For convenience, during the proof, we write $L^p_k$ short for $L^p_k(Y\times (T-1,T+1))$.

(1) Denote $(a,b)=(A_\rho,\Phi_\rho)-(A_0,\Phi_0)$, then for any $\xi\in L^p_2$, we have $$\|d^0_{(A_\rho,\Phi_\rho)}\xi-\DZ\xi\|_{L^p_1}\leq C(\|[a,\xi]\|_{L^p_1}+\|[b,\xi]\|_{L^p_1})\leq C(\|a\|_{L^p_1}+\|b\|_{L^p_1})\|\xi\|_{L^p_2}.$$ If $\xi\in \Ker\;\DZ$, we obtain $$\|d^0_{(A_\rho,\Phi_\rho)}\xi\|_{L^p_k}\leq C(\|a\|_{L^p_1}+\|b\|_{L^p_1})\|\xi\|_{L^p_2}.$$ For $T$ is large enough,  $C(\|a\|_{L^p_1}+\|b\|_{L^p_1})$ is smaller than the operator norm of $d^0_{(A_0,\Phi_0)}$, which implies $\xi=0$.

(2) If the inequality is not true, then there exists a sequence $\{\xi_n\}$ with $$\|\xi_n\|_{L^p_1}=1,\;\lim_{n\rightarrow \infty}\|\DZ\xi_n\|_{L^p_1}=0,$$ which implies $\|\xi_n\|_{L^p_2}$ is bounded.

Then, $\xi_n$ weak converges to $\xi_{\infty}$ in $L^p_2$ and strongly converges to $\xi_{\infty}$ in $L^p_1$, which implies $\|\DZ\xi_{\infty}\|_{L^p_1}=0$ and $\|\xi_{\infty}\|_{L^p_1}=1$. As $\Ker\;\DZ=0$, we have $\xi_{\infty}=0$, contradicting $\|\xi_\infty\|_{L^p_1}=1.$

\qed

We have the following corollary:

\begin{corollary}
If over the cylindrical end, $(A_0,\Phi_0)$ converges in $L^p_1$ norm to $(A_\rho,\Phi_\rho)$ and $(A_\rho,\Phi_\rho)$ is an irreducible flat $\slc$ connection, then $\Ker\;d^0_{(A_0,\Phi_0)}=0.$
\label{irreducible1}
\end{corollary}

\proof
As before, we denote $d^0:=d^0_{(A_0,\Phi_0)}$. By Kato inequality, we know for $\xi\in\Omega^0$, we have the pointwise estimate $$d|\xi|\leq|d_{A_0}\xi|\leq|d^0(\xi)|.$$
Therefore, $\xi\in \Ker d^0$ implies $|\xi|$ is a constant. In addition, by Lemma \ref{irreducibletubular}, we know $\xi=0$ over $Y\times (T-1,T+1)$ when $T$ is large enough, therefore, we have $\xi$ is identically zero.
\qed

\begin{proposition}
If the limiting connection $\rho$ is irreducible, then for $T$ is large enough, there is a constant C such that for any section $\xi\in\gpp$, with $d^0\xi\in L^p(Y\times [T,+\infty))$ and $\nabla_0d^0\xi\in L^p(Y\times [T,+\infty))$, we have 

(1) $|\xi|\rightarrow 0$ at the cylindrical end and $\sup|\xi|\leq C(\|d^0\xi\|_{L_1^p(Y\times [T,+\infty))}),$

(2) Either $|g(x)-1|\rightarrow 0$ or $|g(x)+1|\rightarrow 0$ as $x$ tends to infinity in X.

\end{proposition}
\proof 

(1) For an integer $k>T+1$, over a band $B_k:=Y\times(k-1,k+1)$, we have $$\|\xi\|_{C^0(B_k)}\leq \|\xi\|_{L^p_2(B_k)}\leq C\|d^0\xi\|_{L^p_1(B_k)}.$$ The statement follows immediately.

(2) Denote by $g_k$ the restriction of $g$ to the band $B_k$. After identifying different bands with $Y^3\times (-1,1)$, we can consider $\{g_k\}$ a sequence of gauge transformation over $Y^3\times (-1,1)$ with $\|\nabla_0d^0g_k\|_{L^p},\;\|\nabla_0g_k\|_{L^p}$ converging to zero. As the pointwise norm of g is always 1, by Rellich lemma, $g_k$ strongly converges to $g_{\infty}$ in $L_1^p$ which implies $d^0g_{\infty}=0$. By our assumption, we have $g_{\infty}=\pm 1.$
\qed

Now, we can define a framed quotient space as follows:
\begin{definition}
We define $\FBS$ as the following weighted framed moduli space:
$$\FBS=\FCS/\FGGG.$$
\end{definition}
In addition, we have the following definition of the moduli space:

\begin{definition}
The framed moduli space $\MM_{p,\lambda}^{fr,\rho}(X)$ is defined as follows:
\begin{equation}
\MM_{p,\lambda}^{fr,\rho}(X)=\{(A,\Phi)\in\FCS|KW(A,\Phi)=0\}/\FGGG.
\label{modulispace1}
\end{equation}
\end{definition}

We have the following basic properties of the framed moduli space:
\begin{proposition}
(1) For any $(A,\Phi)\in \MM_{p,\lambda}^{fr,\rho}$ satisfies the Nahm pole boundary condition.

(2) Any $(A,\Phi)\in \MM_{p,\lambda}^{fr,\rho}$ converges to $\rho$ in $L^p_1$ norm.
\end{proposition}
\proof
(1) is an immediate consequences of Proposition \ref{propositionframe}. (2) is a consequence of the definition of $\FCS.$ 
\qed

\end{subsection}

\begin{subsection}{Slicing Theorem}
Now we study the local properties of the moduli space and we will assign a suitable norm to the Kuranishi complex (\ref{complex}).

Define $(\Omega^0,\lambda,k,p)$ as follows:
$$(\Omega^0,\lambda,k,p):=y^{\lambda} H^{k,p}_0(\Omega^0(\gpp)).$$
Here the notation $y^{\lambda} H^{k,p}_0(\Omega^0(\gpp))$ means the completion of the smooth sections of $\Omega^0(\gpp)$ in the norm $y^{\lambda} H^{k,p}_0$. 
Similarly, we define
$$(\Omega^1\times\Omega^1,\lambda,k,p):=y^{\lambda} H^{k,p}_0(\Omega^1(\gpp)\times\Omega^1(\gpp))$$ and $$(\Omega^2\times\Omega^0,\lambda,k,p):=y^{\lambda}H^{k,p}_0(\Omega^2(\gpp)\times\Omega^0(\gpp)). $$

Now we rewrite the Kuranishi complex (\ref{complex}) at the point $(A_0,\Phi_0)$  with respect to the new norm as follows:

\begin{equation}
0\rightarrow (\Omega^0,\lambda+1+\frac{1}{p},2,p)\xrightarrow{d^0_{(A_0,\Phi_0)}}(\Omega^1\times\Omega^1,\lambda+\frac{1}{p},1,p)\xrightarrow{\mathcal{L}_{(A_0,\Phi_0)}}(\Omega^2\times\Omega^0,\lambda+\frac{1}{p}-1,0,p)\rightarrow 0,
\label{newcomplex}
\end{equation}
Here we only considered $\lambda\in[1-\frac{1}{p},1).$

We have the following Proposition for the operator $d^{0}_{(A_0,\Phi_0)}$:

\begin{proposition}
The operator $\DZ$:
$$\DZ:y^{\lambda+1+\frac{1}{p}} H^{2,p}_0(\Omega^0(\gpp))\rightarrow y^{\lambda+\frac{1}{p}} H^{1,p}_0(\Omega^1(\gpp)\times\Omega^1(\gpp))$$
is a closed operator.
\label{DZclose}
\end{proposition}
\proof see Appendix 2.
\qed

\begin{corollary}
$y^{\lambda+\frac{1}{p}} H^{2,p}_0(\Omega^1(\gpp)\times\Omega^1(\gpp))=\mathrm{Im\;} \DZ\oplus (\mathrm{Ker\;} \DZS \cap y^{\lambda+\frac{1}{p}} H^{1,p}_0)$
\label{DZsplit}
\end{corollary}
\proof Let $x=(x_1,x_2)\in \Omega^1(\gpp)\times\Omega^1(\gpp)$, by definition of $\DZ$, we have $$\lan \DZ\xi,x \ran =\lan d_{A_0}\xi,x_1\ran +\lan [\Phi_0,\xi],x_2\ran $$ where $\lan\;,\;\ran $ means the $L^2$ inner product. 

Integrating by parts, we have $$\lan d_{A_0}\xi,x_1\ran =\lan \xi,d_{A_0}^{\star}x_1\ran -\int_{\partial X}tr(\xi\wedge\star x_1).$$ 

As $\xi\in \LamP H^{2,p}_0(\Omega^0(\gpp))$ and $x_1\in y^{\lambda+\frac{1}{p}} H^{1,p}_0(\Omega^1(\gpp)\times\Omega^1(\gpp))$, 
we have $x_1|_{\partial X}=0$ and $\xi|_{\partial X}=0$. Therefore, $$\lan \DZ\xi,x\ran =\lan \xi,\DZS x\ran .$$

Suppose $x\in \Coker\;\DZ$, then for $\forall \xi\in \LaPP H^{2,p}_0$, we obtain $\lan \DZ\xi,x\ran =0$. As $\lambda> -1$, integrating by parts, we have $\lan \xi,\DZS x\ran =0$. Thus $\DZS x=0.$ Combining this with Proposition \ref{DZclose}, we finish the proof.

\qed

Fixe a reference connection pair $(A_0,\Phi_0)\in \FCS$ and $\epsilon> 0$. We set:
\begin{equation}
T_{(A,\Phi),\epsilon}: =\{(a,b)\in\Omega^1(\gpp)\times\Omega^1(\gpp)\;|\;\DZZS(a,b)=0,\|(a,b)\|_{\LamP H^{1,p}_0}< \epsilon\}.
\end{equation}

Thus we have a natural map $p:T_{(A,\Phi),\epsilon}\rightarrow \FBS$, which is induced by the inclusion of $T_{(A,\Phi),\epsilon}$ into $\FCS$ composed with quotienting by the gauge group $\FGGG$.

We have the following slicing theorem for the moduli space $\FBS$. 

\begin{theorem}
Given a point $(A,\Phi)\in \FCS$, denote by $[(A,\Phi)]\in \FBS$ the equivalence class under the projection map. For small $\epsilon> 0$, 

(1) if $(A,\Phi)$ is irreducible, then $T_{(A,\Phi),\epsilon}$ is a homeomorphism to a neighborhood of $[(A,\Phi)]$ in $\FBS$.

(2) if $(A,\Phi)$ is reducible, then $T_{(A,\Phi),\epsilon}/\Gamma_{(A,\Phi)}$ is a homeomorphism to a neighborhood of $[(A,\Phi)]$ in $\FBS$.
\end{theorem}
\proof Consider the map
\begin{equation}
\begin{split}
    &S:T_{(A,\Phi),\epsilon}\times \GP/\{\pm 1\}\rightarrow \CP,\\
    &S(A+a,\Phi+b,g)=g(A+a,\Phi+b).
\end{split}
\end{equation}
The map has derivative at $a=0,b=0,g=1$ as:
\begin{equation}
\begin{split}
    DS:&\Ker\; \DZSZ \times \Omega^0(\gpp)\rightarrow \Omega^1(\gpp)\times\Omega^1(\gpp),\\
    &(a,b,\xi)\rightarrow (a,b)+\DZZ(\xi).
\end{split}
\label{slicederivative}
\end{equation}
By Corollary \ref{DZsplit}, we know $DS$ is always surjective.

(1) If $(A,\Phi)$ is irreducible, then $DS$ is injective, by the implicit function theorem, we know for $\epsilon$ small enough, $S$ is a homeomorphism.

(2) If $(A,\Phi)$ is reducible, then $DS$ has kernel $H_{(A,\Phi)}^0.$ Let $H_{(A,\Phi)}^{0\perp}$ be the orthogonal of $H_{(A,\Phi)}^0$ with respect to the $L^2$ inner product. Then this time the restriction map $$S:T_{(A,\Phi),\epsilon}\times \exp(H_{(A,\Phi)}^{0\perp})/\{\pm 1\}\rightarrow \CP$$
is a local diffeomorphism. In addition, the multiplication map $\Gamma_{(A,\Phi)}\times \exp(H_{(A,\Phi)}^{0\perp})\rightarrow \GP$ at the identity will have derivative $1$. Thus for $g\in\GP$ close to $1$, there exist $l\in\Gamma_{(A,\Phi)}$ and $m\in \exp(H_{(A,\Phi)}^{0\perp})$ such that $g=ml$ and the splitting is unique. Therefore, we get a homeomorphism from $T_{(A,\Phi),\epsilon}/\Gamma_{(A,\Phi)}$ to a neighborhood of $[(A,\Phi)]$ in $\FBS$.

\qed
\end{subsection}

\begin{subsection}{Kuranishi Model}
Given $(A_0,\Phi_0)$ a solution to the Kapustin-Witten equations, all the other solutions within the slice $T_{(A,\Phi),\epsilon}$ are given by the set $Z(\Psi)$ of zeros of the map
\begin{equation}
\begin{split}
    &T_{(A,\Phi),\epsilon}\xrightarrow{\Psi} \LamM L^p(\Omega^2(\gpp)\times\Omega^0(\gpp)),\\
    &\Psi(a,b)=KW(A_0+a,\Phi_0+b)=\ML^1(a,b)+\{(a,b),(a,b)\}.
\end{split}
\end{equation}
By Theorem \ref{FredholmOperator}, $D\Psi$ is a Fredholm operator and for the homology associated to the Kuranishi complex, we have the following identification:
\begin{equation}
\begin{split}
&H^2_{(A,\Phi)}\cong \rm{Ker} \mathcal{L}^{\star}\cap \LamM L^p,\\
&H^1_{(A,\Phi)}\cong \rm{Ker} \mathcal{L} \cap (\rm{Ker} \DZS \cap y^{\lambda+\frac{1}{p}} H^{1,p}_0).  
\end{split}
\end{equation}

Therefore, we have the following Kuranishi picture of the moduli space:
\begin{proposition}{\rm{\cite{donaldson1983application} Proposition 8}}
For any solution $(A,\Phi)$ with $KW(A,\Phi)=0$, for $\epsilon$ sufficiently small, there is a map $\rho$ from a neighborhood of the origin in the harmonic space $H^1_{(A,\Phi)}$ to the harmonic space $H^2_{(A,\Phi)}$ such that if $(A,\Phi)$ is irreducible, a neighborhood of $[(A,\Phi)]\in \FMS$ is carried by a diffeomorphism onto $$Z(\rho)=\rho^{-1}(0)\subset H^1_{(A,\Phi)}$$ 
and if $(A,\Phi)$ is reducible, then a neighborhood of $[(A,\Phi)]$ is modelled on
$$Z(\rho)/\GAS.$$
\label{Kuranishimodel}
\end{proposition}

\end{subsection}

\end{section}

\begin{section}{Exponential Decay}
In this section, we will prove the exponential decay over the cylindrical ends which is identified with $Y^3\times (0,+\infty)$ with the convergence assumption. We denote $y$ as the coordinate of $(0,+\infty)$. As before, let $\Vol_Y$ be a given volume form of $Y$. We denote $\star_4$ by the 4-dimensional Hodge star operator with respect to the Volume form $\Vol_Y\wedge dy$ and denote $\star$ by the 3 dimension Hodge star operator with respect to $\Vol_Y$. 

\begin{theorem}
Let $(A,\Phi)$ be a solution to the Kapustin-Witten equations over $Y^3\times (0,+\infty)$, and $Y_T:=Y^3\times \{T\}\in Y^3\times (0,+\infty)$ is a slice. For a non-degenerate flat $\slc$ connection $(A_{\rho},\phi_{\rho})$ corresponding to the representation $\rho$. 

Suppose for some $p>2$,  $\lim_{T\rightarrow+\infty}\|(A,\Phi)-(A_{\rho},\phi_{\rho})\|_{L^p_1(Y_T)}=0$, then there exists a positive number $\delta$, such that $\|(A,\Phi)-(A_{\rho},\phi_{\rho})\|_{C^{\infty}(Y\times [T,+\infty))}\leq Ce^{-\delta T}.$
\label{exponentialdecay}
\end{theorem}

We follow the ideas in \cite{donaldson2002floer}.

Over slice $Y_t$, denote $(A,\Phi)\vert_{Y_t}=(A(t),\phi(t)+\phi_y(t)dy)$, denote $\gamma(t):=(A(t),\phi(t),\phi_y(t))$. Recall the gradient of the extended Chern-Simons function is denoted as $\nabla\ECS$, then the flow equations (\ref{BigGradientFlow}) can be rewrote as 
\begin{equation}
\frac{d}{dy}\gamma(t)+\nabla\ECS(\gamma(t))=0.
\label{600}
\end{equation}

\begin{definition}
The analytic energy is defined as $\Eng:=\int_{Y_t}|\frac{d}{dy}\gamma(t)|^2+|\nabla\ECS(\gamma(t))|^2 d\Vol.$
\end{definition}

Now we introduce some basic computation related to the analytic energy $\Eng$. Let $\ECS(\rho)$ to be $\ECS(A_\rho,\phi_\rho,0)$ and $\ECS(T)$ to be $\ECS(\gamma(T))$, then we have the following proposition:

\begin{proposition}
Over $Y\times (0,+\infty)$, under the assumption of Theorem \ref{exponentialdecay}, denote $(A,\Phi=\phi+\phi_ydy)$ and denote $J(T):=\int^{+\infty}_{T}\Eng dt$, then $J(T)=2(\ECS(T)-\ECS(\rho))$.
\label{JT}
\end{proposition}
\proof

We have the following computations for the gradient flow equations (\ref{600}):
\begin{equation}
\begin{split}
    0=&\int_{Y\times [T,+\infty)}|\frac{d}{dy}\gamma(t)+\nabla\ECS(\gamma(t))|^2\\
    =&\int_{Y\times [T,+\infty)}|\frac{d}{dy}\gamma(t)|^2+|\nabla\ECS(\gamma(t))|^2+2\lan\frac{d}{dy}\gamma(t),\nabla\ECS(\gamma(t))\ran.
\end{split}
\end{equation}

By definition of the gradient, we obtain
\begin{equation}
\int_{Y\times [T,+\infty)}2\lan\frac{d}{dy}\gamma(t),\nabla\ECS(\gamma(t))\ran=2(\ECS(\rho)-\ECS(T)).
\end{equation}

Therefore, $J(T)=2(\ECS(T)-\ECS(\rho)).$
\qed

Let $(A_t,\Phi_t=\phi_t+(\phi_y)_tdy)$ be the restriction of the solution $(A,\Phi)$ to the slice $Y_t$. Take $(a_t,b_t,c_t)=(A_t,\Phi_t,(\phi_y)_t)-(A_\rho,\phi_\rho,0)$, denote $\MH_\rho(a_t,b_t,c_t):=\MH_{(A_\rho,\phi_\rho,0)}(a_t,b_t,0,c_t)$ and $\MEH_\rho:=\MEH_{(A_\rho,\phi_\rho,0)}$. Then we have the following lemma:
\begin{lemma}\label{Guudlemma}
Under the assumption of Theorem \ref{exponentialdecay}, then for $t$ large enough, there exist positive constant $C_1$, $C_2$ such that we have the following estimates:
\begin{equation}
\Eng\geq C_1\|\MEHR(a_t,b_t,c_t)\|^2_{L^2(Y_t)},
\end{equation}
\begin{equation}
\ECS(A_\rho+a_t,\phi_\rho+b_t,c_t)-\ECS(A_\rho,\phi_\rho,0)\leq C_2\|\MEHR(a_t,b_t,c_t)\|^2_{L^2(Y_t)},
\end{equation}
\end{lemma}
where $C_2$ is a number depends on the smallest absolute eigenvalue of $\MEHR$.
\proof
As we assume $\rho$ is nondegenerate, by Proposition \ref{kernelhess}, we have the following estimate $\|(a_t,b_t,c_t)\|_{L^2_1(Y_t)}\leq C\|\MEHR(a_t,b_t,c_t)\|_{L^2(Y_t)}$ over the slice $Y_t$. 

In addition, by the convergence assumption, $(a_t,b_t,c_t)$ will have small $L^p_1$ norm when $t$ is large enough, by Proposition $\ref{KWG}$, in each slice $ Y_t$, we can make $(A,\Phi)$ in the Kapustin-Witten gauge relative to $(A_{\rho},\Phi_{\rho})$: $$\LGF_{(A_{\rho},\Phi_{\rho})}(a_t,b_t)=0.$$ To be explicit, we have \begin{equation}
d_{A_{\rho}}^{\star}a_t-\star[\Phi_\rho,\star b_t]=0.
\label{ConditionB}
\end{equation}

Under this gauge, we have 
\begin{equation}
\MH_\rho(a_t,b_t,c_t)=\MEHR(a_t,b_t,c_t).
\end{equation}

Now, we have enough preparation for proving the estimate.

\begin{equation}
\begin{split}
\Eng&\geq \int_{Y_t}|\nabla\ECS(A_t,\phi_t,(\phi_y)_t)|^2 \\
&\geq \int_{Y_t}|\MHR(a_t,b_t,c_t)+\{(a_t,b_t,c_t)^2\}|^2\\
&\geq \int_{Y_t}|\MHR(a_t,b_t,c_t)|^2+|\{(a_t,b_t,c_t)^2\}|^2+2\lan\MHR(a_t,b_t,c_t),\{(a_t,b_t,c_t)^2\} \ran.
\end{split}
\label{6969}
\end{equation}

Over 3 dimensional manifold, we have the Sobolev embedding $L^2_1(Y_t)\rightarrow L^r(Y_t)$ for $r\leq 6$.

Then we have $$\|\{(a_t,b_t,c_t)^2\}^2\|_{L^2(Y_t)}^2\leq C\|(a_t,b_t,c_t)\|_{L^4(Y_t)}^4 \leq C\|\MEHR(a_t,b_t,c_t)\|_{L^2(Y_t)}^4,$$
$$\int_{Y_t}\lan\MHR(a_t,b_t,c_t),\{(a_t,b_t,c_t)^2\} \ran\leq C\|\MEHR(a_t,b_t,c_t)\|^3_{L^2(Y_t)}.$$

By the convergence assumption, we know $\lim_{t\rightarrow+\infty}\|\MEHR(a_t,b_t,c_t)\|_{L^2(Y_t)}=0$. Thus, for inequality (\ref{6969}), $\|\{(a_t,b_t,c_t)^2\}^2\|_{L^2(Y_t)}^2$ and $\int_{Y_t}\lan\MHR(a_t,b_t,c_t),\{(a_t,b_t,c_t)^2\} \ran$ can be absorb by the first term and by choosing $t$ large enough, we get the estimate we want.

For the statement (2), under the gauge fixing condition (\ref{ConditionB}), we have the following estimate
\begin{equation}
\begin{split}
&\int_{Y_t}\lan (a_t,b_t,c_t),\MHR(a_t,b_t,c_t) \ran\\
=&\int_{Y_t}\lan (a_t,b_t,c_t),\MEHR(a_t,b_t,c_t) \ran\\
\leq&\frac{1}{\delta}\|\MEHR(a_t,b_t,c_t)\|_{L^2(Y_t)}^2,
\end{split}
\end{equation}
where $\delta$ is smallest absolute eigenvalue of the operator $\MEHR$ and by the non-degenerate assumption, $\Ker\MEHR=0$ and $\delta$ is bounded blow away from $0$. 

Thus, we have
\begin{equation}
\begin{split}
    &\ECS(A_\rho+a_t,\phi_\rho+b_t,c_t)-\ECS(A_\rho,\phi_\rho,0)\\
    \leq& -\frac{1}{2}\int_{Y_t}\lan (a_t,b_t,c_t),\MHR(a_t,b_t,c_t) \ran-\int_{Y_t}\{(a,b,c)^3\}\\
    \leq& C\|\MEHR(a_t,b_t,c_t)\|_{L^2(Y_t)}^2+C\|\MEHR(a_t,b_t,c_t)\|^3_{L^2(Y_t)}\\
    \leq& C_2\|\MEHR(a_t,b_t,c_t)\|_{L^2(Y_t)}^2.
\end{split}
\end{equation}
\qed

Now we obtain the following proposition:
\begin{proposition}
With the assumption above, if $\Eng$ is bounded, then there exists a constant C such that 
$$J(T)\leq Ce^{-\delta t}.$$
\end{proposition}
\proof By Lemma \ref{Guudlemma}, we have the following:
\begin{equation}
\begin{split}
    J(t)&=\ECS(A_\rho+a_t,\phi_\rho+b_t,c_t)-\ECS(A_\rho,\phi_\rho,0)\\
    &\leq C_2\|\MEHR(a_t,b_t,c_t)\|^2_{L^2(Y_t)}\\
    &\leq \frac{C_2}{C_1}\Eng\\
    &\leq -\frac{C_2}{C_1}\frac{d}{dt}J(t).
\end{split}
\end{equation}
Thus, take $\delta=\frac{C_1}{C_2}$, we have:
$$\delta J(t)+\frac{dJ(t)}{dt}\leq 0.$$
From here we get that $J(t)\leq Ce^{-\delta t}.$
\qed

Using these corollaries, we can give the following estimate of the decay of solutions.
\begin{proposition}
For all $T$ is large enough that we have $\|(A_t,\Phi_t)-(A_{\rho},\Phi_{\rho})\|_{L^2_k(Y\times [T,+\infty))}^2\leq Ce^{-\delta T}$.
\end{proposition}
\proof
Fixing a Kapustin-Witten gauge for $(a_t,b_t,c_t):=(A_t,\phi_t,(\phi_y)_t)-(A_{\rho},\phi_{\rho},0)$. By the non-degenerate assumption, we have $\|(a_t,b_t,c_t)\|_{L^2_1(Y_t)}\leq C\|\MEHR(a_t,b_t,c_t)\|_{L^2(Y_t)}$.

In addition, by Lemma \ref{Guudlemma}, for $T$ is large enough, we have  $$\|\MEHR(a_t,b_t)\|^2_{L^2(Y_t)}\leq C\Eng.$$ 
Therefore, we compute
\begin{equation}
\begin{split}
    \|(A,\Phi)-(A_{\rho},\Phi_{\rho})\|_{L^2_1(Y\times [T,+\infty))}    = &\int^{+\infty}_T\|(a_t,b_t)\|_{L^2_1(Y_t)}dt\\
    \leq &C\int^{+\infty}_T \|\MEHR(a_t,b_t)\|_{\LTY}dt\\\leq&  C\int^{+\infty}_T\Eng dt\\    \leq& C J(t).
\end{split}
\end{equation}

By the exponential decay of $J(t)$, we proved the result for k=1. Take bootstrapping method, we get that the $L^2_k$ norm exponentially decays.
\qed
$$$$
\textit{Proof of Theorem \ref{exponentialdecay}}: 
We only need to show that for every integer $k$, we have $\|(a_t,b_t)\|_{\mathcal{C}^k}\leq Ce^{-\delta t}$. By the Sobolev embedding for $L^2_{k'}$ and $C^k$, the result follows immediately.
\qed
\end{section}

\begin{section}{Constructing Solutions}
In this section, we will prove the gluing theorem for the Kapustin-Witten equations with Nahm pole boundary condition.

For $i=1,2$, consider $X_i$ to be a 4-manifold with boundary $Z_i$ and infinite cylindrical end identified with $Y_i\times (0,+\infty)\subset X_i$, let $P_i$ to denote a $SU(2)$ bundle and $(A_i,\Phi_i)\in \mathcal{C}_{P_i}$ be a solutions to the Kapustin-Witten equations over bundle $P_i$ which approach to a flat connection $\rho_i$ and satisfies the Nahm pole boundary condition on the boundary $Z_i$.

If $Y_1=Y_2$, we can define a new family of 4-manifolds $\XT$. To be precise, we fix an isometry between $Y_1$ and $Y_2$, we first delete the infinite portions $Y_1\times [2T,+\infty)$ $Y_2\times$, $[2T,+\infty)$ from the two ends, and then identify $(y,t)\in Y_1\times (T,2T)$ with $(y,2T-t)\in Y_2\times (T,2T)$. This is in Figure \ref{fig:text4463} and Figure \ref{fig:text4621}:

\begin{figure}[H]
    \centering
    \includegraphics[width=0.8\textwidth]{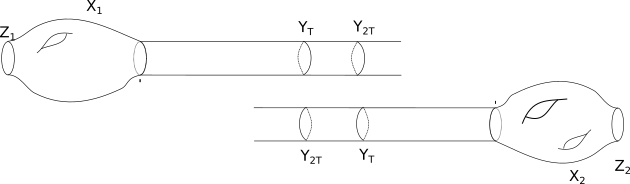}
    \caption{Two cylindrical-end manifold $X_1$ and $X_2$ with boundary}
    \label{fig:text4463}
\end{figure}

\begin{figure}[H]
    \centering
    \includegraphics[width=0.8\textwidth]{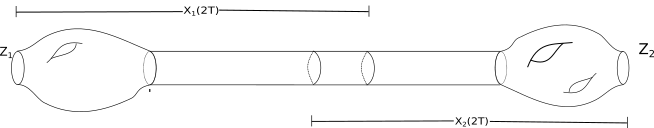}
    \caption{$X_1$, $X_2$ glued together to form $\XT$}
    \label{fig:text4621}
\end{figure}

In addition, if the limit flat connections coincide, we denote  $\rho=\rho_1=\rho_2$, we can fix an identification of these flat bundles and get a new bundle $\PT$ with a natural connection $(A^{\sharp},\Phi^{\sharp})$, which we will explicitly define in subsection 7.1.

Now we restate our theorem as follows:
\begin{theorem}
Under the hypotheses above, if

(a) $\lim_{T\rightarrow +\infty}\|(A_i,\Phi_i)-(A_{\rho_i},\Phi_{\rho_i})\|_{L^{p_0}_1(Y_i\times\{T\})}=0$ for some $p_0>2$,

(b) $\rho$ is an acyclic $\slc$ flat connection,

then for $p\geq 2$ and $\lambda\in[1-\frac{1}{p},1)$, we have:

(1) for some constant $\delta$, there exists a $y^{\lambda+\frac{1}{p}}H^{1,p}_0$ pair $(a,b)\in \Omega^1_{\XT}(\gpp)\times\Omega^1_{\XT}(\gpp)$ with $$\|(a,b)\|_{y^{\lambda+\frac{1}{p}-1}L^p_1}\leq Ce^{-\delta T},$$

(2) there exists an obstruction class $h\in H^2_{(A_1,\Phi_1)}(X_1)\times H^2_{(A_2,\Phi_2)}(X_2)$ such that $h=0$ if and only if $(A^{\sharp}+a,\Phi^{\sharp}+b)$ is a solution to the Kapustin-Witten equations (\ref{KW}).
\label{gluingobstrucion}
\end{theorem}

We break the proof of this theorem into several parts. 
\begin{subsection}{Approximate Solutions}
Denote by $(a_i,\phi_i):=(A_i,\Phi_i)-(A_{\rho_i},\Phi_{\rho_i}),$ the difference between our solution and the limit flat connections. 

Define a new pair $(A_i',\Phi_i')=(A_{\rho_i},\Phi_{\rho_i})+\chi(t)(a_i,\phi_i)$, here $\chi(t)$ is a cut off function which equals $0$ on the complement of  $Y_i\times(T+1,+\infty)$ and $1$ on $X_i(T):=X_i\backslash Y_i\times (T,+\infty)$. From this construction, we know that $KW(A_i',\Phi_i')$ is supported on $Y\times (T,T+1)$.

As $(A_1',\Phi_1')$ and $(A_2',\Phi_2')$ agree on the end, we can glue then together to get an approximate solution $(A^{\sharp},\Phi^{\sharp})$ on $\XT$, and we denote the new bundle as $P^{\sharp}$. Using the above process, we can define a map: \begin{equation}
\begin{split}
    &I:\mathcal{C}_{P_1}\times \mathcal{C}_{P_2}\rightarrow \mathcal{C}_{P^{\sharp}}\\
    &I((A_1,\Phi_1),(A_2,\Phi_2)):=(A^{\sharp},\Phi^{\sharp}).
\end{split}
\label{approximatemap}
\end{equation} and this map depends on choice of $T$ and the cut-off function we choose.

We have the following estimate for the approximate solution:

\begin{proposition}
$\|KW(A^{\sharp},\Phi^{\sharp})\|_{\LamM L^p(\XT)}\leq Ce^{-\delta T}.$
\label{small}
\end{proposition}
\proof
By construction, we know $KW(A^{\sharp},\Phi^{\sharp})$ is only supported on a compact subset of $\XT$ and in this area the weight function is bounded. Combining this with the exponential decay result: Theorem $\ref{exponentialdecay}$, we get the estimate we want.
\qed

\end{subsection}
\begin{subsection}{Gluing Regular Points}
First, we assume that $H^2_{(A_1,\Phi_1)}=0$ and $H^2_{(A_2,\Phi_2)}=0$ and the limit flat connection is irreducible.

We use the previous notation from (\ref{KWO}). Recall that $\ML_{(A_i,\Phi_i)}$ is the linearization of the Kapustin-Witten equations, we  denote $\ML_i:=\ML_{(A_i,\Phi_i)}$. Recall $\LGF_i:=\LGF_{(A_i,\Phi_i)}$ is the Kapustin-Witten gauge fixing operator (\ref{KWGF1}). Now, we denote $\MD_i:=\ML_i\oplus \LGF_i$. By Theorem \ref{FredholmOperator}, we get a Fredholm operator over $X_i$:
$$\MD_i:\LamP H^{1,p}_0(X_i)\rightarrow \LamM L^p(X_i).$$

By assumption, we know $\MD_i$ is surjective, then there exists a right inverse
\begin{equation}
Q_i:\LamM L^p(X_i)\rightarrow \LamP H^{1,p}_0(X_i),
\label{Qi}
\end{equation} such that $\MD_iQ_i=Id$. Therefore, after restricting the domain of $Q_i$ to the image of $\ML_i$, we get a right inverse for $\ML_i$ and for simplicity, we still denote the right inverse as $Q_i$ and we obtain $\ML_iQ_i=Id.$ 

Take $\phi_i$ to be a cut off function supported in $X_i(2T)$, with $\phi_i(x)=1$ on $X_i(T)$ and $\phi_1+\phi_2=1$ on $\XT$. The graph of cut-off function $\phi_1$ is in the following Figure \ref{fig:text4871}:

\begin{figure}[H]
    \centering
    \includegraphics[width=0.8\textwidth]{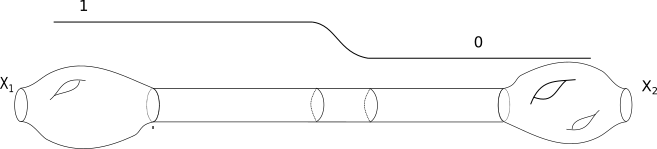}
    \caption{The graph of cut-off function $\phi_1$}
    \label{fig:text4871}
\end{figure}

By definition, we chose $\phi_i$ with the estimate $\|\nabla\phi_i\|_{L^{\infty}(X_i)}\leq \epsilon(T)\leq \frac{C}{T}$.

Take $\xi\in \LamM L^p(\XT)$, denote by $\xi_i$ the restriction of $\xi$ on $X_i(2T)$, then we define a new approximate inverse operator $\hat{Q}(\xi):=\phi_1Q_1(\xi_1)+\phi_2Q_2(\xi_2)$, which can be written as \begin{equation}
\hat{Q}:\LamM L^p(X_i)\rightarrow \LamP H^{1,p}_0(X_i).
\label{hatQ}
\end{equation}

Denoting $\LS:=\mathcal{L}_{(A^{\sharp},\Phi^{\sharp})}$ as follows:

$$\LS:\LamP H^{1,p}_0(X_i)\rightarrow \LamM L^p(X_i).$$
After these preparations, we have the following relationship between these two operators:
\begin{lemma}
\label{inverse1}
For $\LS$, $\hat{Q}$ as above, and for $\forall \xi\in \LamM L^p(\XT)$, we have 
$$\|\LS\hat{Q}(\xi)-\xi\|_{\LamM L^p(\XT)}\leq \epsilon(T)\|\xi\|_{\LamM L^p(\XT)}.$$
\end{lemma}
\proof
For $\LS\hat{Q}(\xi)$, by definition, we have the following computation:
\begin{equation}
    \begin{split}
        \LS \hat{Q}(\xi)&=\LS(\phi_1 Q_1(\xi_1)+\phi_2 Q_2(\xi_2))\\
        =&\nabla\phi_1\star Q_1(\xi_1)+\nabla\phi_2\star Q_2(\xi_2)+\phi_1\LS Q_1(\xi_1)+\phi_2\LS Q_2(\xi_2)\\
        =&\nabla\phi_1\star Q_1(\xi_1)+\nabla\phi_2\star Q_2(\xi_2)+\phi_1\ML_1 Q_1(\xi_1)+\phi_2\ML_2 Q_2(\xi_2)\\
        &+\phi_1(\LS-\ML_1)Q_1(\xi_1)+\phi_2(\LS-\ML_2)Q_2(\xi_2).
    \end{split}
\end{equation}
For the term $\nabla\phi_1\star Q_1(\xi_1)+\nabla\phi_2\star Q_2(\xi_2)$, we know $\|\nabla\phi\|_{L^{\infty}}< \epsilon(T).$ 

Therefore, we obtain $$\|\nabla\phi_1\star Q_1(\xi_1)\|_{\LamM L^p(\XT)}\leq\epsilon(T)\|Q_1(\xi_1)\|_{\LamM L^p(\XT)}.$$

By Proposition \ref{spacelooklike}, $\LamP H^{1,p}_0(\XT)\subset \LamP L^p(\XT)$, and by Proposition \ref{embeddinglemma}, we obtain $\LamP L^p(\XT)\subset \LamM L^p(\XT)$, therefore, $\LamP H^{1,p}_0(\XT)\subset \LamM L^p(\XT)$. In addition, by (\ref{Qi}), we know $$\|Q_1(\xi_1)\|_{\LamP H^{1,p}_0(\XT)}\leq C\|\xi_1\|_{\LamM L^p(\XT)}.$$

Therefore, we obtain
\begin{equation}
\begin{split}
    \|\nabla\phi_1\star Q_1(\xi_1)\|_{\LamM L^p(\XT)}
    \leq &\epsilon(T)\|Q_1(\xi_1)\|_{\LamM L^p(\XT)}\\
    \leq &\epsilon(T)\|Q_1(\xi_1)\|_{\LamP H^{1,p}_0(\XT)}\\
    \leq &\epsilon(T)C\|\xi\|_{\LamM L^p(\XT)} (\text{Here the constant $C$ is independent of $T$}).
\end{split}
\end{equation}
Similarily, we have the same estimate for $Q_2$:
\begin{equation}
\|\nabla\phi_2\star Q_2(\xi_2)\|_{\LamM L^p(\XT)}\leq \epsilon(T)C\|\xi\|_{\LamM L^p(\XT)} (\text{Here the constant $C$ is independent of $T$}).
\end{equation}

For the term $\phi_1(\LS-\ML_1)Q_1(\xi_1)+\phi_2(\LS-\ML_2)Q_2(\xi_2)$, by Theorem $\ref{exponentialdecay}$, we know that the operators $\LS-\ML_1$ and  $\LS-\ML_2$ are order zero and the operator norm will exponentially decay as $T\rightarrow\infty$ . Therefore, we have
\begin{equation}
\|\phi_1(\LS-\ML_1)Q_1(\xi_1)+\phi_2(\LS-\ML_2)Q_2(\xi_2)\|_{\LamM L^p(\XT)}< \epsilon(T)\|\xi\|_{\LamM L^p(\XT)}.
\end{equation} 

For the remaining terms, we have
\begin{equation}
    \begin{split}
        &\phi_1\ML_1Q_1(\xi_1)+\phi_2\ML_2 Q_2(\xi_2)=\phi_1\xi_1+\phi_2\xi_2=\xi.
    \end{split}
\end{equation}

Combining all the discussion above, we get the estimate we want.

\qed

\begin{proposition}
\label{inverseoperator}
There exists an operator $\QS$ with $$\QS:\LamM L^p(\XT)\rightarrow \LamP H^{1,p}_0(\XT),$$ such that $\LS \QS=Id$. In addition, there exists a constant $C$ independent of $T$ such that for $\forall \xi\in \LamM L^p(X_i), $ we have
$$\|\QS(\xi)\|_{\LamP H^{1,p}_0(\XT)}\leq C\|\xi\|_{\LamM L^p(\XT)}.$$
\end{proposition}

\proof
Take $R(\xi):=\LS \hat{Q}(\xi)-\xi$, by Proposition \ref{inverse1}, we know when $T$ is large enough, the operator norm of $R$ will be very small. Therefore, $R+Id$ is invertible.

Take $\QS:=\hat{Q}(Id+R)^{-1}$ then by definition, we have that $\QS$ is an operator from $\LamM L^p(\XT)$ to $\LamP H^{1,p}_0(\XT)$ and $\LS\QS=Id$. 

 The operator norm of $R+Id$ is less than three and the operator norm of $\hat{Q}$ is dominated by the operator norm of $Q_1$ plus the operator norm of $Q_2$. Thus, the operator norm of $\QS$ is independent of $T$.

\qed

Given an approximate solution $(A^{\sharp},\Phi^{\sharp})$, for any connection $(A,\Phi)$, write $$(A,\Phi)=(A^{\sharp},\Phi^{\sharp})+(a,b),$$ we hope to find suitable $(a,b)$ such that $KW(A,\Phi)=0$.

By Proposition \ref{structure}, we have the following quadratic expansion: $$KW(A,\Phi)=KW(A^{\sharp},\Phi^{\sharp})+\LS(a,b)+\{(a,b),(a,b)\}.$$
We will solve the equations 
\begin{equation}
KW(A^{\sharp},\Phi^{\sharp})+\LS(a,b)+\{(a,b),(a,b)\}=0.
\label{gluingregularity}
\end{equation}

Take $\eta:=-KW(A^{\sharp},\Phi^{\sharp})$ and replace $(a,b)$ by $\QS(\alpha)$, then the quadratic expansion becomes
\begin{equation}
    \eta=\alpha+\{\QS(\alpha),\QS(\alpha)\}.
    \label{quadraticequaion}
\end{equation} 

Now our target is to solve this equation for some $\alpha\in \LamM L^p(\XT)$.

Take $S(\alpha):=\{\QS(\alpha),\QS(\alpha)\}$, we have the following proposition for the operator $S$:

\begin{proposition}
\label{SS}
For any $\lambda_0\in[1-\frac{1}{p},1)$, S is an operator: $$S:\LamM L^p(\XT)\rightarrow \LamM L^p(\XT)$$ satisfying
$S(0)=0$ and for two elements $\alpha,\beta\in L^p(\XT)$, there exists a constant $k$ independent of $T$ such that $$\|S(\alpha)-S(\beta)\|_{\LamM L^p}\leq k(\| \alpha \|_{\LamM L^p}+\|\beta\|_{\LamM L^p})(\|\alpha-\beta\|_{\LamM L^p}).$$\end{proposition}
\proof First, we prove $S$ is the suitable operator. As $\QS:\LamM L^p(\XT)\rightarrow \LamP H^{1,p}_0(\XT)$, using the Sobolev embedding $\LamP H^{1,p}_0(\XT)\hookrightarrow \LamM L^p_1(\XT)$(Corollary \ref{Sobolevembedding}), we can consider $\QS$ as an operator from $\LamM L^p(\XT)$ to $\LamM L^p_1(\XT)$. 

Denote $\alpha=(\alpha_1,\alpha_2)$, from the definition of $S$, we have
\begin{equation}
\begin{split}
    S(\alpha)&=\{(\alpha_1,\alpha_2),(\alpha_1,\alpha_2)\}\\
    &=\left(
    \begin{array}{l}
        \QS(\alpha_1)\wedge \QS(\alpha_1))-\QS(\alpha_2)\wedge \QS(\alpha_2)+\star[\QS(\alpha_1),\QS(\alpha_2)] \\
        -\star[\QS(\alpha_1),\star\QS(\alpha_2)]
    \end{array}
    \right).
\end{split}
\end{equation}

As the terms appearing in $S(\alpha)$ are quadratic terms, using the H\"{o}lder inequality for $\lambda\geq 1-\frac{1}{p}$ such that $\|fg\|_{\LamM L^p(\XT)}\leq \|f\|_{\LamM L^p_1(\XT)}\|g\|_{\LamM L^p_1(\XT)}$(Corollary \ref{Sobolevembedding}), we have that $S$ is an operator $$S:\LamP L^p(\XT)\rightarrow \LamP L^p(\XT).$$ 

Now, we will show that $S$ has the desired estimate. Denote $\beta=(\beta_1,\beta_2)$, from the definition of $S$, we have the following computation:
\begin{equation}
\begin{split}
    &S(\alpha)-S(\beta)\\
    =&\{(\alpha_1,\alpha_2),(\alpha_1,\alpha_2)\}-\{(\beta_1,\beta_2),(\beta_1,\beta_2)\}\\
    =&\left(
    \begin{array}{l}
        \QS(\alpha_1)\wedge \QS(\alpha_1)-\QS(\alpha_2)\wedge \QS(\alpha_2)+\star[\QS(\alpha_1),\QS(\alpha_2)]-\star[\QS(\beta_1),\QS(\beta_2)] \\
        -\star[\QS(\alpha_1),\star\QS\alpha_2]+\star[\QS(\beta_1),\star\QS\beta_2]
    \end{array}
    \right).
\end{split}
\end{equation}
Now we make estimates for each term appearing in $S(\alpha)-S(\beta)$. 
\begin{equation}
\begin{split}
    &\|\QS(\alpha_1)\wedge \QS(\alpha_1)-\QS(\beta_1)\wedge \QS(\beta_1)\|_{\LamM L^p(\XT)}\\
    \leq&\|\QS(\alpha_1)\wedge (\QS(\alpha_1)-\QS(\beta_1))\|_{\LamM L^p(\XT)}+\|(\QS(\alpha_1)-\QS(\beta_1))\wedge \QS(\beta_1)\|_{\LamM L^p(\XT)}\\
    \leq&\|\QS(\alpha_1)\|_{\LamM L^p_1(\XT)}\|\QS(\alpha_1)-\QS(\beta_1)\|_{_{\LamM L^p_1(\XT)}}\\
    &+\|\QS(\beta_1)\|_{\LamM L^p_1(\XT)}\|\QS(\alpha_1)-\QS(\beta_1)\|_{_{\LamM L^p_1(\XT)}}\\
    \leq&C\|\alpha_1\|_{\LamM L^p(\XT)}\|\alpha_1-\beta_1\|_{\LamM L^p(\XT)}+C\|\beta_1\|_{\LamM L^p(\XT)}\|\alpha_1-\beta_1\|_{\LamM L^p(\XT)}\\
    \leq&C(\|\alpha_1\|_{\LamM L^p(\XT)}+\|\beta_1\|_{\LamM L^p(\XT)})(\|\alpha_1-\beta_1\|_{\LamM L^p(\XT)}).
\end{split}
\end{equation}
For another term, we have 
\begin{equation}
\begin{split}
    &\|\star[\QS(\alpha_1),\QS(\alpha_2)]-\star[\QS(\beta_1),\QS(\beta_2)]\|_{\LamM L^p(\XT)}\\
    \leq&\|[\QS(\alpha_1),\QS(\alpha_2)]-[\QS(\alpha_1),\QS(\beta_2)]+[\QS(\alpha_1),\QS(\beta_2)]-[\QS(\beta_1),\QS(\beta_2)]\|_{\LamM L^p(\XT)}\\
    \leq&\|[\QS(\alpha_1),\QS(\alpha_2-\beta_2)]\|_{\LamM L^p(\XT)}+\|[\QS(\alpha_1-\beta_1),\QS(\beta_2)]\|_{\LamM L^p(\XT)}\\
    \leq&C(\| \alpha_1 \|_{\LamM L^p(\XT)}\|\alpha_2-\beta_2\|_{\LamM L^p(\XT)}+\|\beta_2\|_{\LamM L^p(\XT)} \|\alpha_1-\beta_1\|_{\LamM L^p(\XT)})\\
    \leq& C(\|\alpha\|_{\LamM L^p(\XT)}+\|\beta\|_{\LamM L^p(\XT)})(\|\alpha-\beta\|_{\LamM L^p(\XT)}).
\end{split}
\end{equation}
Similarily, we have the following estimate:
\begin{equation}
\begin{split}
    &\|-\star[\QS(\alpha_1),\star\QS(\alpha_2)]+\star[\QS(\beta_1),\star\QS(\beta_2)]\|_{\LamM L^p(\XT)}\\
    \leq&C(\|\alpha\|_{\LamM L^p(\XT)}+\|\beta\|_{\LamM L^p(\XT)})(\|\alpha-\beta\|_{\LamM L^p(\XT)}).
\end{split}
\end{equation}
Combining the previous computations, we have the result we want.
\qed

We have the following lemma about the operator $S$:
\begin{lemma}{\rm{(\cite{donaldson1990geometry} Lemma $7.2.23$})}\label{LS}
Let $B$ be a Banach space and let $\|\,\|_B$ be the norm on B. Let $S:B\rightarrow B$ be a smooth map on the Banach space B with $S(0)=0$ and $\|S\xi_1-S\xi_2\|_B\leq k(\|\xi_1\|_B+\|\xi_2\|_B)(\|\xi_1-\xi_2\|_B)$, for some constant $k> 0$ and all $\xi_1,\xi_2$ in B,then for each $\eta\in B$ with $\|\eta\|_B< \frac{1}{10k}$, there exists a unique $\xi$ with $\|\xi\|_B\leq \frac{1}{5k}$ such that $$\xi+S(\xi)=\eta.$$
\end{lemma}
We now can complete the proof of Theorem 1.1.

\leftline{\textit{Proof of Theorem 1.1}: 
Recall we hope to solve the equation ($\ref{quadraticequaion}$), which is}
$$\eta=\alpha+S(\alpha).$$
By Proposition \ref{SS}, in Lemma \ref{LS}, if we take the Banach space $B$ as $\LamM L^p(\XT)$, we know that the operator $S$ satisfies the assumption in Lemma \ref{LS}. Therefore, there exists an solution $\alpha$ to equation ($\ref{quadraticequaion}$) with $\alpha\in \LamM L^p(\XT)$. 

Let $(a,b):=\QS(\alpha)$ where $(a,b)\in \Omega^1\times\Omega^1$, then $(A^{\sharp},\Phi^{\sharp})+(a,b)$ is a solution to the Kapustin-Witten equations ($\ref{KW}$). 

Now we will prove the regularity statement of Theorem 1.1. By Proposition \ref{inverseoperator}, we have 
$$\QS:\LamM L^p(\XT)\rightarrow \LamP H^{1,p}_0(\XT).$$ Therefore,  we know $(a,b)\in \LamP H^{1,p}_0(\XT)$.

As $S$ satisfies $$\|S(\alpha)-S(\beta)\|_{\LamM L^p(\XT)}\leq k(\| \alpha \|_{\LamM L^p(\XT)}+\|\beta\|_{\LamM L^p(\XT)})(\|\alpha-\beta\|_{\LamM L^p(\XT)}).$$
Take $\beta=0$, we have $\|S(\alpha)\|_{\LamM L^p(\XT)}\leq k\|\alpha\|_{\LamM L^p(\XT)}^2.$ By equation ($\ref{quadraticequaion}$), we have the following estimate:
\begin{equation}
\begin{split}
    \|\alpha\|_{\LamM L^p(\XT)}\leq &\|\eta\|_{\LamM L^p(\XT)}+\|S(\alpha)\|_{\LamM L^p(\XT)}\\
    \leq &\|\eta\|_{\LamM L^p(\XT)}+k\|\alpha\|^2_{\LamM L^p(\XT)}.
\end{split}
\end{equation}

WLOG, we can assume $1-k\|\alpha\|_{\LamM L^p(\XT)}\geq \frac{1}{2}$ and we obtain 
\begin{equation}
\|\alpha\|_{\LamM L^p(\XT)}\leq 2\|\eta\|_{\LamM L^p(\XT)}.
\label{decayestimate2}
\end{equation}
 As in Proposition \ref{SS}, we use the estimate $\|\QS(\alpha)\|_{\LamM L^p_1(\XT)}\leq C\|\alpha\|_{\LamM L^p(\XT)}$, we have \begin{equation}
\|(a,b)\|_{\LamM L^p_1(\XT)}\leq C\|\alpha\|_{\LamM L^p(\XT)}\leq 2C\|\eta\|_{\LamM L^p(\XT)}.
\label{decayestimate1}
\end{equation}

Applying Proposition \ref{small}, we get the estimate we want.\qed

We can say more about the regularity of solutions we get. Using the equations (\ref{gluingregularity}), we have the following proposition.

\begin{proposition}
For $p> 2$ and $T$ is large enough, suppose $(a,b)$ satisfies the equations (\ref{gluingregularity}) over $\XT$, then $(a,b)$ is smooth in the interior of $\XT$.
\label{regularityofgluing}
\end{proposition}
\proof
Fix a interior open set $U\subset \XT$. By (\ref{decayestimate1}), for any given constant $C$, we can choose $T$ is large enough such that $\|(a,b)\|_{\LamM L^p_1(\XT)}\leq C.$
Applying Theorem \ref{gaugefixingboundary} over $U$, we get a gauge fixing condition for $(a,b)$. Combing this with equations (\ref{gluingregularity}) and using the bootstrapping method, we get the regularity we want.
\qed

\begin{corollary}
Under the assumption as Theorem $\ref{gluingobstrucion}$, if $H^2_{(A_i,\Phi_i)}=0$ and the limiting flat connection is irreducible, then for $T$ is large enough, there exists a solution to the Kapustin-Witten equations (\ref{KW}).
\end{corollary}

\end{subsection}

\begin{subsection}{Gluing Singular Points in Moduli Space}\label{GlueSingular}
In this subsection, we will deal with the singular points $(A_i,\Phi_i)$ with $H^2_{(A_i,\Phi_i)}\neq 0$. As before, we take the norm $\LamP H^{1,p}_0(X_i)$ on $\Omega^1_{X_i}(\gpp)$ and $\LamM L^p(X_i)$ on $\Omega^0_{X_i}(\gpp)\oplus\Omega^2_{X_i}(\gpp)$. 

As before, we denote $H^2_{(A_i,\Phi_i)}:=(\Omega^2(\gpp)\times\Omega^0(\gpp))/\Im\mathcal{L}_i$. For any $\tau\in \Omega^2(\gpp)\times\Omega^0(\gpp)$, we denote by $[\tau]\in H^2_{(A_i,\Phi_i)}$ the equivalence class of $\tau$. 

We have the following lemma for this cohomology group:
\begin{lemma} Given any bounded open set $U\subset X_i$, for any $\alpha \in H^2_{(A_i,\Phi_i)}$, there exist a $\beta\in\Omega^2(\gpp)\times\Omega^0(\gpp)$, such that $\beta$ is supported in $U$ and $[\beta]=\alpha$ as cohomology class.
\label{finitesupport}
\end{lemma}
\proof
As the range of $\ML_i$ is closed in $\LamM L^p(\Omega^2(\gpp)\times \Omega^0(\gpp)$, we have the following splitting: 
\begin{equation}
\LamM L^p(\Omega^2(\gpp)\times \Omega^0(\gpp))=\Im \ML_i\oplus (Ker \mathcal{L}_i^{\star}\cap \LamM L^p),
\label{identification111}
\end{equation}

where $\ML_i^{\star}$ is the $L^2$ adjoint of $\ML_i$.

Thus we have the identification $H^2_{(A_i,\Phi_i)}\cong \Ker \mathcal{L}_i^{\star}\cap \LamM L^p$. By the classicial unique continuation property of an elliptic operator on the interior \cite{aronszajn1956unique}, for any $\alpha\in \Ker \ML_i^{\star}$, we have $\alpha$ nonvanishing on any interior open set. Denote $l=\dim H^2_{(A_i,\Phi_i)}$, then for an integer $j$, $0 \leq j\leq l$, there exist a basis $\{a_j\}\in H^2_{(A_i,\Phi_i)}$. In addition, we can choose $\{a_j\}$ orthogonal to each other w.r.t the $L^2$ inner product. 

In order to prove the lemma, we only need to prove the statement for one of the base $a_j$. We claim that for any fixed $a_j$, there exists a differential form $f\in \Omega^2(\gpp)\times\Omega^0(\gpp)$ such that $\lan f,a_j\ran\neq 0$ and $f$ vanishes over the boundary, $f|_{\partial U}=0.$ If not, for any $f\in C^{\infty}_0(\Omega^2(\gpp)\times\Omega^0(\gpp))$, we have $\lan f,a_j \ran=0.$ This will imply $a_j$ is identically $0$ over an interior open set which contradicts $a_j\in \Ker\ML_i^{\star}$.

By the Gram–Schmidt process and rescaling, we can find a function $g$ which vanishes over $\partial U$, $\lan g,a_j\ran=1$ and for $s\neq j$ $\lan g,a_s\ran=0$. By the splitting (\ref{identification111}), we know there exist a $g_0\in \Im\ML_i$, such that $g=g_0+a_j.$
\qed

By the previous lemma, we know there exist linear operators $\sigma_i$,  $$\sigma_i:H^2_{(A_i,\Phi_i)}\rightarrow \Omega^2_{X_i}\oplus \Omega^0_{X_i},$$
such that the operators $$\ML_i\oplus \sigma_i:\Omega^1_{X_i}\oplus\Omega^1_{X_i}\oplus H^2_{(A_i,\Phi_i)}\rightarrow \Omega^2_{X_i}\oplus\Omega^0_{X_i}$$ are surjective. By Theorem \ref{FredholmOperator}, we know that $H^2_{(A_i,\Phi_i)}$ is finite dimensional, therefore, we can take the image of $\sigma_i$ to be supported in $X_i(T)$ for $T$ is large enough. 

In the notation above, take $H=H^2_{(A_1,\Phi_1)}\oplus H^2_{(A_2,\Phi_2)}$, we can define a map $\sigma$:
\begin{equation*}
    \begin{split}
        \sigma=\sigma_1+\sigma_2:H\rightarrow \Omega^2_{\XT}\oplus \Omega^0_{\XT}.
    \end{split}
\end{equation*}

As $\ML_i\oplus \sigma_i$ is surjective, there exists an operator $Q_i$, such that $(\ML_i\oplus \sigma_i)Q_i=Id$,
\begin{equation}
Q_i:\Omega^2_{X_i}\oplus\Omega^0_{X_i}\rightarrow H^2_{(A_i,\Phi_i)}\oplus \Omega^1_{X_i}\oplus\Omega^1_{X_i}.
\label{QII}
\end{equation}

Composing $Q_i$ with the projection map into different part of the image, we get operators $\pi_i$ and $P_i$. To be explicit, $Q_i:=\pi_i\oplus P_i$ where $$\pi_i:\Omega^2_{X_i}\oplus\Omega^0_{X_i}\rightarrow H^2_{(A_i,\Phi_i)},$$ and $$P_i:\Omega^2_{X_i}\oplus\Omega^0_{X_i}\rightarrow \Omega^1_{X_i}\oplus\Omega^1_{X_i}.$$ 
Therefore, by definition, for $\forall\xi\in \Omega^2_{X_i}\oplus \Omega^0_{X_i}$, we have 
$$\xi=\ML_iP_i(\xi)+\sigma_i\pi_i(\xi).$$

As before, we take $\phi_i$ be a cut off function supported in $X_i(2T)$ as in Figure \ref{fig:text4871}, with $\phi_i(x)=1$ on $X_i(T)$ and $\phi_1+\phi_2=1$ on $\XT$. We have the estimate $\|\nabla\phi_i\|_{L^{\infty}(X_i)}\leq \epsilon(T)$. 

Given $\xi\in \Omega^2_{\XT}\oplus\Omega^0_{\XT}$, denote by $\xi_i$ the restriction of $\xi$ to $X_i(2T)$, we can define two approximate inverse operators as follows:

Let $\hat{P}(\xi):=\phi_1P_1(\xi_1)+\phi_2P_2(\xi_2)$, $\hat{\pi}(\xi):=\phi_1\pi_1(\xi_1)+\phi_2\pi_2(\xi_2)$. Similarly, we take $\LS:=\mathcal{L}_{(A^{\sharp},\Phi^{\sharp})}$, we have the following lemma:
\begin{lemma}
$\|\LS \hat{P}(\xi)+\sigma\hat{\pi}(\xi)-\xi\|_{\LamM L^p(\XT)} \leq \epsilon(T)\| \xi \|_{\LamM L^p(\XT)}.$
\label{inverse2}
\end{lemma}
\proof Compared to Lemma \ref{inverse1}, we have some additional terms in computing $\LS\hat{P}$.

We have the following computation:
\begin{equation}
\begin{split}
    \LS\hat{P}(\xi)&=\LS(\phi_1 P_1 \xi)+\LS(\phi_2 P_2 \xi)\\
    &=\nabla \phi_1\star P_1(\xi)+\nabla\phi_2\star P_2(\xi)+\phi_1\LS P_1(\xi)+\phi_2\LS P_2(\xi)\\
    &=\nabla \phi_1\star P_1(\xi)+\nabla\phi_2\star P_2(\xi)+\phi_1\ML_1P_1(\xi)+\phi_2\ML_2 P_2(\xi)\\
    &+\phi_1(\LS-\ML_1)P_1(\xi)+\phi_2(\LS-\ML_2) P_2(\xi).
    \label{computeinverse}
\end{split}
\end{equation}

For the terms $\phi_1\ML_1 P_1(\xi)+\phi_2\ML_2 P_2(\xi)$, we have
\begin{equation}
    \begin{split}
        &\phi_1\ML_1 P_1(\xi)+\phi_2\ML_2 P_2(\xi)\\
        =&(\phi_1+\phi_2)(\xi)+\phi_1\sigma_1\pi_1(\xi)+\phi_2\sigma_2 \pi_2(\xi)\\
        =&\xi+\sigma \hat{\pi(\xi)}.
    \end{split}
\end{equation}

For the other terms in the final step of (\ref{computeinverse}), the estimates are exactly the same as Lemma \ref{inverse1} and is bounded by $\epsilon(T)\|\xi\|_{\LamM L^p(\XT)}.$

Combining all the arguement above, we get the estimate we want.
\qed

Now we can construct the inverse of the operator $\LS$.
\begin{corollary}
For $T$ is large enough, there exist operators $\PS:\Omega^2_{\XT}\oplus\Omega^0_{\XT}\rightarrow \Omega^1_{\XT}\oplus\Omega^1_{\XT}$ and 
$\RS:\Omega^2_{\XT}\oplus\Omega^0_{\XT}
\rightarrow H $
 such that 
$\forall\xi\in \Omega^2_{\XT}\oplus\Omega^0_{\XT}$, we have 

\begin{equation}
\xi=\LS\PS(\xi)+\SIS\RS(\xi).
\label{identity2}
\end{equation} 

In addition, the operator norm of $\PS$ and $\RS$ is bounded independent of $T$.

\end{corollary}
\proof
By Lemma \ref{inverse2}, denoting $R:=(\LS\oplus\sigma)(\hat{P}\oplus\hat{\pi})-Id$, we know that when $T$ is large enough, $R$ has operator norm small. Therefore, $Id+R$ is invertible and $Q=(\hat{P}\oplus\hat{\pi})(1+R)^{-1}$ will be the right inverse of $\LS\oplus\sigma$. As the image of $Q$ is $H^2_{(A^{\sharp},\Phi^{\sharp})}\oplus \Omega^1_{\XT}\oplus\Omega^1_{\XT}$, we can take $\PS$ to be the projection to the $\Omega^1\oplus\Omega^1$ part of image of $Q$ and $\RS$ to be the projection of $H^2$ part of image of $Q$, then by definition, we have $$\xi=\LS\PS(\xi)+\SIS\RS(\xi).$$

By classical functional analysis, we know the operator norm of $Id+R$ can be choose to be smaller than 3 and the operator norm of $\hat{P}\oplus\hat{\pi}$ is dominated by $Q_i$ (\ref{QII}). Therefore, the operator norm is independent of $T$.
\qed

For a pair $(\xi,h)$ with $\xi\in\Omega^2_{\XT}\oplus\Omega^0_{\XT}$ and $h\in H^2$, consider the perturbation equation
\begin{equation}
KW((A^{\sharp},\Phi^{\sharp})+\PS(\xi))+\SIS(h)=0.
\label{perturbationequation}
\end{equation}

Therefore, we have
\begin{equation}
\begin{split}
    &KW(A^{\sharp},\Phi^{\sharp})+\LS\PS(\xi)+\{\PS(\xi),\PS(\xi)\}+\SIS(h)=0,\\
&KW(A^{\sharp},\Phi^{\sharp})+\xi-\SIS \RS(\xi)+\{\PS(\xi),\PS(\xi)\}+\SIS(h)=0.(\text{Applying (\ref{identity2})})
\end{split}
\end{equation}

Take $h=\RS(\xi)$, we obtain
\begin{equation}
KW(A^{\sharp},\Phi^{\sharp})+\xi+\{\PS(\xi),\PS(\xi)\}=0
\end{equation}

which is the equation (\ref{quadraticequaion}) and it has solution $\xi$. As $\PS$ is an operator mapping $\Omega^2_{\XT}\oplus\Omega^0_{\XT}$ to $\Omega^1_{\XT}\oplus\Omega^1_{\XT}$, we can define $(a,b)\in \Omega^1_{\XT}\oplus\Omega^1_{\XT}$ by $(a,b):=\PS(\xi)$. Then if we denote $(A,\Phi):=(A^{\sharp}+a,\Phi^{\sharp}+b)$, $(A,\Phi)$ will solve the equation $$KW(A,\Phi)+\sigma(h)=0.$$ 

By the previous arguments, we get the following corollary which completes the proof of the second part of Theorem 1.1.
\begin{corollary}
For any interior open set $U$, there exists $(a,b)\in \Omega^1_{\XT}\oplus\Omega^1_{\XT}$ and $h\in H^2_{(A_1,\Phi_1)}\oplus H^2_{(A_2,\Phi_2)}$ solve the equation $KW(A^{\sharp}+a,\Phi^{\sharp}+b)+\sigma(h)=0$ and satisfy
(1) $(A^{\sharp}+a,\Phi^{\sharp}+b)$ is a solution to the Kapustin-Witten equations over $\XT$ if and only if h=0.

(2) We have the estimate: $$\|(a,b)\|_{\LamM L^p_1}\leq Ce^{-\delta T},\;\|\sigma(h)\|\leq Ce^{-\delta T}.$$This two constants depend on the choice of the open set and $\delta$ is the positive constant in Proposition \ref{small}.

(3) $\sigma(h)$ is supported in $U$.
\label{perturbationequations}
\end{corollary}
\proof The first statement is obvious. For the second statement, by definition, we have 
\begin{equation}
\begin{split}
\|(a,b)\|_{\LamM L^p_1}&=\|\PS(\xi))\|_{\LamM L^p_1}\\
&\leq C\|\xi\|_{\LamM L^p_1} (\text{$\PS$ is bounded})\\
&\leq C\|KW(A^{\sharp},\Phi^{\sharp})\|_{\LamM L^p_1} (\text{By (\ref{decayestimate2})})\\
&\leq Ce^{-\delta T}.(\text{By Proposition \ref{small}})
\end{split}
\end{equation}

Similarly, 
\begin{equation}
\begin{split}
    \|\sigma(h)\|_{\LamM L^p_1}&=\|\sigma(\RS(\xi))\|_{\LamM L^p_1}\\
    &\leq C\|\xi\|_{\LamM L^p_1}(\text{$\RS$ and $\sigma$ are bounded})\\
    &\leq C e^{-\delta T}(\text{By Proposition \ref{small}})\\   
\end{split}
\end{equation}

The third statement is a direct corollary of lemma \ref{finitesupport}.
\qed
\end{subsection}

Given $(A_i,\Phi_i)$, denote by $\Gamma_i$ the isotropy group of $(A_i,\Phi_i)$, $\Gamma_i=\{g|g(A_i,\Phi_i)=(A_i,\Phi_i)\}$. By Corollary \ref{irreducible1}, we know $\Gamma_i=1$.
We will combine the Kuranishi descriptin in Proposition \ref{Kuranishimodel} with the previous construction. Let $N_i\subset H^1_{(A_i,\Phi_i)}$ be a set parametrize a neighborhood of $(A_i,\Phi_i)$ in the moduli space of Nahm pole solutions. If we denote $N:=N_1\times N_2$, then we have the following proposition:
\begin{proposition}
For large enough $T$ and small enough $N_i$, given $n\in N$ then we have

(1) A family of $\LamP H^{1,p}_0$ connections $(A(n),\Phi(n))+(a(n),b(n))$ parametrized by $N$. 

(2) There exist a map $\Psi:N\rightarrow H^2_{(A_1,\Phi_1)}\times H^2_{(A_2,\Phi_2)}$, such that $(A(n),\Phi(n))+(a(n),b(n))$ satisfies the Kapustin-Witten equations if and only if $\Psi(n)=0$.

(3) Let $\MM_{\XT}$ be the moduli space of Nahm pole solutions to the Kapustin-Witten equations over $X$, then there exists a map $\Theta$, whose image is the moduli space $\MM_{\XT}$:
\begin{equation}
\begin{split}
    \Theta:\Psi^{-1}(0)&\rightarrow \MM_X\\
n&\rightarrow (A(n),\Phi(n))+(a(n),b(n)).\\
\end{split}
\end{equation}
\label{97111}
\end{proposition}

\begin{subsection}{Gluing for Non-degenerate Limit}
In this section, we will build the gluing theorem for the reducible connection. For simplicity, in this subsection, we only consider the case $H^2_{(A_i,\Phi_i)}(X_i)=0$. For the $H^2$ non-vanishing case, the result will follows similarly as in subsection 7.3.

As before, we are dealing with manifolds $X_1$, $X_2$ with cylindrical ends and boundaries as in Figure \ref{fig:text4463}. We constructed $\XT$, identified the connecting region with $Y\times(-\frac{T}{2},\frac{T}{2})$. This will be more precisely shown in Figure
\ref{fig:path3720}.
\begin{figure}[H]
    \centering
    \includegraphics[width=0.8\textwidth]{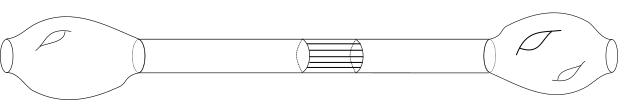}
    \caption{The shaded part is $Y\times(-\frac{T}{2},\frac{T}{2})$}
    \label{fig:path3720}
\end{figure}

For a positive real number $\alpha$, take a smooth weighted function $W_T=e^{\alpha(\frac{T}{2}-|t|)}$ and over a neighborhood of the boundary of $\XT$, let $W_T$ be the distance function to the boundary.

Over the manifolds with boundary and cylindrical ends $X_1$ and $X_2$, we have fixed weighted functions $W_1$ and $W_2$, such that in the connected area $W_1=e^{\alpha(\frac{T}{2}+t)}$, $W_2=e^{\alpha(\frac{T}{2}-t)}$ and in the neighborhood of the boundary, $W_1$ and $W_2$ are the distance functions to the boundaries. It is easy to get that in the common area $W_1$, $W_2$ and $W_T$ dominated each other.

On $1$-forms of $\XT$, use the norm $\LamP H^{1,p}_{0,\alpha}(\XT)$ given by the weighted norm given by $H^{1,p}_0(\XT)$ and weight function $W_T$. On the $2$-forms of $\XT$, use the norm $\LamM L^p_{\alpha}(\XT)$ given by $L^p(\XT)$ and weighted function $W_T$. Respectively, we get $\LamP H^{1,p}_{0,\alpha}(X_i)$ and $\LamM L^p_{\alpha}(X_i)$ for $X_i$.

By these constructions, we get the following estimate for the approximate solution:

\begin{proposition}
$\|KW(A^{\sharp},\Phi^{\sharp})\|_{_{\LamM L_{\alpha}^p(\XT)}}\leq C(e^{(\alpha-\delta)T}).$
\label{lsmall}
\end{proposition}
\proof
By Theorem $\ref{exponentialdecay}$, we know the $C^{\infty}$ norm will decays as $e^{-\delta t}$. In addition, we have the weighte function that equals to $e^{\alpha t}$ in the end. Therefore, we get the decay rate we want. 
\qed

Therefore, we can take $\alpha< \delta$ such that the approximate term exponentially decays as $T\rightarrow \infty$.

For $i=1,2$, denoting $\ML_i:=\ML_{(A_i,\Phi_i)}$, we can regard the operator as
$$\ML_{i,\alpha}:\LamP H^{1,p}_{0,\alpha}(X_i)\rightarrow \LamM L^p_{\alpha}(X_i).$$

For the approximate solution $(A^{\sharp},\Phi^{\sharp})$, we also have the Fredholm operator $\LS_{\alpha}$ for the weighted norm
$$\LS_{\alpha}:\LamP H^{1,p}_{0,\alpha}(\XT)\rightarrow \LamM L^p_{\alpha}(\XT).$$

By our assumption $H^2_{(A_i,\Phi_i)}(X_i)=0$, we know there exists a right inverse $Q_i$: $$Q_i:\LamM L^p_{\alpha}(X_i)\rightarrow \LamP H^{1,p}_{0,\alpha}(X_i),$$ such that $\ML_{i,\alpha}Q_i=Id$.

As before, we take $\phi_i$ be a cut off function supported in $X_i(2T)$ as in Figure \ref{fig:text4871}, with $\phi_i(x)=1$ on $X_i(T)$ and $\phi_1+\phi_2=1$ on $\XT$ and we can have the estimate $\|\nabla\phi_i\|_{L^{\infty}(X_i)}\leq \epsilon(T)$. 

Take $\xi\in \LamM L^p_{\alpha}(\XT)$, denote by $\xi_i$ the restriction of $\xi$ on $X_i(2T)$, then we define a new approximate inverse operator $\hat{Q}_{\alpha}(\xi):=\phi_1Q_1(\xi_1)+\phi_2Q_2(\xi_2)$, which can be written as $$\hat{Q}:\LamM L_{\alpha}^p(X_i)\rightarrow \LamP H^{1,p}_{0,\alpha}(X_i).$$

\begin{subsubsection}{Right Inverse}
Similarily, we have the following estimate for the operator $\LS_{\alpha}$.

\begin{lemma}\label{inverse3}
For $\LS_{\alpha}$, $\hat{Q}_{\alpha}$ as above, and for $\forall \xi\in \LamM L^p_{\alpha}(\XT)$, we have 
$$\|\LS_{\alpha}\hat{Q}_{\alpha}(\xi)-\xi\|_{\LamM L^p_\alpha(\XT)}\leq \epsilon(T)\|\xi\|_{\LamM L^p_\alpha(\XT)}.$$
\end{lemma}
\proof
After we choose $\alpha< \delta$, we still get the exponential decay result and the proof is exactly the same as Lemma \ref{inverse1}.
\qed
\begin{proposition}
There exists an operator $\QS_{\alpha}$,$$\QS_{\alpha}:\LamM L^p_{\alpha}(\XT)\rightarrow \LamP H^{1,p}_{0,\alpha}(\XT)$$ 
such that $\LS_{\alpha}\QS_{\alpha}=Id.$ In addition, the operator norm of $QS_{\alpha}$ is independent of T.
\end{proposition}
\proof
By Lemma \ref{inverse3}, we know $\LS_{\alpha}\hat{Q}_{\alpha}$ has an inverse and we just take $\QS_{\alpha}:=\hat{Q}_{\alpha}(\LS_{\alpha}\hat{Q}_{\alpha})^{-1}$. By definition, we get the inverse we want. For the independence of $T$ from the operator norm, the arguement is exactly the same as Proposition \ref{inverseoperator}.
\qed
\end{subsubsection}

\begin{subsubsection}{Existence Theorem}
Over $\XT$, for arbitary $(A,\Phi)$, we denote $(a,b):=(A,\Phi)-(A^{\sharp},\Phi^{\sharp})$. We have the following expansion for the Kapustin-Witten map.
\begin{equation}
KW(A,\Phi)=KW(A^{\sharp},\Phi^{\sharp})+\LS_{\alpha}(a,b)+\{(a,b),(a,b)\}.
\label{qe2}
\end{equation}

Take $S_{\alpha}(a,b)=\{(a,b),(a,b)\}$, then we have the following proposition.
\begin{proposition}\label{SSS}
For any $\lambda_0\in[1-\frac{1}{p},1)$, $S_{\alpha}$ is an operator $S_{\alpha}:\LamM L^p(\XT)\rightarrow \LamM L^p(\XT)$ satisfying
$S_{\alpha}(0)=0$ and for two elements $\beta,\gamma\in \LamM L_{\alpha}^2(\XT)$, there exists a constant $k$ such that $$\|S_{\alpha}(\beta)-S_{\alpha}(\gamma)\|_{\LamM L_{\alpha}^p}\leq k(\| \beta \|_{\LamM L_{\alpha}^p}+\|\gamma\|_{\LamM L_{\alpha}^p})(\|\beta-\gamma\|_{\LamM L_{\alpha}^p}).$$
\end{proposition}
\proof
The proof is basically the same as the proof of Proposition \ref{SS}. We only need to check the Sobolev inequality is still true in the weighted case and this is proved in Proposition \ref{SoI}.
\qed
\end{subsubsection}

Now, we have a parallel theorem to Theorem 1.1:

\begin{theorem}
Under the gluing hypotheses in the beginning of the chapter, if

(a) $\lim_{T\rightarrow +\infty}\|(A_i,\Phi_i)-(A_{\rho_i},\Phi_{\rho_i})\|_{L^{p_0}_2(Y_i\times\{T\})}=0$ for some $p_0>2$,

(b) $\rho$ is a non-degenerate $\slc$ flat connection,

then for $\lambda_0\in[1-\frac{1}{p},1)$, there exists a real number $\alpha> 0$ such that we have

(1) for some constant $\delta$, there exists a $y^{\lambda_0+\frac{1}{2}}H^{1,p}_{0,\alpha}$ pair $(a,b)\in \Omega^1_{\XT}(\gpp)\times\Omega^1_{\XT}(\gpp)$ with $$\|(a,b)\|_{y^{\lambda_0-\frac{1}{2}}L_{1,\alpha}^p}\leq Ce^{(\alpha-\delta) T},$$

(2) there exists an obstruction class $h\in H^2_{(A_1,\Phi_1)}(X_1)\times H^2_{(A_2,\Phi_2)}(X_2)$ such that $h=0$ if and only if $(A^{\sharp}+a,\Phi^{\sharp}+b)$ is a solution to the Kapustin-Witten equations (\ref{KW}).
\label{glureducible}
\end{theorem}
\proof
For the case that $H^2$ vanishes, by Proposition \ref{SSS}, we know that the opertor $S_{\alpha}$ satisfies the assumption for Lemma \ref{LS}. 
By Proposition \ref{lsmall}, we know we can choose $\eta$ small enough satifying Lemma \ref{LS}. Therefore, by Lemma \ref{LS}, there exists a solution to the equation (\ref{qe2}) and we get a solution to the Kapustin-Witten equations (\ref{KW}). The regularity statements of the connections in the theorem will follows by the same way as in Chapter 8.1.

Similarly, we can follow exactly the same as Chapter 8.3 and prove the second statement of the theorem.
\end{subsection}
\end{section}

\begin{section}{Local Model for Gluing Picture}
In this section, we will give a Kuranishi description of the gluing construction for the Kapustin-Witten equations.

For the description for the anti-self-dual equations, see \cite{donaldson2002floer}, \cite{donaldson1990geometry}, \cite{taubes1982self}. In this section, we assume $p\in (2,4)$, $\lambda\in[1-\frac{1}{p},1)$ and denote by $q$ the real number satisfying the relationship $1+\frac{4}{q}=\frac{4}{p}.$ 

\begin{subsection}{Gauge Fixing Problem}

For $i=1,2$, let $\MM_i$ be the moduli space of Nahm pole solutions to the Kapustin-Witten equations over $X_i$ defined in (\ref{modulispace1}). Let $N_i$ be pre-compact subsets of the moduli space $\MSI$ such that any element of $N_i$ is regular in the moduli space $\MSI$. To be more explicit, for any $(A_i,\Phi_i)\in N_i$, we have $H^0_{(A_i,\Phi_i)}=0$ and $H^2_{(A_i,\Phi_i)}=0$. By Proposition $\ref{97111}$, we know there exists a map $\Theta_T$ defined as follows:
\begin{equation}
\Theta_T:N_1\times N_2\rightarrow \MSZ
\end{equation}

We have the following proposition on the map $\Theta$:
\begin{proposition}
There exists a $T_0$, such that for any $T> T_0$, we have:

(1) For $(A_i,\Phi_i)\in N_i$, let $(A,\Phi):=\Theta_T((A_1,\Phi_1),(A_2,\Phi_2))$, we have $H^2_{(A,\Phi)}=0,$

(2) $\Theta_T$ is a diffeomorphism to its image.
\label{91}
\end{proposition}
\proof (1) Let $(A^{\sharp},\Phi^{\sharp})$ be the approximate solution, let $(a,b):=(A,\Phi)-(A^{\sharp},\Phi^{\sharp})$. By Theorem $\ref{gluingobstrucion}$, we have $\|(a,b)\|_{\LamM L^p_1}\leq Ce^{-CT}$. Let $\mathcal{L}^1_{(A^{\sharp},\Phi^{\sharp})}$ ($\mathcal{L}^1_{(A,\Phi)}$ ) be the linearization operator of $(A^{\sharp},\Phi^{\sharp})$ ($\mathcal{L}_{(A,\Phi)}$). By Proposition \ref{inverseoperator}, there exists an operator $\mathcal{Q}^{\sharp}:\LamM L^p\rightarrow \LamP H^{1,p}_0$ such that $\mathcal{L}^1_{(A^{\sharp},\Phi^{\sharp})}\mathcal{Q}^{\sharp}=Id.$ 
Therefore, we can choose $T$ big enough such that  $\|\mathcal{L}^1_{(A,\Phi)}\mathcal{Q}^{\sharp}-\mathcal{L}^1_{(A^{\sharp},\Phi^{\sharp})}\mathcal{Q}^{\sharp}\|_{\LamM L^p}\leq \frac{1}{2}$. This implies that $\mathcal{L}^1_{(A,\Phi)}$ has a right inverse and it is surjective.

(2) By the assumption that $N_i$ is regular, we have $\dim N_i= \Ind P_i$. By Proposition $\ref{glueindex}$, we have $\Ind P=\Ind P_1+ \Ind P_2$. Let $\Im(\Theta)$ to be the image of $\Theta$. we have $\dim(\Im(\Theta))=\dim N_1+\dim N_2.$ Therefore, in order to prove $\Theta_T$ is a diffeomorphism, we only need to prove $d\Theta$ is injective. Choose an open subset $U\subset X_1$ which is away from the gluing part. By Proposition \ref{regularityofgluing}, we know over $U$, $(A,\Phi)$ is $C^1$ close to $(A_1,\Phi_1)$ thus proves that $d\Theta$ is injective. 
\qed

Now we will characterize the Nahm pole solutions we found by our gluing construction. Given solutions $(A_i,\Phi_i)$. Let $(A^{\sharp},\Phi^{\sharp})$ be the approximate solution. Let $d^{\lambda}_q$ be the metric on the space $\mathcal{B}$ given by
\begin{equation}
d^{\lambda}_q([(A_1,\Phi_1)],[(A_2,\Phi_2)])=\inf_{u\in\mathcal{G}}\|(A_1,\Phi_1)-u(A_2,\Phi_2)\|_{\LamM L^q}.
\end{equation}

Then, we can define an open neighborhood $U(\epsilon)$ of $(A^{\sharp},\Phi^{\sharp})$ by 
\begin{equation}
U_{(A^{\sharp},\Phi^{\sharp})}(\epsilon)=\{(A,\Phi)\in \mathcal{B}|d_q^{\lambda}((A,\Phi),(A^{\sharp},\Phi^{\sharp}))|< \epsilon,\; \|KW(A,\Phi)\|_{\LamM L^p}< \epsilon\}.
\label{neighborhooda}
\end{equation}

Then we have the following theorem
\begin{theorem} For $\star=0,1,2$, if $H^{\star}_{(A_i,\Phi_i)}=0$, then for small enough $\epsilon$, any point $(A,\Phi)\in U(\epsilon)$ can be represented by the following form
$(A,\Phi)=(A^{\sharp},\Phi^{\sharp})+\MQ\phi$, where 
$\|\phi\|_{\LamM L^p}\leq C\epsilon$ and $\MQ$ is the right inverse operator defined in Proposition \ref{inverseoperator}.
\label{gaugefixing}
\end{theorem}

We prove Theorem \ref{gaugefixing} by the method of continuation. We need a new interpretation of the operator.

Given $(A_i,\Phi)$ satisfying the assumption of Theorem \ref{gaugefixing}, let $(A^{\sharp},\Phi^{\sharp})$ to be the approximate solution over $\XT$. In this section, for simplification, we denote $\ML$ the linearization operator of $(A^{\sharp},\Phi^{\sharp})$ and let $Q$ be the right inverse of $\ML$. Combining this with the embedding $\LamP H^{1,p}_0\hookrightarrow \LamM L^q$, we have
\begin{equation}
\begin{split}
    &\ML:\LamP H^{1,p}_0(\Omega^1(\gpp)\times\Omega^1(\gpp))\rightarrow \LamM L^p(\Omega^2(\gpp)\times\Omega^0(\gpp)),\\
    &Q:\LamM L^p(\Omega^2(\gpp)\times\Omega^0(\gpp))\rightarrow \LamM L^q(\Omega^1(\gpp)\times\Omega^1(\gpp)).
\end{split}
\end{equation}

Let $B\in U_{(A^{\sharp},\Phi^{\sharp})}(\epsilon)$, then WLOG, we assume $B=(A^{\sharp},\Phi^{\sharp})+(a,b)$ and consider $B_t$ which is a path of connection pairs defined as follows:
$$B_t:=(A^{\sharp},\Phi^{\sharp})+t(a,b)$$ and we can define the following set $S$:

\begin{definition}
Given $\delta$ small enough, define $S\subset [0,1]$ to be the interval of all $t\in [0,1]$ such that there exists gauge transform $u:[0,t]\rightarrow \mathcal{G}$ and $\phi:[0,t]\rightarrow \Omega^2(\gpp)\times\Omega^0(\gpp)$ such that 

$(1)\; \phi(0)=0,\;u(0)=1,$

$(2)\; u_t(B_t)=(A^{\sharp},\Phi^{\sharp})+Q(\phi_t)$ with $\|\phi_t\|_{\LamM L^p}< \delta.$
\label{SetS}
\end{definition}

Our target is to prove $S=[0,1]$. By definition of $S$, we have the following Proposition:
\begin{proposition}
$S$ is non empty.
\label{nonempty}
\end{proposition}
\proof
As $B_0=(A^{\sharp},\Phi^{\sharp})$, take $\phi_0=0$ and $u(0)=1$, we know $0\in S.$
\qed

Now, we are going to prove $S$ is an open set and before the proving, we will need some preparations.

Let $d^0$ to be $d^0_{(A^{\sharp},\Phi^{\sharp})}$ in the Kuranishi complex 
($\ref{complex}$), where for $\xi\in\Omega^0(\gpp)$, $d^0(\xi)=(-d_{A^{\sharp}}\xi,[\xi,\Phi^{\sharp}])$. For any $\xi\in\Omega^0(\gpp)$ and $\phi\in\Omega^2(\gpp)\times\Omega^0(\gpp)$, define the operator 
\begin{equation}
\begin{split}
    \Pi:&\Omega^0(\gpp)\times\Omega^2(\gpp)\times\Omega^0(\gpp)\rightarrow \Omega^1(\gpp)\times\Omega^1(\gpp),\\
    &(\xi\;,\;\phi)\rightarrow d^0(\xi)+Q(\phi).
\end{split}
\end{equation}

Let $V_1$ to be a norm over $\Omega^0(\gpp)\times\Omega^2(\gpp)\times\Omega^0(\gpp)$ defined as follows: $$\|(\xi,\phi)\|_{V_1}=\|d^0(\xi)\|_{\LamM L^q}+\|\phi\|_{\LamM L^p}.$$ For $(a,b)\in \Omega^1(\gpp)\times\Omega^1(\gpp)$, we define another norm $V_2$ as
$$\|(a,b)\|_{V_2}=\|(a,b)\|_{\LamM L^q}+\|\ML(a,b)\|_{\LamM L^p}.$$

Then we have the following Proposition:
\begin{proposition}
Considering $\Pi$ as operator from $V_1$ to $V_2$:$$\Pi:V_1\rightarrow V_2,$$ we have

(1) $\Pi$ is a bounded operator from $V_1$ to $V_2$,

(2) There exist a constant $C$ independent of $T$ such that 
$\|(\xi,\phi)\|_{V_1}\leq C\|\Pi(\xi,\phi)\|_{V_2}$.
\label{longpeop}
\end{proposition}
\proof
(1) We have the following computation for the operator $\Pi$:
\begin{equation}
\begin{split}
    &\|\Pi(\xi,\phi)\|_{B_2}\\
\leq&\|d^0(\xi)+Q(\phi)\|_{\LamM L^q}+\|\mathcal{L}\circ d^0(\xi)+\mathcal{L}\circ Q(\phi)\|_{\LamM L^p}\\
&(\mathrm{Here\; we \;use\; \ML\circ d^0(\xi)=[KW(A^{\sharp},\Phi^{\sharp}),\xi]\;and\; \ML\circ Q=Id })\\
\leq& \|d^0(\xi)\|_{\LamM L^q}+\|\phi\|_{\LamM L^p}+\|[KW(A^{\sharp},\Phi^{\sharp}),\xi]\|_{\LamM L^p}+\|\phi\|_{\LamM L^p}\\
\leq& \|d^0(\xi)\|_{\LamM L^q}+\|\phi\|_{\LamM L^p}(\rm{Here\; we\; use\; Proposition\; \ref{small}} \;and\;|\xi|_{C^0}\leq\|d^0(\xi)\|_{\LamM L^q}).
\end{split}
\end{equation}

(2) Take $\alpha=d^0(\xi)+Q(\phi)$, then we have $\ML \alpha=[KW(A^{\sharp},\Phi^{\sharp}),\xi]+\phi$. We have the following estimate:
\begin{equation}
\begin{split}
    \|d^0(\xi)\|_{\LamM L^q}&\leq \|\alpha -Q\phi\|_{\LamM L^q}\\
    &\leq \|\alpha\|_{\LamM L^q}+\|Q\phi\|_{\LamM L^q}\\
    &\leq \|\alpha\|_{\LamM L^q}+\|\phi\|_{\LamM L^p}.
    \label{86}
\end{split}
\end{equation}

In addition, by the relation $\ML \alpha=[KW(A^{\sharp},\Phi^{\sharp}),\xi]+\phi$, we have
\begin{equation}
\begin{split}
    \|\phi\|_{\LamM L^p}&\leq\|\ML\alpha\|_{\LamM L^p}+\|[KW(A^{\sharp},\Phi^{\sharp}),\xi]\|_{\LamM L^p}\\
    &\leq\|\alpha\|_{V_2}+\epsilon\|\xi\|_{C^0}(Here\; we\; use\; Proposition\; \ref{small})\\
    &\leq\|\alpha\|_{V_2}+\epsilon \|d^0(\xi)\|_{{\LamM L^q}}\\
    &=\|\alpha\|_{V_2}+\epsilon \|\alpha-Q\phi\|_{\LamM L^q}\\
    &\leq \|\alpha\|_{V_2}+\epsilon \|\alpha\|_{\LamM L^q}+\epsilon \|\phi\|_{\LamM L^p}.
    \label{87}
\end{split}
\end{equation}
By taking $\epsilon$ small enough, we get $\|\phi\|_{\LamM L^p}\leq C\|\alpha\|_{V_2}.$

By definition, $\|(\xi,\phi)\|_{V_1}=\|d^0(\xi)\|_{\LamM L^q}+\|\phi\|_{\LamM L^p}.$ Combining equations (\ref{86}), (\ref{87}), we obtain 
$$\|(\xi,\phi)\|_{V_1}\leq C\|\alpha\|_{V_2}=C\|\Pi(\xi,\phi)\|_{V_2}.$$
\qed.

By this estimate, we get an immediate corollary:
\begin{corollary}
$\Pi$ is an injective operator.
\end{corollary}

\begin{proposition}
For $\star=0,1,2$, if $H_{(A_i,\Phi_i)}^{\star}=0$, the operator $\Pi$ is a surjective operator from $V_1$ to $V_2$
\label{surjectiveopen}
\end{proposition}
\proof
As $Q$ is the inverse of $\ML$, by the assumption $H_{(A_i,\Phi_i)}^{\star}=0$, we know $\Ind\;\Pi=-\Ind\;\DKW=0.$ By Proposition \ref{longpeop}, we know $\Pi$ is injective, thus $\Pi$ is surjective. 
\qed

\begin{proposition}
$S$ is an open set in $[0,1]$.
\label{openset}
\end{proposition}
\proof
By Proposition \ref{surjectiveopen}, $\Pi$ is surjective. By the implicit function theorem, we get the result immediately. 
\qed

Now, we hope to prove that the set $S$ is a closed set. To begin with, we prove that the condition (2) in Definition \ref{SetS} is a closed condition:

\begin{lemma}
For suitable $\delta$ and $\epsilon$, we have $\|\phi_t\|_{\LamM L^p}\leq \frac{1}{2}\delta$.
\label{closecondition}
\end{lemma}
\proof
By the relation $u_t(B_t)=(A^{\sharp},\Phi^{\sharp})+Q(\phi_t)$, we have:
\begin{equation}
KW(u_t(B_t))=KW(A^{\sharp},\Phi^{\sharp})+\phi_t+\{Q(\phi_t),Q(\phi_t)\}.
\end{equation}
Therefore, we have
\begin{equation}
\begin{split}
  \|\phi_t\|_{\LamM L^p}&\leq \|KW(A^{\sharp},\Phi^{\sharp})\|_{\LamM L^p}+\|KW(B_t)\|_{\LamM L^p}+\|Q(\phi_t)\|^2_{\LamM L^q}\\
  &(Here\; we\; use\; Proposition\; \ref{small}\;and\;definition\;(\ref{neighborhooda}) )\\
  &\leq \epsilon(T)+\epsilon+C^2\|\phi_t\|^2_{\LamM L^q}
\end{split}
\end{equation}
For $\delta< \frac{1}{2C^2}$, $\epsilon(T)\leq \frac{1}{4}\delta$ and $\epsilon\leq\frac{1}{4}\delta$, we have $\|\phi_t\|\leq \frac{1}{2}\delta$, so the open condition is also closed.
\qed

\begin{proposition}
For $\delta$ small enough and suitable parameter $T$ and $\epsilon$, $S$ is a closed set in $[0,1]$.
\label{closeset}
\end{proposition}
\proof
Now is routine to prove the set $S$ is closed. Let assume a sequence $t_i\in S$ with $t_i\rightarrow t_0$. For simplification, we denote $B_i:=B_{t_i}$ and $\phi_i:=\phi_{t_i}$. By the definition of $S$, we have the relationship $u_t(B_i)=(A^{\sharp},\Phi^{\sharp})+Q(\phi_i)$.

By Lemma \ref{closecondition}, we have the closed condition $\|\phi_t\|_{\LamM L^p}\leq \frac{1}{2}\delta$. By definition of $B_i$, we have $B_i=(A^{\sharp},\Phi^{\sharp})+t_i(a,b)$ and $(a,b)\in \LamP H^{1,p}_0\subset \LamM L^p_1$. We know $B_i$ strongly converges in $\LamM L^p_1$.

By the uniform bound on the $\phi_{i}$, the $\phi_{i}$ converges to a limit $\phi_0$ weakly in $\LamM L^p$. Define $A_i=(A^{\sharp},\Phi^{\sharp})+Q(\phi_i)$. $A_i$ is uniformly bounded in $\LamP H^{1,p}_0\hookrightarrow \LamM L^p_1$. Therefore, $A_i$ converges weakly in $\LamM L^{p}_1$. 

As $u_i$ is a gauge transformation, by the relation $u_i(B_i)=A_i$, we have $du_i=u_iA_i-B_iu_i.$ By the boundedness of $A_i$ and $B_i$, we know $u_i$ weakly converges to $u_0$ in $\LamM L^p_2$. Therefore, by the Sobolev embedding theorem, $u_i$ strongly converges in $\LamM L^p_1$ to $u_0$. Therefore, we have the relationship $u_0(B_0)=A_0$ which imply $t_0\in S$.
\qed
\end{subsection}

We get an immediate corollary from Proposition \ref{nonempty}, Proposition \ref{openset} and Proposition \ref{closeset}:

\begin{corollary}
For the set $S$ in definition \ref{SetS}, we have $S=[0,1]$.
\end{corollary}
The proof of Theorem \ref{gaugefixing} follows immediately.

\begin{subsection}{Local Model for Regular Moduli Space}
Now, we are able to construct a local model for the gluing picture in the acyclic case without the assumption on $H^1$. 

Denote $n_i=\Ind(P_i)$ and we don't assume $n_i=0$. Denote $\MSIS$ ($\MSPS$) to be the moduli space which only consists of solutions to the Kapustin-Witten equations over $X_i$ ($\XT$) which have $H^2=0.$

For $i=1,2$, given two solutions $(A_i,\Phi_i)\in \MSIS$, there exists a open neighborhood $U_i$ such that we can find functions $$\chi:U_i\subset \MSIS\rightarrow \mathbb{R}^{n_i}$$ which give local coordinates around $(A_i,\Phi_i)$ in the moduli spaces $\MSIS$. Denote 
\begin{equation}
U_P(\epsilon)=\{(A,\Phi)\in \mathcal{B}|\exists (A_0,\Phi_0)\in \MSPS\; with\; d_q^{\lambda}((A,\Phi),(A_0,\Phi_0))< \epsilon,\; \|KW(A,\Phi)\|_{\LamM L^p}< \epsilon\}.
\label{closeneighborhood}
\end{equation}

Then by the exponential decay result (Theorem \ref{exponentialdecay}), we know that by choosing suitable compact sets $G_i\subset X_i$ and cut-off functions, we have a natural inclusion $U_i$ into $\MSIS$. Choose $y_i\in Im(\chi_i(U_i))$ and define the cut-down moduli space
$$L=\chi_1^{-1}(y_1)\cap \chi_2^{-1}(y_2)\cap\MSPS \subset U_P(\epsilon)$$ has virtual dimensional $0$.

For $T$ is large enough, recall $I:\MC_{P_1}\times\MC_{P_2}\rightarrow \MC_P$ is the operator defined in (\ref{approximatemap}) that constructs the approximate solution. Denote by $(A_0,\Phi_0):=I(\chi_1^{-1}(y_1),\chi_2^{-1}(y_2))$ the approximate solution constructed by $\chi_1^{-1}(y_1)$ and $\chi_2^{-1}(y_2)$.
Then we have the following Proposition. Compare this to Theorem \ref{gaugefixing}:

\begin{proposition}
For $\epsilon$ small enough, there exists a unique solution $(A',\Phi')$ in $L$ such that $U_{(A_0,\Phi_0)}(\epsilon)\cap L=(A',\Phi').$
\end{proposition}

Now, we will define a distance to make a comparision between connection pairs $(A_0,\Phi_0)$ over $\XT$ and $(A_i,\Phi_i)$ over $X_i$. 

We can define the norm $d$ as 
\begin{equation}
d((A^{\sharp},\Phi^{\sharp});(A_1,\Phi_1),(A_2,\Phi_2))=inf_{u\in\GP}\|(A_0,\Phi_0)-I((A_1,\Phi_1),(A_2,\Phi_2))\|_{L^q(\XT)}.
\end{equation}

where the $I$ is the operator that constructs the approximate solutions defined in (\ref{approximatemap}).

Summarizing Proposition \ref{91} and Theorem \ref{gaugefixing}, we obtain the following statement:
\begin{theorem}
Denote by $U_i$ the compact sets of regular points in the moduli space $\MSIS$. There exist $T_0$, $\epsilon_0$ such that for $T> T_0$ and $\epsilon< \epsilon_0$, there exist open neighborhoods $N_i$ of $U_i$ and a map
$$\Theta:N_1\times N_2\rightarrow \MSPS,$$ such that

$(1)\;\Theta$ is a diffeomorphism to its image, and the image contains regular points,

$(2)$ $d(\Theta((A_1,\Phi_1),(A_2,\Phi_2));(A_1,\Phi_1),(A_2,\Phi_2)))\leq \epsilon$ for any $(A_i,\Phi_i)\in N_i$,

$(3)$ Any connection $(A^{\sharp},\Phi^{\sharp})\in\MSPS$ with $d((A^{\sharp},\Phi^{\sharp});(A_1,\Phi_1),(A_2,\Phi_2))\leq \epsilon$ for some $(A_i,\Phi_i)\in N_i$ lies in the image of $\Theta$.
\label{RegularLocalModel}
\end{theorem}
\end{subsection}

Now we will have a brief discussion of the local gluing picture in the general case the $H^2$ is non-vanishing. For $(A_i,\Phi_i)\in \MSI$ with $H^2(A_i,\Phi_i)$ non-vanishing, we can do the trick as in Section \ref{GlueSingular} by adding some finite dimensional linear space as the obstruction class and have a similar obstruction type statement as in Theorem $\ref{RegularLocalModel}$. We will precise by state the theorem in general in the next subsection.

\begin{subsection}{Conclusions}
Now, we can summarize what we have proved and state the following theorem

\begin{theorem}
Let $(A_i,\Phi_i)$ be connections pairs over manifolds $X_i$ with Nahm poles over $Z_i$, for sufficiently large $T$, there is a local Kuranishi model for an open set in the moduli space over $\XT$: 

$(1)$ There exists a neighborhood $N$ of $\{0\}\subset H^1_{(A_1,\Phi_1)}\times H^1_{(A_2,\Phi_2)}$ and a map $\Psi$ from $N$ to $H^2_{(A_1,\Phi_1)}\times H^2_{(A_2,\Phi_2)}$.

$(2)$ There exists a map $\Theta$ which is a homeomorphism from $\Psi^{-1}(0)$ to an open set $V\subset \mathcal{M}_{X^{\sharp}}^{\star}.$

\end{theorem}

\end{subsection}
\end{section}

\begin{section}{Some Applications}
In this section, we will introduce some applications of the gluing theorem \ref{gluingobstrucion}.

As for the model solution in Section 2, we don't know whether the obstruction class vanish or not and right now we don't have any transversality result for the Kapustin-Witten equations. We just consider the obstruction class as a perturbation to the equation. See \cite{donaldson1983application} for the obstruction perturbation for ASD equations.

Consider a compact 4-dimensional manifold $X^4$ with a cylindrical end which is identified with $Y^3\times [0,+\infty)$, given any $\slc$ representation $\rho$ of $\pi_1(X^4)$:
$$\rho:\pi_1(X^4)\rightarrow \slc,$$ denote by $(A_\rho,\Phi_\rho)$ the $\slc$ flat connection associated to $\rho$. Then we know $(A_\rho,\Phi_\rho)$ satisfies the following equations:
\begin{equation}
\begin{split}
    &F_{A_{\rho}}-\Phi_{\rho}\wedge\Phi_{\rho}=0,\\
    &d_{A_{\rho}} \Phi_{\rho}=0,\\
    &d_{A_{\rho}}^{\star} \Phi_{\rho}=0.
\end{split}
\end{equation}

Obviously, $(A_\rho,\Phi_\rho)$ is a solution to the Kapustin-Witten equations (\ref{KW}).

By gluing the suitable $\slc$ flat connection, we have the following theorem: 
\begin{theorem}
Consider a smooth compact 4-manifold $M$ with boundary $Y$. Assume $Y$ is $S^3$, $T^3$ or any hyperbolic 3-manifold. For $Y$ is hyperbolic, we assume the inclusion of $\pi_1(Y)$ into $\pi_1(M)$ is injective. For a real number $T_0$, we can glue $M$ with $Y\times (0,T_0]$ along $\partial M$ and $Y\times \{T_0\}$ to get a new manifold, which denote as $M_{T_0}:=Y\times (0,T_0)\cup M$. For $T_0$ large enough, there exists an 
$SU(2)$ bundle $P$ and its adjoint bundle $\gpp$ over $M_{T_0}$ such that given any interior non-empty open neighborhood $U\subset M$, we have:

(1) There exist $h_1\in\Omega^2_{M_{T_0}}(\gpp),\;h_2\in\Omega^0_{M_{T_0}}(\gpp)$ supported on $U$,

(2) There exists a connection $A$ over $P$ and a $\gpp$-valued 1-form $\Phi$ such that $(A,\Phi)$ satisfies the Nahm pole boundary condition over $Y\times \{0\}\subset M_{T_0}$ and $(A,\Phi)$ is a solution to the following obstruction perturbed Kapustin-Witten equations over $M_{T_0}$: 
\begin{equation}
\begin{split}
    F_A-\Phi\wedge\Phi+\star d_A\Phi&=h_1,\\
    d_A^{\star}\Phi&=h_2.
\end{split}
\end{equation}
\end{theorem}
\proof By Example \ref{NahmpoleS3}, \ref{NahmpoleT3} and \ref{NahmpoleY3}, we know we have model Nahm pole solutions for $Y\times (0,+\infty)$ when Y is $S^3$, $T^3$ or any hyperbolic manifold. Denote the limit of the model solution as $\rho$ which is a flat $\slc$ connection. Here the model solution for hyperbolic manifold has limit in cylindrical end to a irreducible flat $\slc$ connection. 

Let $M_{\infty}=Y\times (0,+\infty)\cup M$, choose the flat connection $\rho$ and this will give a solution to the Kapustin-Witten equations over $M_{\infty}$. For $Y$ hyperbolic, we use the assumption $\pi_1$ injective in order to obtain a flat $\slc$ connection over $M_{\infty}$ with limit the irreducible $\slc$ connection over the cylindrical end coming from the hyperbolic metric. 

Applying Theorem \ref{gluingobstrucion} and Theorem \ref{glureducible}, we can glue these two solutions together and by Corollary \ref{perturbationequations}, we prove the statement for $h_1$, $h_2$. 
\qed

In addition, by gluing the model solutions $(A_0,\Phi_0)$ on Example \ref{NahmpoleS3}, \ref{NahmpoleT3} and \ref{NahmpoleY3} with themselves, we get the following corollary:

\begin{corollary}
For a 3-manifold $Y^3$ equals to $S^3$, $T^3$ or any hyperbolic 3-manifold, for T large enough, there exists a solutions $(A,\Phi)$ over $Y^3\times (-T,T)$ to the twisted Kapustin-Witten equations
$$KW(A,\Phi)+h=0.$$ 

Here $(A,\Phi)$ satisfies the Nahm pole boundary condition over $Y^3\times \{-T\}$ and $Y^3\times \{T\}$ and $h$ can be choosen to be support on any interior open set.
\end{corollary}

\end{section}

\section*{Appendix 1}
In this Appendix, we will give a brief introduction to the Fredholm theory of uniformly degenerate  elliptic operators that is developed in \cite{rafe1991elliptic} and \cite{mazzeo2013nahm}. We use the notation from \cite{mazzeo2013nahm} for most of the definitions in this Appendix. 

Let $M$ be a compact smooth 4-manifold with 3-manifold boundary $Y$ and choose coordinates $(\ox,y)$ near the boundary where $y\geq 0$ and $\ox=(x_1,x_2,x_3)\in \mathbb{R}^3.$ A differential operator $\MD_0$ is called uniformly degenerate if for $\alpha=(\alpha_1,\alpha_2,\alpha_3)$, in any coordinate chart near the boundary, it has the form
\begin{equation}
\MD_0=\sum_{j+|\alpha|\leq m}A_{j\alpha}(\ox,y)(y\partial_y)^j(y\partial_x)^{\alpha},
\end{equation}
where $(y\partial_x)^{\alpha}=(y\partial_{x_1})^{\alpha_1}(y\partial_{x_2})^{\alpha_2}
(y\partial_{x_3})^{\alpha_3}.$ 

We define the leading term of $\MD_0$ in this coordinate chart as
\begin{equation}
\MD_0^m:=\sum_{j+|\alpha|=m}A_{j\alpha}(\ox,y)(y\partial_y)^j(y\partial_x)^{\alpha}.
\label{leadingterm}
\end{equation}

The operator $\MD_0$ is called uniformly degenerate elliptic if $\MD_0$ is elliptic at the interior point and if in a neighborhood of the boundary and for (\ref{leadingterm}), we replace each $y\partial_{x_i}$ and $y\partial_y$ by variables $\sqrt{-1}k_i$ and $\sqrt{-1}k_4$ and it is invertible when $(k_1,\cdots,k_4)\neq 0.$

There is a model operator over $\mathbb{R}^4_+$, called the indicial operator
\begin{equation}
I(\MD_0)=\sum_{j\leq m}A_{j0}(\vec{x},0)\lambda^j. 
\end{equation}

The indicial root of $I(\MD_0)$ is the set of complex numbers $\lambda$ such that $s^{-\lambda}I(\MD_0)s^\lambda$ is not invertible.

In \cite{rafe1991elliptic}, Mazzeo works in the class of pseudodifferential operators on $M$ adapted to some particular type of singularity which includes the Nahm pole bounary condition. The class is called 0-pseudodifferential operators. Denote by $\Psi^{\star}_0(M)$ the elements which are described by the singularity structure of their Schwartz kernels. 

Given a pseudodifferential operator $A$, we denote the Schwartz kernel of $A$ as $\Mk_A(y,\ox,y,\ox')$ which is a distribution over $M^2:=M\times M$. We allow $\Mk_A$ to have the standard singularity of pseudodifferential operator along the diagonal $\{y=y',\ox=\ox'\}$ and we will require some special behavior over the boundary of $M^2$, which in coordinates is described as $\{y=0,y'=0\}$ and over the intersection of diagonal with the boundary, $\{y=0,y'=0,\ox=\ox'\}$. 

Let $M_0^2$ be a real blow-up of $M^2$ at the boundary of diagonal, which is constructed by replacing each point in $\{y=0,y'=0,\ox=\ox'\}$ with its inward-pointing normal sphere-bundle. We can describe it in polar coordinates:
$$R=(y^2+(y')^2+|\ox-\ox'|^2)^\frac{1}{2},\omega=(\omega_0,\omega_0^{'},\hat{\omega})=(\frac{y}{R},\frac{y'}{R},\frac{ \ox-\ox'}{R}).$$ Each point at $R=0$ is replaced by a quarter-sphere and $(R,\omega,x')$ can be regarded as a full set of coordinates. $M_0^2$ is a manifold with corners, we call the surface corresponding to $R=0$ the front face. The surfaces corresponding to $\omega_0=0$ and $\omega_0'=0$ are called its left and right faces. We have an obviously blow-down map $\pi:M_0^2\rightarrow M^2.$ We say $A\in \Psi^{\star}_0$ if $\Mk_A$ is the push forward of a distribution on $M_0^2$ by the blow-down map  $\pi$. 

Take a cut-off function $\chi$ over $M_0^2$ which is equals $1$ over a small neighborhood of the diagonal set $\{\omega_0=\omega_0',\hat{\omega}=0\}$ and 0 outside of a larger neighborhood. Then $\Mk_A=\Mk_A'+\Mk_A''$, where $\Mk_A'=\chi\Mk_A$ and $\Mk_A''=(1-\chi)\Mk_A$. Here $\Mk_A'$ supported away from the left and the right faces and has a pseudodifferential singularity of order $m$ along the lift diagonal area. If we factor $\Mk_A'=R^{-4}\hat{\Mk}_A'$, then $\hat{\Mk}_A$ extends smoothly of the front face of $M_0^2$ along the conormal diagonal singularity and $R^{-4}$ only depends on the manifold's dimension that corresponds to the determinant of the blow-down map $\pi$ and $\Mk_A''$ is smooth over the diagonal singularity.

Now we have the following definition of space $\Psi^{m,s,a,b}_0(M)$: 
\begin{definition}
For any real number $s,a,b$, we denote a psedudifferential operator $A\in \Psi^{m,s,a,b}_0(M)$ if its Schwartz kernel $\Mk_A$ has polyhomogeneous expansion with the terms $R^{-4+s}$ at the front face, $\omega_0^a$ at the left face and $\omega_0^b$ at the right face. 

We denote $A\in \Psi^{-\infty,a,b}(M)$ if its Schwartz kernels are smooth in the interior and polyhomogeneous at two hypersurfaces ($y=0$ and $y'=0$) of $M^2$.
\end{definition}

In this setting, the identity operator $Id\in \Psi_0^{0,0,\emptyset,\emptyset}$, has zero order over the diagonal and its Schwartz kernel $\delta(y-y')\delta(\ox-\ox')$ is supported over the diagonal which has a trivial expansion at the left and right faces. In polar coordinates, we have the following identification 
\begin{equation}
\delta(y-y')\delta(\ox-\ox')=R^{-4}\delta(\omega_0-\omega_0')\delta(\hat{\omega}),
\end{equation}
and this corresponds to zero in the second superscript.

Now, suppose $P$ is an $SU(2)$ bunlde over $M$ and let $(A,\Phi)\in\CP$ be a Nahm pole solution to the Kapustin-Witten equations. For simplification, let $\MD:=\DKW$. We denote $\MD_0=y\MD$. As pointed out in \cite{mazzeo2013nahm}, $\MD_0$ is a uniformly degenerate operator of order 1. Choose $p\geq 2$ and $q$ satisfying $\frac{1}{p}+\frac{1}{q}=1$.

In \cite{mazzeo2013nahm} Section 5.3, Mazzeo and Witten prove the following result:
\begin{theorem}{\rm{\cite{mazzeo2013nahm}}}
There exists operators $\MS\in \Psi_0^{-1,1,\bar{\lambda},b}(M)$, $R_1\in \Psi^{-\infty,\bar{\lambda},b}(M)$ and $R_2\in\Psi^{-\infty,b,\bar{\lambda}}$ such that 
\begin{equation}
\MD\circ\MS=Id-R_1,\;\;\MS\circ\MD=Id-R_2,
\end{equation}
where $\bar{\lambda}=1$ for the case this paper considered and $b\geq 1$. 
\end{theorem}

In \cite{rafe1991elliptic}, Mazzeo prove the following lemma about the distribution $\Psi_0^{m,s,a,b}.$
\begin{lemma}{\rm{\cite{rafe1991elliptic} Theorem 3.25, Remark above Proposition 3.28}}

$\bullet$ For any real number $\delta$ and $\delta'$, take $A\in \Psi^{-\infty,s,a,b}$ and let $A'=y^{\delta'}A {y'}^{-\delta}$ be its conjugation, then we have $A'\in \Psi^{-\infty,s+\delta-\delta',a-\delta',b+\delta}$.

$\bullet$ For $A\in \Psi^{-\infty,s,a,b}$, if $a>  -\frac{1}{p}$, $b>  -\frac{1}{q}$, $s\geq 0$, and $u\in L^p$, $v\in L^q$, we have $\vert\lan Au,v\ran _{L^2}\vert\leq\|u\|_{L^p}\|v\|_{L^q}$, which implies that $A$ is a bounded operator from $L^p$ to $L^p$.
\label{importantlemma}
\end{lemma}

We have the following proposition whose proof is slightly modified from \cite{rafe1991elliptic} Section 3 due to R. Mazzeo.

\begin{proposition}
For any real number $\lambda$ and any $p> 1$, the operator $\MD$,
\begin{equation}
\MD:\LamP H_0^{1,p}(M)\rightarrow \LamM L^p(M)
\end{equation} is a bounded linear operator. 
\end{proposition}
\proof As $\MD$ is a differential operator, the result follows immediately from the definition of $\MD$. 
\qed

For the operator $\MS$, we have the following proposition:
\begin{proposition}
For $\lambda\in (-1,1)$, the operator $\MS\in\Psi_0^{-1,1,1,b}(M)$:
\begin{equation}
\MS:\LamM L^p(M)\rightarrow \LamP H^{1,p}_0(M)
\end{equation}
is a bounded linear operator.
\label{bdstatement}
\end{proposition}
\proof 
Denote the Schwartz kernel of $\MS$ as $\Mk_S$.  By choosing a cut off function over the diagonal,  $\Mk_S=\Mk_S'+\Mk_S''$ where $\Mk_S'$ supported away from the left and the right faces and has a pseudodifferential singularity of order $m$ along the lift diagonal area and $\Mk_S''$ is smooth over the diagonal. Denote $S'$ ($S''$) to be the operator corresponds to the Schwartz kernel $\Mk_S'$($\Mk_S''$). 

We first prove that $S':\LamM L^p(M)\rightarrow \LamP H^{1,p}_0(M)$ is bounded. It is sufficient to prove that $S'y^{-1}$ is bounded operator from $\LamP L^p(M)$ to $\LamP H^{1,p}_0(M)$. We denote $A':=S'y^{-1}$ and now the $A'$ is dilation invariant. Choose a Whitney decomposition of $M$ into a union of boxes $B_i$ whose diameter in $x$ and $y$ directions is comparable to the distance to $\partial M$. For each $B_i$, we can choose an affine map $p_i$ which identifies a standard box $B$ with $B_i$. For $f\in \LamM L^p(M)$, denote by $f_i$ its restriction to $B_i$. Then $\|f\|_{\LamP H^{1,p}_0(M)}$ and $\sum_i y_i^{-\lambda-\frac{1}{p}}\|f_i\|_{H^{1,p}_0(B_i)}$ are comparable to each other where $y_i$ can be the $y$ coordinate of any points in $B_i$ and same for $\|f\|_{\LamM L^p}$ and $\sum_i y_i^{-\lambda-\frac{1}{p}}\|f_i\|_{L^p(B_i)}$. We denote $A_i'$ to be the restriction of $A'$ over $B_i$ then we have $p_i^{\star}(A'f)_i=A'_i(p^{\star}_if_i).$ By the approximate dilation invariance, we know $A_i'$ are a uniformly bounded family of psedodifferential operators. Then we have $$\|A'f\|_{\LamP H^{1,p}_0(M)}\leq C\sum_iy_i^{-\lambda-\frac{1}{p}}\|(A'f)_i\|_{H^{1,p}_0(B_i)}\leq C\sum_iy_i^{-\lambda-\frac{1}{p}}\|p_i^{\star}(A'f)_i\|_{H^{1,p}_0(B_i)}.$$ The classical $L^p$ theory about pseudodifferential operators of order $1$ in every box \cite{taylor1996pseudodifferential} gives $$\|p_i^{\star}(A'f)_i\|_{H^{1,p}_0(B_i)}=\|A'_i p_i^{\star}f_i\|_{L^p_1(B)}\leq C\|p_i^{\star}f_i\|_{L^p(B)}=C\|f_i\|_{L^p(B_i)}.$$

Summarizing the discussion above, we get $$\|A'f\|_{\LamP H_0^{1,p}(M)}\leq C\|f\|_{\LamP L^p(M)},$$ thus we get $\|S'f\|_{\LamM H_0^{1,p}(M)}\leq \|f\|_{\LamP L^p(M)}$.

Now, let's consider the operator $S''$, for some $b> 1$. For any integer $k$ and $k'$, we 
will show $S'':\LamM H_0^{k,p}(M)\rightarrow \LamP H_0^{k',p}(M)$ is a bounded operator. As $S''$ is an infinite smoothing operator over the diagonal, we only need to prove $S''$ is bounded from $\LamM L^p$ to $\LamP L^p$. Denote $A''={\LamP}S''y^{1-\frac{1}{p}-\lambda}$, then after the shifting, on the left faces, $A''$ will be polyhomogenous with leading order $b+\lambda+\frac{1}{p}-1$. In order to get bounds of the Schwartz kernel, we require that $b+\lambda+\frac{1}{p}-1> -\frac{1}{q}$. When $\lambda> -1$, this is automatically satisfied as $b> 1$. The leading order on the right faces will be $\bar{\lambda}-\lambda-\frac{1}{p}$, as $\lambda< 1$ and $\bar{\lambda}=1$, we automatically get $\bar{\lambda}-\lambda-\frac{1}{p}< -\frac{1}{p}$. By applying the second bullet of lemma $\ref{importantlemma}$, we get $A''$ is bounded from $L^p$ to $L^p$ which implies that $S''$ is bounded from $\LamM L^p$ to $\LamP L^p.$ 
\qed

For the operator $R_1$, $R_2$, we have the following proposition:
\begin{proposition}
For $\lambda\in(-1,1)$, $i=1,2$, and any $\lambda'\leq 1$, the operator $R_i$
\begin{equation}
R_i:y^{\lambda+\frac{1}{p}}H_0^{k,p}(M)\rightarrow y^{\lambda'+\frac{1}{p}}H_0^{k',p}(M)
\end{equation}
is a bounded for any $k,k'$. 
In addition, 
$$R_i:\LamP L^p\rightarrow \LamP L^p$$ is a compact operator.
\label{cptstate}
\end{proposition}
\proof We first prove the bounded statement. As $R_1\in \Psi^{-\infty,1,b}$, it is smooth over the diagonal, we only need to prove that $R_1$ is a bounded operator from $\LamP L^p$ to $y^{\lambda'+\frac{1}{p}}L^p$. Using the same trick as the previous proposition, we denote $R_1'=y^{\delta'}R_1y^{-\lambda-\frac{1}{p}}$. In order to get $C^0$ bound of the Schwartz kernel, now we require  $\delta'< \frac{1}{p}+1$ on the left face, which implies $\lambda'\leq 1$, same argument works for $R_2$. 

By Arzela-Ascoli theorem, we get that $R_i$ is compact operator.
\qed

\section*{Appendix 2}
Let $X$ be a manifold with boundary $Z$ and cylindrical end with a fixed limit $\slc$ flat connection, then for any connection pairs $(A_0,\Phi_0)$ satisfying the Nahm boundary condition over $Z$ and converges to $\slc$ flat connection over the cylindrical end in $L^p_1$ norm for some $p>2$, we will prove the closeness property of the operator $\DZ$. In this appendix, we assume $k\geq 0$ and $\lambda\geq -1$.  

We have the following lemma about bounded linear operators between Banach spaces:
\begin{lemma}{\rm{\cite{wehrheim2004uhlenbeck} Appendix E, Lemma E.3}}\label{1000}
Let $D:X\rightarrow Y$ be a bounded operator between Banach spaces.

(i) The following are equivalent:

$\bullet$ $D$ has a finite dimensional kernel and its image is closed.

$\bullet$ There exists a compact operator $K:X\rightarrow Z$ to another Banach space $Z$ and a constant $C$ such that 
\begin{equation}
\|u\|_{X}\leq C(\|Du\|_Y+\|Ku\|_Z)\;\forall u\in X.
\label{101}
\end{equation}

(ii) The following are equivalent:

$\bullet$ $D$ is injectie and its image is closed.

$\bullet$ There exists a constant C such that
\begin{equation}
\|u\|_X\leq C\|Du\|_Y\;\forall u\in X.
\label{102}
\end{equation}
\end{lemma}
In particular, if a bounded linear operator satisfies (\ref{101}) and it is injective, then it satisfies (\ref{102}).

Consider the operator $d_{A_0}$ associated with the following norms defined as:
\begin{equation}
    d_{A_0}:y^{\lambda+\frac{1}{p}+1} H^{2,p}_0(\Omega^0(\gpp))\rightarrow y^{\lambda+\frac{1}{p}} H^{1,p}_0(\Omega^1(\gpp)).
\end{equation}
By the definition of the norm, $d_{A_0}$ is a bounded linear operator. 

Let $\Omega^{odd}(\gpp)$ be the direct sum of odd differential forms and let $\Omega^{even}(\gpp)$ be the direct sum of even differential forms. Consider the following operator:
\begin{equation}
\begin{split}
    \MK:y^{\lambda+\frac{1}{p}+1} H^{2,p}_0(\Omega^{even}(\gpp))\rightarrow \LamP H^{1,p}_0(\Omega^{edd}(\gpp)).\\
\end{split}
\end{equation}
We denote $\MK_0=y\MK$ and we will study the semi Fredholm property of operators $\MK$ and $\MK_0$.

\begin{proposition}
$\MK_0$ is a uniformly degenerate elliptic operator and 0 is the only indicial root.
\end{proposition}
\proof The statement of uniformly degenerate elliptic operator is obvious. The indicial operator of $\MK_0$ is $I(\MK_0,\lambda)=A_{10}\lambda$ where $A_{10}$ is an invertible matrix. Thus $I(\MK_0,\lambda)$ is not invertible if and only if $\lambda =0.$
\qed

In \cite{rafe1991elliptic} Theorem 6.1, Mazzeo proves the following semi Fredholm theory of uniformly degenerate operator:
\begin{theorem}{\rm{\cite{rafe1991elliptic} Theorem 6.1}}
For any $\lambda>0$, there exist operators $G$ and $P$ such that 
$$G\MK_0=Id-P.$$
Here $G$ is a bounded operator $G:y^{\lambda+\frac{1}{p}}H_0^{1,p}\rightarrow y^{\lambda+\frac{1}{p}}H_0^{2,p}$ and $P$ is a compact operator. 
\end{theorem}

\begin{remark}
As there is only one indicial root, the $\bar{\lambda}$ in the original statement has to be 0. The bounded operator statement and compact operator statement can be proved in a similar way as Proposition \ref{bdstatement} and Proposition \ref{cptstate}.
\end{remark}

An immediately corollary of this theorem is that 
\begin{corollary}
When $\lambda>0$, $\MK_0$ has finite dimensional kernel and closed range. 
\end{corollary}
\proof
By the previous theorem, as $P$ is a compact operator and if $f\in \Ker \MK_0$, we have $Pf=f$. Therefore, the kernel of $P$ is finite dimensional. By Lemma \ref{1000} and the boundness property of $G$, we know $\MK_0$ has closed range.
\qed

We have the following proposition:
\begin{proposition}
For $\lambda>-1$, the $d_{A_0}:y^{\lambda+\frac{1}{p}+1} H^{2,p}_0(\Omega^0(\gpp))\rightarrow y^{\lambda+\frac{1}{p}} H^{1,p}_0(\Omega^1(\gpp))$ has finite dimensional kernel and closed range.
\label{dA0close}
\end{proposition}
\proof 
WLOG, we can assume $\Ker d_{A_0}$ is zero and prove the closed range statement. As $\Omega^0(\gpp)$ is a closed subset of $\Omega^{even}(\gpp)$, the restriction of $\MK_0$ over $\Omega^0(\gpp)$ which is $yd_{A_0}$ also has closed image. By Lemma $\ref{1000}$, we have the following inequality
\begin{equation}
\|u\|_{y^{\lambda+\frac{1}{p}+1} H^{2,p}_0(\Omega^0(\gpp))}\leq C\|yd_{A_0}u\|_{y^{\lambda+\frac{1}{p}+1} H^{1,p}_0(\Omega^1(\gpp))},
\end{equation}
which implies
\begin{equation}
\|u\|_{y^{\lambda+\frac{1}{p}+1} H^{2,p}_0(\Omega^0(\gpp))}\leq C\|d_{A_0}u\|_{y^{\lambda+\frac{1}{p}} H^{1,p}_0(\Omega^1(\gpp))}.
\end{equation}

Thus $d_{A_0}$ has closed range.
\qed
\begin{remark}
If $\Ker$ $d_{A_0}$=0, we have $\|\xi\|_{y^\lambda+\frac{1}{p}+1}\leq C\|d_{A_0}\xi\|_{y^\lambda+\frac{1}{p}}$
which is a gauge theory version of the $L^p$ Hardy inequality over $\mathbb{R}^4_+$ for compact supported functions $u$ and $s=p\lambda+p+1$:
\begin{equation}
(\int_{\RF}y^{p-s}|\partial_y u|^p)^{\frac{1}{p}}\geq \frac{n-1}{p}(\int_{\RF}y^{-s}|u|^p)^{\frac{1}{p}}.
\end{equation}
\end{remark}

Now we have the following proposition:
\begin{proposition}
The operator 
$$\DZ:y^{\lambda+1+\frac{1}{p}} H^{2,p}_0(\Omega^0(\gpp))\rightarrow y^{\lambda+\frac{1}{p}} H^{1,p}_0(\Omega^1(\gpp)\times\Omega^1(\gpp))$$
is a closed operator with finite dimensional kernel.
\end{proposition}
Recall the definition of $\DZ$ is $\DZ(\xi)=(d_{A_0}(\xi),[\Phi_0,\xi])$. Therefore, we obtain $\Ker\; \DZ\subset \Ker\; d_{A_0}$ and by Proposition \ref{dA0close}, we know $\Ker\;d^0_{(A_0,\Phi_0)}$ has finite dimension.

Without loss of generality, we assume $\Ker\;\DZ=0$. By Proposition \ref{dA0close}, there exists a constant such that $\|u\|_{y^{\lambda+\frac{1}{p}+1} H^{2,p}_0(\Omega^0(\gpp))}\leq C\|d_{A_0}u\|_{y^{\lambda+\frac{1}{p}} H^{1,p}_0(\Omega^1(\gpp))}.$ By adding a positive term on the right hand side of the inequality, we have 
\begin{equation*}
    \begin{split}
        &\|u\|_{y^{\lambda+\frac{1}{p}+1} H^{2,p}_0(\Omega^0(\gpp))}\\
        \leq& C(\|d_{A_0}u\|_{y^{\lambda+\frac{1}{p}} H^{1,p}_0(\Omega^1(\gpp))}+\|[\Phi_0,u]\|_{y^{\lambda+\frac{1}{p}} H^{1,p}_0(\Omega^1(\gpp))})\\=&C\|\DZ u\|_{\LamP H_0^{1,p}}. 
    \end{split}
\end{equation*}
Applying Lemma \ref{1000}, $\DZ$ is a closed operator.
\qed

\section*{Acknowledgements}
The author greatly thanks Ciprian Manolescu, Jianfeng Lin, Rafe Mazzeo and Thomas Walpuski for their kindness and helpful discussions. Part of this work was done when the author was visiting Stanford University and the author is grateful to Rafe Mazzeo for his hospitality.
\medskip

\bibliographystyle{unsrt}
\bibliography{references}
\end{document}